\documentclass[12pt]{article}
\usepackage{amsmath, amssymb, xcolor, hyperref}

\newlength{\abstractwidth}
\renewcommand{\thefootnote}{\fnsymbol{footnote}}
\renewcommand{\thanks}[1]{\footnote{#1}} 
\newcommand{\starttext}{ \setcounter{footnote}{0}
\renewcommand{\thefootnote}{\arabic{footnote}}}

\newcommand{\be}{\begin{equation}}
\newcommand{\bea}{\begin{eqnarray}}
\newcommand{\eea}{\end{eqnarray}} \newcommand{\ee}{\end{equation}}

\def\ba{\begin{eqnarray}}
\def\ea{\end{eqnarray}}

\def\cM{{\cal M}}

\def\o{\omega}
\def\Re{{\rm Re}}
\def\Im{{\rm Im}}
\def\tr{{\rm tr}}
\def\det{{\rm det}}

\def\log{\,{\rm log}\,}

\def\o{\omega}

\def\e{\varepsilon}

\def\o{\omega}

\def\na{\nabla}

\def\Z{{\bf Z}}

\def\R{{\bf R}}
\def\C{{\bf C}}

\def\p{\prod}

\def\cM{{\cal M}}

\def\na{{\nabla}}

\def\[{{\bf [}}
\def\]{{\bf ]}}

\def\p{\partial}

\begin{document}
\starttext \baselineskip=15pt \setcounter{footnote}{0}
\newtheorem{theorem}{Theorem}
\newtheorem{lemma}{Lemma}
\newtheorem{definition}{Definition}
\newtheorem{proposition}{Proposition}
\newtheorem{corollary}{Corollary}
\setcounter{tocdepth}{2}

\begin{center}
{\Large \bf  GEOMETRIC FLOWS FOR THE TYPE IIA STRING}
\footnote{Work supported in part by the National Science Foundation Grants DMS-1855947 and DMS-1809582.}
\bigskip\bigskip

\centerline{Teng Fei, Duong H. Phong, Sebastien Picard,
and Xiangwen Zhang}

\begin{abstract}

{\footnotesize A geometric flow on $6$-dimensional symplectic manifolds is introduced which is motivated by supersymmetric compactifications of the Type IIA string. The underlying structure turns out to be SU(3) holonomy, but with respect to the projected Levi-Civita connection of an almost-Hermitian structure. The short-time existence is established, and new identities for the Nijenhuis tensor are found which are crucial for Shi-type estimates. The integrable case can be completely solved, giving an alternative proof of Yau's theorem on Ricci-flat K\"ahler metrics. In the non-integrable case, models are worked out which suggest that the flow should lead to optimal almost-complex structures compatible with the given symplectic form.}

\end{abstract}

\end{center}


\section{Introduction}
\setcounter{equation}{0}

There has been a remarkable confluence in recent years between high energy physics, more specifically unified string theories, and special geometry. The earliest and particularly influential development was the 1985 recognition by Candelas, Horowitz, Strominger, and Witten \cite{CHSW} of Calabi-Yau manifolds as supersymmetric compactifications of the heterotic string. The importance of geometry in physical laws at their most fundamental level has of course been long recognized with electromagnetism, general relativity, and gauge theories. The key new feature here is the requirement that the $6$-dimensional internal manifold have a special geometric structure, in this case a complex structure together with a holomorphic section of the canonical bundle. This requirement is equivalent to the manifold having SU(3) holonomy, and can be traced back to supersymmetry. Since then, the class of Calabi-Yau solutions has been enlarged in many directions. On one hand, the Calabi-Yau condition can be extended to the Hull-Strominger system for conformally balanced metrics \cite{H, S, FY1, FY2, PPZ3, PPZ4, GF, FeiY, FHP}. On the other hand, it emerged from the mid 1990's that string theories, which are $10$-dimensional and of which there are five, can be unified themselves into another theory, namely M Theory, one of whose limits is $11$-dimensional supergravity \cite{HW, Town}. The compactification of an $11$-dimensional space-time to a more familiar $4$-dimensional space time results in a $7$-dimensional internal space, and the role of Calabi-Yau manifolds is assumed in this case by manifolds with ${\rm G}_2$
or Spin(7) holonomy \cite{BBS, Gu}.

\medskip
While the full string theories have been conjectured to merge ultimately into a single M Theory,  this requires highly non-trivial dualities,  and their low-energy approximations and geometric settings can be quite different. A common feature of their supersymmetric compactifications is a metric satisfying a curvature condition as well as a cohomological condition. In K\"ahler geometry, these are characteristic features of the notion of canonical metric, of which the Calabi-Yau condition is the prime example. Thus the general case can be viewed as a search for canonical metrics in non-K\"ahler geometry. The compactifications discussed above arise from three of the string theories, which are the Type I theory and the two heterotic string theories.
The other two string theories are the Type IIA and the Type IIB theories. There is an immense literature on their supersymmetric compactifications, but some attractive mathematical formulations can be found in Grana et al. \cite{Getal} and Tomasiello \cite{T}, and the study of the most basic examples was begun in Tseng and Yau \cite{TY1,TY2,TY3}. The geometric realm for the Type IIB equation is that of complex geometry, albeit non-K\"ahler, and we described a geometric flow approach to it in \cite{FPPZ}. The main goal of the present paper is to present a geometric flow approach to the Type IIA equation.

\medskip
The Type IIA equation is of particular interest, because of all string theory compactifications, its geometric setting is unique in being that of symplectic geometry instead of complex geometry. More specifically, let $M$ be a compact $6$-dimensional manifold, equipped with a symplectic form $\o$, that is, a closed non-degenerate $2$-form. Recall that on any oriented $6$-manifold, Hitchin \cite{H} had shown how to associate to a non-degenerate $3$-form $\varphi$ an almost-complex structure $J_\varphi$. In the Type IIA equation, a symplectic form $\o$ is given, so it makes sense to consider the condition of primitiveness for $\varphi$ with respect to $\o$. Explicitly, this is the condition $\Lambda\varphi=0$, where
\bea
\Lambda:A^{k}(M)\to A^{k-2}(M)
\eea
is the usual Hodge operator of contracting with $\o$. As shown in \S 4.1 below, the symplectic form $\o$ is invariant with respect to the almost-complex structure $J_\varphi$ when $\varphi$ is primitive. We obtain then a Hermitian form
\bea
g_\varphi(X,Y)=\o(X,J_\varphi Y).
\eea
which becomes a Hermitian metric under the open condition that it be strictly positive. Thus we obtain an almost-K\"ahler $3$-fold $(M,\o,J_\varphi,g_\varphi)$ with the additional requirement that $\varphi$ be primitive and $J_\varphi$ arise from $\varphi$ by the above construction. When $\varphi$ is also closed, we shall refer to such a structure as a ``Type IIA geometry".

\smallskip
Let $\rho_A$ be now the Poincar\'e dual of a given finite linear combination of special Lagrangians calibrated by $\varphi$. Then the Type IIA equation is the following system of equations for a real-valued primitive $3$-form $\varphi$
\bea
\label{eq:IIA}
d\Lambda d(|\varphi|^2\star\!\varphi)=\rho_A,
\qquad d\varphi=0,
\qquad g_\varphi>0.
\eea
Here $\star$ is the Hodge star operator and $|\varphi|$ the norm of $\varphi$ with respect to the metric $g_\varphi$.


\medskip
As in the case of the other string theories, the Type IIA equations as written in (\ref{eq:IIA}) involve, besides the open condition $g_\varphi>0$, a curvature-type equation and the cohomological constraint $d\varphi=0$. In order to enforce this cohomological constraint without invoking any particular Ansatz, we introduce the following geometric flow of $3$-forms
$\varphi$,
\bea
\label{fl:IIA}
\p_t\varphi=d\Lambda d(|\varphi|^2\star\!\varphi)-\rho_A,
\eea
for any closed and primitive initial data $\varphi_0$ with $g_{\varphi_0}>0$. Since the right hand side is closed, the flow preserves the closedness condition. It can also be verified to preserve the primitiveness property of $\varphi$.
Thus it is a flow of Type IIA geometries, whose stationary points would give solutions of the Type IIA equation without recourse to any Ansatz. We shall refer to (\ref{fl:IIA}) as the Type IIA flow. The idea of preserving the closedness of a form by a flow was introduced by Bryant \cite{B} in the Laplacian flow for $G_2$ structures. More recently, it was applied in \cite{PPZ1, PPZ2, FPPZ, FP, Phong} to the construction of geometric flows which preserve the conformally balanced condition in the Hull-Strominger system and the Type IIB equation. The geometric flow approach was particularly appropriate there, as it allowed to bypass the absence of a $\p\bar\p$-lemma in non-K\"ahler geometry.

\medskip
The main goal of this paper is to start an in-depth study of the Type IIA flow. Except for the original formulation, we shall restrict to the most basic source-free case $\rho_A=0$.
Despite its very simple formulation (\ref{fl:IIA}), the flow turns out to be highly non-trivial and to present many new difficulties specific to symplectic geometry:

\medskip
$\bullet$ The first is that the Type IIA equation is not elliptic. This difficulty was well-recognized in the works of Tseng and Yau \cite{TY3} and Tseng and Wang \cite{TW}, and led them to consider instead some 4th-order equations. However, 4th-order equations are complicated, and the closedness constraint on the $3$-form $\varphi$ would have to be imposed separately. Thus it appears still preferable to confront the specific difficulties of the Type IIA flow. They originate in any case from the geometric assumption of a given symplectic structure, which is fundamental in symplectic geometry.

\smallskip
$\bullet$ The second may be appreciated by comparing the flow of almost-complex structures $J_\varphi$ in the Type IIA flow with the gradient flow of the Blair-Ianus functional
on a symplectic manifold. The Blair-Ianus functional is the $L^2$ norm of the Nijenhuis tensor \cite{BI}. Its gradient flow was called the anti-complexified Ricci flow by L\^e and Wang \cite{LW}, who also established its short-time existence. However, this flow has proved to be difficult to use, because neither the corresponding  Nijenhuis tensor nor curvature evolves there by parabolic equations. For the Type IIA flow to be viable, it has to overcome such difficulties.

\smallskip
$\bullet$ The third difficulty is more technical, but still serious. The Type IIA flow of $3$-forms $\varphi$ induces a flow of metrics $g_\varphi$ which will be one of the main tasks of this paper to determine explicitly. The simplest case is when the initial almost-complex structure is integrable. It turns out that the Type IIA flow preserves the integrability condition, and becomes equivalent to the dual Anomaly flow introduced in \cite{FeiPicard}. Using the techniques there as well as in \cite{PPZ5, PZ}, it gives a new proof of Yau's \cite{Y} celebrated theorem on the existence of K\"ahler Ricci-flat metrics. Thus the difficult case is the case of non-integrable almost-complex structures. There we shall see that the flow of metrics in the Type IIA flow is conformally equivalent to a perturbation of the Ricci flow by first-order terms involving the Nijenhuis tensor. In this respect, it is analogous to Bryant's $G_2$ Laplacian flow, which was shown relatively recently by Lotay and Wei \cite{LoW} to be a perturbation of the Ricci flow by first-order terms involving the torsion tensor. However, the long-time behavior of the $G_2$ flow remains at this moment a subject of extensive research \cite{Lo, BV}.

\medskip

Despite these difficulties, we shall find that the Type IIA flow is particularly rich, reflecting its unique position at the crossroads of symplectic geometry, complex geometry, and unified string theories. This rich structure will be much in evidence in the results described below. Furthermore, examples suggest that even when the flow develops singularities, it may be possible in some cases to continue the flow of the Nijenhuis tensor. Thus, besides the original motivation from string theory and interest in its stationary points, the flow should also be useful in finding optimal almost-complex structures.

\section{Main results}
\setcounter{equation}{0}

We describe now our main results. Throughout this section, $M$ is a compact $6$-dimensional manifold equipped with a symplectic form $\o$. Given a primitive $3$-form $\varphi$, we denote by $J_\varphi$ the almost-complex structure defined by Hitchin \cite{H}, and by $g_\varphi$ the corresponding Hermitian form, which we assume is a metric.

\subsection{A Laplacian flow formulation}

Our first result is an alternative formulation of the Type IIA flow:

\begin{theorem}
\label{th:Laplacian}
The Type IIA flow defined in (\ref{fl:IIA}) can be rewritten as the following flow
\bea
\p_t\varphi=-dd^\dagger (|\varphi|^2\varphi)+2\,d(|\varphi|^2N^\dagger\cdot\varphi)-\rho_A
\eea
where $d^\dagger$ is the adjoint of the operator $d$ with respect to the metric $g_\varphi$, and $N^\dagger$ is the operator from $\Lambda^3(M)$ to $\Lambda^2(M)$ defined by
\bea \label{Ndagger}
(N^\dagger\cdot\varphi)_{kj}=N^\mu{}_j{}^\lambda\varphi_{\mu k\lambda}-N^\mu{}_k{}^\lambda\varphi_{\mu j\lambda}.
\eea
Here $N^m{}_{\gamma\beta}$ is the Nijenhuis tensor of $J_\varphi$, and indices are raised using the metric $g_\varphi$.
\end{theorem}

We note that, up to the factor of $|\varphi|^2$, the first terms on the right-hand side of the Type IIA  flow are the same as in the standard heat equation. Up to sign, they are also reminiscent of Bryant's $G_2$ Laplacian flow. However, the terms involving the Nijenhuis tensor are also of leading order and account for a wide range of different phenomena.

\smallskip
Henceforth, we assume that the source $\rho_A$ is $0$, unless stated explicitly otherwise.

\subsection{The short-time existence of the Type IIA flow}

When a flow is not strictly parabolic, even its short-time existence can be a difficult question. Two powerful tools developed over the years for this issue have been the reparametrization method of DeTurck \cite{DeT} and the Hamilton-Nash-Moser theorem \cite{Ha1}, a combination of which has been applied successfully to many important flows, such as the Ricci flow \cite{Ha1}, the $G_2$ Laplacian flow \cite{B,BX}, and the anti-complexified Ricci flow \cite{LW}. The fundamental new difficulty in the Type IIA flow is that there is a given symplectic form $\o$. It is not hard to see that reparametrizations by symplectomorphisms do not improve the parabolicity of the flow, while more general reparametrizations lead to a {\it coupled flow} of both metrics and symplectic forms. Thus a first major task in this paper is to establish the following theorem:

\begin{theorem}
\label{th:short}
Let $(M,\o)$ be a compact $6$-dimensional symplectic manifold. Then for any $\varphi_0$ which is a smooth positive, primitive, and closed $3$-form, the source-free Type IIA flow (\ref{fl:IIA}) with initial value $\varphi_0$ admits a unique and smooth solution on some time interval $[0,T)$ with $T>0$. Furthermore, $\varphi$ continues to be positive, primitive and closed at all times.
\end{theorem}

While the theorem deals only with the short-time existence of the flow, the proof requires a rather deep probe of the structure of the flow and several new elements which are also useful elsewhere:

\smallskip

The first element is the behavior of the coupled flows mentioned above. It turns out that these coupled flows admit natural parabolic regularizations, which reduce to the desired flow for primitiveness initial data. While primitiveness was a requirement in the solution $\varphi$ of the Type IIA equation, it may not have been anticipated that it would play such a central role for the very existence of the flow.

The second feature permeates the rest of the paper, and is the underlying Type IIA geometry. In the present context, it allows us to recapture the flow of the forms $\varphi$ from the flow of the metrics $\tilde g_\varphi=|\varphi|^2g_\varphi$, even though, pointwise, there is an ambiguity in determining $\varphi$ from $\tilde g_\varphi$. Since we shall have to analyze in detail the flow of metrics in order to obtain Shi-type estimates and long-time existence criteria, it is simplest to deduce the uniqueness part of Theorem \ref{th:short} from the corresponding uniqueness theorems for the flow of metrics. The flow of $\check g_\varphi=|\varphi|^{-2}g_\varphi$ turns out to be a perturbation of the Ricci flow by terms of first order. In general, integrability operators are not stable under first-order terms perturbations. However, the underlying Type IIA geometry is what allows us to modify the Bianchi operator used by Hamilton \cite{Ha1} for the Ricci flow into an integrability operator for the flow of $\check g_\varphi$. From there we can establish the uniqueness for $\check g_\varphi$, and from there the uniqueness in Theorem \ref{th:short}.

\subsection{Type IIA geometry}

We have stressed that the underlying structure for the Type IIA flow is Type IIA geometry, as defined in the Introduction, and which is more special than just a symplectic structure on a $6$-manifold. The holonomy, curvature, and Nijenhuis tensor in Type IIA geometry have rather special properties, which play an important role in every aspect of the Type IIA flow.
We list here some properties which are applied repeatedly in the paper and are the easiest to describe, but we expect others to emerge and prove their worth in time:

\begin{theorem}
\label{th:holonomy}
Let $(M,\o,\varphi)$ be a Type IIA geometry, and $g_\varphi$ the corresponding metric.  Set
\bea
\label{tildeg0}
\tilde g_\varphi=|\varphi|^2g_\varphi.
\eea
Let $\frak D$ and $\tilde{\frak D}$ be the projected Levi-Civita connections of $g_\varphi$ and $\tilde g_\varphi$ respectively, $\Omega=\varphi+iJ_\varphi\varphi$, and $|\Omega|_{\tilde g_\varphi}$ the norm of $\Omega$ with respect to $\tilde g_\varphi$. Then

{\rm (a)} $\tilde{\frak D}({\Omega\over |\Omega|_{\tilde g_\varphi}})=0$. Thus $(M,\tilde g_\varphi)$ has holonomy in $SU(3)$, but with respect to the connection $\tilde{\frak D}$.

{\rm (b)} ${\frak D}^{0,1}\Omega=0$, so $\Omega$ is formally holomorphic, even when $J_\varphi$ is not integrable.

{\rm (c)} The Nijenhuis tensor has only $6$ independent components.

\end{theorem}

\subsection{The flow of metrics in the Type IIA flow}

Next, we can describe the flow of the $3$-forms $\varphi$ and metrics $g_\varphi$ in terms of curvature:

\begin{theorem}
\label{th:g-varphi}
Let $\varphi$ be a smooth positive, primitive and closed $3$-form evolving by the source-free Type IIA flow. Set $u=\log |\varphi|^2$.

{\rm (a)} The flow of $\varphi$ is given by
\bea
\p_t\varphi_{iab}=e^{2u}\!\!\!\!\!\sum_{\mathrm{cyc~} i,a,b}\!\!\!\!\left(\!\varphi_{sab}(\tilde{\frak D}_i+u_i)u^s+2\varphi_{sta}(N^{st}{}_p\tilde{\mathfrak T}^p{}_{ib}-\frac{u^s}{2}\tilde{\mathfrak T}^t{}_{ib}+(\tilde{\frak D}_i+u_i)N^{st}{}_b)\!\right)\!,\label{ode}
\eea
where $\tilde{\frak T}$ is the torsion tensor of the connection $\tilde{\frak D}$.

{\rm (b)} The flow of $\tilde g_\varphi$ is given by
\bea
(\p_t\tilde g_\varphi)_{ij}&=&e^{2u}\bigg[-2\tilde R_{ij}-2(\tilde\nabla^2u)_{ij}+4u^s(N_{isj}+N_{jsi})+u_iu_j-u_{Ji}u_{Jj}-4(N^2_-)_{ij}\nonumber\\
&&+(|du|^2_{\tilde g_\varphi}+|N|^2_{\tilde g})(\tilde g_\varphi)_{ij}\bigg],\label{tildeg}
\eea
where $\tilde{\nabla}$ is the Levi-Civita connection of $\tilde g_{\varphi}$.
\end{theorem}

Note that $u$ is determined by $\tilde g_\varphi$, so the right hand side of (\ref{tildeg}) involves only tensors determined by $\tilde g$, and the equation is a self-contained flow for $\tilde g_\varphi$.
Furthermore, the equation for $\varphi$ can be viewed as a linear ODE of $\varphi$ whose coefficients are tensors determined by $\tilde g_\varphi$. Thus $\varphi$ is completely determined once $\tilde g_\varphi$ is determined. As we noted before, this is to be contrasted with the problem of having to resolve an ambiguity if we just try to recapture $\varphi$ from $g_\varphi$ at each fixed time. It may be worth observing that
the ambiguity in recapturing $\varphi$ from $g_\varphi$ pointwise in time is reminiscent of the ambiguity in defining the angle in the special Lagrangian equation. It would be interesting to investigate if the angle in the special Lagrangian equation can be recaptured by a mechanism similar to the above Theorem \ref{th:g-varphi} for the Type IIA flow.

\subsection{An integrability operator for the flow of metrics}

The flow (\ref{tildeg}) is reminiscent of the Ricci flow, except for the term $(\tilde\na^2 u)_{ij}$ which can normally be eliminated by a reparametrization. But as we noted in the above discussion of Theorem \ref{th:short}, a reparametrization would create other difficulties since it would change the given symplectic form $\o$. To bypass this difficulty, we make instead a conformal change
\bea
\label{checkg}
\check g_{ij}=|\varphi|^{-2}(g_\varphi)_{ij}
 \eea
and find that $\check g_{ij}$ evolves by
\bea\label{flcheckg}
\p_t\check g_{ij}&=&e^{\frac{3}{2}u}\bigg[-2\check R_{ij}+\frac{3}{2}u_iu_j-u_{Ji}u_{Jj}+4u^k(N_{ikj}+N_{jki})-4(N^2_-)_{ij}\nonumber\\
&&+\frac{1}{2}\left(|du|^2_{\check g}+|N|^2_{\check g}\right)\check g_{ij}\bigg],\label{evolucheck}
\eea
For the purpose of completing Theorem \ref{th:short}, we are particularly interested in the uniqueness of this flow. For the classical Ricci flow, both the short-time existence and uniqueness were established by Hamilton \cite{Ha1} using a version of the Nash-Moser theorem. This version, often referred to as the Hamilton-Nash-Moser theorem, requires an integrability condition, which was provided by the Bianchi identity in the case of the Ricci flow. More precisely, let $L_0$ be the operator defined by
\bea
L_0: \ {\rm Sym}^2(TM)\ni S_{ij}\ \to\ (L_0(S))_j
=2\check g^{ik}\check\na_k S_{ik}-\check g^{ik}\check \na_j S_{ik}\in A^1(M)
\eea
Then $L_0(-2\check R_{ij})=0$, which is the desired integrability condition.
Now the flow (\ref{flcheckg}) differs from the Ricci flow by first-order terms, so the Bianchi identity is no longer applicable as an integrability condition. In general, it is by no means clear whether a given first-order perturbation would still allow an integrability condition. So it is again a manifestation of the deep structure of Type IIA geometry that this can be done in this case:

\begin{theorem}
\label{th:integrability}
Let $S_{ij}$ be the symmetric $2$-tensor defined by writing the flow (\ref{flcheckg}) as
\bea
\label{P}
\p_t\check g_{ij}=e^{{3\over 2}u}(-2\check R_{ij}+S_{ij}).
\eea

{\rm (a)} Define the operator $Z$ by
\bea
Z:\ {\rm Sym}^2(TM)\ni S_{ij}\ \to\ (Z(S))_j
=-{4\over 3}u_j\check g^{ik}S_{ik}+2u^sS_{js}-4N^{st}{}_jS_{st}
\in A^1(M).
\eea
If along the flow, the metric $\check{g}$ arises from a Type IIA geometry, then the following operator $L$ is an integrability operator for the flow (\ref{flcheckg}) in the sense that
\bea
\label{integ}
{\rm Sym}^2(TM)\ni \check g_{ij}\ \to\
L(\check g_{ij}):=(L_0+Z)(e^{{3\over 2}u}(-2 \check R_{ij}+S_{ij}))
\eea
is of order $1$ in $\check g_{ij}$.

\noindent
{\rm (b)} As a consequence, the flow exists and is unique on some interval $[0,T)$ with $T>0$.
\end{theorem}

We stress that this theorem is used only to establish the uniqueness of the Type IIA flow, but not its existence. The reason is that, starting from the flow (\ref{flcheckg}) at an initial Type IIA geometry, it is not a priori known whether the flow will remain a Type IIA geometry. Without this information, it is not known whether the above integrability condition (\ref{integ}) holds.

\subsection{Formation of singularities and Shi-type estimates}

In general, the fact that a flow may admit short-time existence does not imply that important geometric quantities will evolve by parabolic equations. For example, that is the case for the Ricci flow \cite{Ha1}, but it is not the case for the anti-complexified Ricci flow \cite{LW}. The Type IIA flow shares common features with both the Ricci flow and the anti-complexified Ricci flow, and one may wonder which way it will behave when it comes to the evolution of the curvature and of the Nijenhuis tensor. It is a very attractive feature of the Type IIA flow that, in this particular respect, it is closer to the Ricci flow. Thus we find

\begin{theorem}
\label{th:Nijenhuis}
Consider the source-free Type IIA flow with a smooth, positive, closed, and primitive initial value $\varphi_0$.

{\rm (a)} Then the Nijenhuis tensor evolves by
\bea \label{evol-norm-N}
(\p_t - e^u\Delta) |N|^2 &=& e^u \bigg[ - 2 |\nabla N|^2 + (\nabla^2
  u) * N^2 + Rm * N^2 \nonumber\\
  && + N * \nabla N * (N+\nabla u) + N^4 + N^3 * \nabla u + N^2 *
  (\nabla u)^2\bigg]
\eea
  {\rm (b)} The Riemann curvature tensor evolves by
\bea \label{evol-Rm-3}
 (\p_t - e^u \Delta) |Rm|^2 &=& e^u \bigg[ - 2 |\nabla Rm|^2 + (\nabla
  Rm + \nabla^3 u + \nabla^2 N) *
  \mathcal{O}(Rm, \nabla u,N) \nonumber\\
  &&+ (\nabla N * \nabla N + \nabla^2 u * \nabla^2 u +1) * \mathcal{O}(Rm, \nabla u,N) \bigg].
\eea
Here, $*$ denote the bilinear pairings (not to be confused with the Hodge star operator) and $\mathcal{O}(\nabla u, Rm, N)$ indicates terms which
only depend on $\nabla u, Rm$ and $N$.

\end{theorem}

Other geometric quantities satisfy similar flows, which are written in greater detail in \S 8. Using these flows, we can establish the following Shi-type estimates and criteria for extending the flow:

\begin{theorem}
\label{th:Shi}
Assume that we have a solution of the source-free Type {\rm IIA} flow on some interval $[0,T)$, and that the bound
\bea
|u|+|Rm|\leq A
\eea
holds for some finite constant $A$. Here $Rm$ denotes the Riemann curvature tensor of the metric $g_\varphi$. Then for any multi-index $\alpha$, we have
\bea
|\na^\alpha\varphi|\leq C(A,\alpha,T,\varphi(0))
\eea
for some constant $C(A,\alpha,T,\varphi(0))$. In particular, the Type {\rm IIA} flow can be continued to an interval $[0,T+\e)$ for some $\e>0$.
\end{theorem}

It may be worth noting that, in this estimate, the estimates for the gradient $|\na u|$ are rather special, and make essential use of the underlying Type IIA geometry.

\subsection{The stationary points in the case of no source}

In the case $\rho_A=0$, the stationary points of the flow can be identified, once we have developed Type IIA geometry:

\begin{theorem}
\label{th:stationary}
A primitive and closed $3$-form $\varphi$ is a stationary point of the flow if and only if the corresponding almost-complex structure $J_\varphi$ is integrable and the norm $|\varphi|$ is constant. Thus $(M,J_\varphi,\o)$ is then a K\"ahler manifold, and the metric $g_\varphi$ is K\"ahler and Ricci-flat.
\end{theorem}

\subsection{The integrable case and the Calabi conjecture}

The simplest case is that of integrable almost-complex structures, and a complete description of the behavior of the Type IIA flow in this case is provided by the following theorem:

\begin{theorem}
\label{th:integrable}
Assume that the initial value $\varphi_0$ of the source-free Type {\rm IIA} flow is a positive, closed and primitive $3$-form, and that the corresponding almost-complex structure $J_{\varphi_0}=:J_0$ is integrable. Then the source-free Type {\rm IIA} flow exists for all time, the almost-complex structure $J_\varphi$ corresponding to $\varphi$ remains integrable along the flow, and the flow converges in $C^\infty$ to a $3$-form corresponding to a K\"ahler Ricci-flat metric.
\end{theorem}

In fact, the corresponding flow of metrics $g_\varphi$ turns out to reduce by diffeomorphisms to the dual Anomaly flow introduced in \cite{FeiPicard} which applies in all dimensions and gives another proof of the Calabi conjecture. This reduction of the Type IIA flow to the dual Anomaly flow, which was itself motivated by duality considerations for the Type IIB flow, can be viewed as a manifestation of the duality between the Type IIA and the Type IIB theories.
We also observe that there are by now several known proofs of the Calabi conjecture, giving each a different sequence converging to the K\"ahler Ricci-flat metric:
besides Yau's original proof, there are for example the proof by the K\"ahler-Ricci flow \cite{C}, by the Anomaly flow \cite{PPZ5}, by the dual Anomaly flow \cite{FeiPicard}, by the inverse Monge-Amp\`ere flow \cite{CK, CHT}, and by more general parabolic Monge-Amp\`ere flows \cite{PZ}. Nevertheless, new proofs from an independent geometric set-up remain of considerable interest, as they can detect different types of obstructions. Such a scenario is nicely illustrated in \cite{CHT}.

\subsection{Examples}

As we can see from the induced flow of metrics, the Type IIA flow is complicated. However, besides the integrable case which was solved above, there are non-integrable, geometrically interesting cases that can also be worked out completely and which exhibit varied and interesting behaviors. They suggest the possibility of a general phenomenon, namely that in all cases, the Type IIA flow leads to an optimal almost-complex structure with respect to the given symplectic form.
A first example is provided by the torus:

\begin{theorem}
\label{th:Model2}
Consider the source-free Type IIA flow on the torus ${\bf R}^6/{\bf Z}^6$ with the symplectic form $\o$ as described in \S \ref{ex:torus}. Consider the Type IIA flow with non-integrable initial data of the form in (\ref{ex:torus}). Then the flow exists for all time, and $\varphi$ converges as $t\to\infty$ to a positive harmonic form.
\end{theorem}

A rich class of models is the special generalized Calabi-Yau (or SGCY) manifolds introduced by de Bartholomeis \cite{dB, dBT}. These are manifolds of Type IIA geometry with the additional property that $|\varphi|$ is constant. They are also sometimes referred to as symplectic half-flat structures. A large subclass is given by nilmanifolds, which are quotients of a nilpotent Lie group by a co-compact lattice, and the other subclass is the solvmanifolds, which are quotients of solvable Lie groups by a co-compact lattice. Details on the models which we consider can be found in \S \ref{ex:half-flat}, and we shall just state here the main conclusions.

\begin{theorem}
\label{th:Model3}
Consider the source-free Type {\rm IIA} flow on the nilmanifold and the solvmanifold described in \S \ref{ex:half-flat},
with the initial data described there.

{\rm (a)} In the case of the nilmanifold, with initial data corresponding to (\ref{data:nil}), the flow exists for all time, and the Nijenhuis tensor tends to $0$ as $t\to\infty$.

{\rm (b)} In the case of the  solvmanifold, with initial data corresponding to (\ref{solvans}),
the flow develops a singularity at a finite time $T$. However, the limit of $J_\varphi$ as $t\to T$ exists, and is a harmonic almost-complex structure.
\end{theorem}

\section{The Type IIA flow as a Laplacian flow}
\setcounter{equation}{0}
\label{s:formal}

We start now the proof of the results described in the previous section. For Theorem \ref{th:Laplacian},
we do not need any detailed information on the metric $g_\varphi$. Rather, we only need the main properties of the corresponding Hermitian connections, and how they differ from the Levi-Civita connection. These have been worked out by Gauduchon \cite{G}, and we begin with a brief review of the results from \cite{G} that we require.

\subsection{Gauduchon's formulas for connections on an almost-complex manifold}

We revert momentarily to the general set-up of a smooth manifold $M$ equipped with a Riemannian metric $g$, a compatible almost-complex structure $J$ (not necessarily integrable), and the associated symplectic form $\o$. This means that $g$ is a positive definite section of the bundle of quadratic forms on $TM$, $\o$ is a $2$-form on $M$, $J$ is a section of the bundle of endomorphisms of $TM$ satisfying $J^2=-\textrm{Id}$,
and $g(X,Y)=\o(X,JY)$ for any two vector fields $X,Y$ on $M$. In local coordinates $x^j$, with $X$ and $Y$ given by their components $X^i$ and $Y^j$,
we shall write
\bea
g(X,Y)=g_{ij}X^iY^j,
\qquad
\o(X,Y)=\o_{ij}X^iY^j,
\quad
(JX)^k=J^k{}_jX^j.
\eea
In particular $g_{ij}$ and $\o_{ij}$ are respectively symmetric and anti-symmetric in $i$ and $j$, and the fact that $J$ is a compatible almost-complex structure translates into
\bea
J^k{}_qJ^q{}_j=-\delta^k{}_j,
\qquad
g_{ij}=\o_{iq}J^q{}_j=\o_{jq}J^q{}_i.
\eea
Note that $\o$ is invariant under the action of $J$, i.e. $\o_{iq}J^q{}_jJ^i{}_\ell=\o_{\ell j}$.

Clearly the structure defined by the triple $(g,\o,J)$ is determined by any two of its components. Here we do not assume that $\o$ is closed. For our purposes, we can assume that $d\omega$ has no $(3,0)+(0,3)$-components. Obviously this condition is satisfied when $\o$ is closed, or conformally closed, or when $J$ is integrable. When such a condition is satisfied, $d^c\o=J^{-1}d\o$ has no $(3,0)+(0,3)$-components either.


Associated to this setup are several important tensors:

\medskip
$\bullet$ The first is the Nijenhuis tensor $N$, defined as
\begin{equation}
    N(X,Y)=\frac{1}{4}([JX,JY]-J[JX,Y]-J[X,JY]-[X,Y]).\label{nijenhuis}
\end{equation}
By construction $N$ is a skew-symmetric $2$-form valued in $TM$. Using the metric, we can lower the superscript to the first slot by $N_{X,Y,Z}=g(X,N(Y,Z))$.

\medskip
$\bullet$ The second is the $3$-form $d^c\o$, where by our convention for differential forms, the form $\o$ is defined from its coefficients $\o_{ij}$ by
\bea
\o={1\over 2}\o_{ij}dx^i\wedge dx^j.\nonumber
\eea
The form $d^c\o$ can also be written as $d^c\o=J^{-1}dJ\o=-Jd\o$ since $d^c=J^{-1}dJ$ and $\o$ is invariant under the action of $J$. In components, we have
\bea
(d\o)={1\over 3!} (\p_i\o_{jk}+\p_j\o_{ki}+\p_k\o_{ij})dx^i\wedge dx^j\wedge dx^k\nonumber
\eea
and hence
\bea
&&
(d\o)_{ijk}=\p_i\o_{jk}+\p_j\o_{ki}+\p_k\o_{ij},
\nonumber\\
&&
(d^c\o)_{abc}=-J^k{}_cJ^j{}_bJ^i{}_a(d\o)_{ijk}
=
-J^k{}_cJ^j{}_bJ^i{}_a(\p_i\o_{jk}+\p_j\o_{ki}+\p_k\o_{ij}).\nonumber
\eea

As shown by Gauduchon \cite{G}, the construction of Hermitian connections associated to the structure $(g,\o,J)$ is clearer if we also view $d^c\o$ as a $2$-form valued in $TM$, i.e. a section of $TM\otimes \Lambda^2T^*M$, in analogy with the Nijenhuis tensor. This can be achieved by raising one index in $d^c\o$, using the metric $g_{ij}$. Unless indicated otherwise, the $TM$-valued $2$-form corresponding to $d^c\o$ is obtained by raising the first index, i.e.,
\bea
(d^c\o)^m{}_{jk}=g^{mi}(d^c\o)_{ijk}.\nonumber
\eea

$\bullet$ The third and the fourth tensors of interest are obtained by decomposing $d^c\o$, viewed as a $TM$-valued $2$-form, into components $U$ and $V$ which are respectively odd and even under the following involution ${\cal M}$ acting on the space of $TM$-valued $2$-forms,
\bea
\label{M0}
({\cal M}\Psi)(X,Y)=\Psi(JX,JY),
\qquad \Psi\in A^2(TM)
\eea
where we have denoted the space of $TM$-valued $2$-forms by $A^2(TM)$. We can then define the $TM$-valued $2$-forms $U$ and $V$ by
\bea\label{UV}
U={1\over 4}(d^c\o+{\cal M}(d^c\o)),
\qquad
V={1\over 4}(d^c\o-{\cal M}(d^c\o)).
\eea
In components,
\bea
\label{M1}
({\cal M}\Psi)^m{}_{bc}=\Psi^m{}_{jk}J^k{}_cJ^j{}_b
\eea
and
\bea
\label{UV0}
U^m{}_{bc}={1\over 4}((d^c\o)^m{}_{bc}+(d^c\o)^m{}_{jk}J^k{}_cJ^j{}_b),
\quad
V^m{}_{bc}={1\over 4}((d^c\o)^m{}_{bc}-(d^c\o)^m{}_{jk}J^k{}_cJ^j{}_b).
\eea
Note that $U^m{}_{bc}$ and $V^m{}_{bc}$ are still anti-symmetric in $b$ and $c$, but if we let $U_{abc}$ and $V_{abc}$ the components of the $T^*M$-valued $2$-form obtained by lowering the index $m$ to an index $a$, then $U_{abc}$ and $V_{abc}$ are not anti-symmetric in $a$ and $b$, unlike the form $(d^c\o)_{abc}$. The tensors $N$, $U$, and $V$ satisfy the following Bianchi-type identities
\bea
&&
N_{ijk}+N_{jki}+N_{kij}=0,\label{nbianchi}\\
&&
U_{ijk}+U_{jki}+U_{kij}=(d^c\o)_{ijk},\quad
V_{ijk}+V_{jki}+V_{kij}=\frac{1}{2}(d^c\o)_{ijk}.
\eea

Given an almost Hermitian structure $(g,\o,J)$, the Gauduchon line of connections is a line of connections preserving all of $(g,\o,J)$ which passes through the Chern connection and the projected Levi-Civita connection. If we denote the Levi-Civita connection by $\na$, since $d^c\omega$ has only type $(2,1)+(1,2)$-components, the Gauduchon line can be parameterized by a real parameter $t$, and the corresponding connection ${\frak D}^t$ is given by
\bea\label{frakD}
{\frak D}_i^tX^m=\na_iX^m+g^{mk}(-N_{ijk}-V_{ijk}+tU_{ijk})X^j.
\eea
Equivalently, if we express a connection $D$ in terms of its connection form $\Gamma(D)^m{}_{ij}$,
\bea
D_iX^m=\p_iX^m+\Gamma(D)^m{}_{ij}X^j\nonumber
\eea
then we have
\bea
\label{Gamma}
\Gamma({\frak D}^t)^m{}_{ij}=\Gamma(\na)^m{}_{ij}+g^{mk}(-N_{ijk}-V_{ijk}+tU_{ijk}).
\eea
Since the torsion of a connection $D$ is given by
\bea
\label{T0}
T(D)^m{}_{ij}=\Gamma(D)^m{}_{ij}-\Gamma(D)^m{}_{ji}\nonumber
\eea
and the Levi-Civita connection has zero torsion, it is readily follows that
\bea
\label{T1}
T(\mathfrak D^t)^m{}_{jk}=N^m{}_{jk}+(t-1)U^m{}_{jk}+2t V^m{}_{jk}.
\eea

The two values of $t$ of particular interest to us are:

\smallskip
$\bullet$  $t=0$ : this is the so-called projected Levi-Civita connection (a.k.a. the first canonical connection), which we shall henceforth denote by just ${\frak D}^0={\frak D}$.

\smallskip

$\bullet$ $t=1$: this is the Chern connection $\na^C$, also characterized by the condition that $\na_{\bar U}^CV=[\bar U,V]^{1,0}$, for any sections $U,V$ of $T^{1,0}M$. Here
we have set ${\bf C}\otimes TM=T^{1,0}M\oplus T^{0,1}M$ and used $J$ to identify $TM$ with $T^{1,0}M$. The expression $[\bar U,V]^{1,0}$ denotes the $(1,0)$-component of $[\bar U,V]$.

\smallskip
The value $t=-1$ gives the Bismut connection, but we shall not need it in this paper.

\subsubsection{A convenient notation for the action of $J$}

For the convenience of later use, we use the following abbreviations
\begin{eqnarray}
(JV)^k=J^k_{~j}V^j=:V^{Jk},\qquad
(JW)_m=W_jJ^j_{~m}=:W_{Jm}.
\end{eqnarray}
For example, the operator ${\cal M}$ acting on a $TM$-valued $2$-form $\Psi$ introduced earlier in (\ref{M0}) and (\ref{M1}) can now be expressed as
\bea
\label{M2}
({\cal M}\Psi)^m{}_{bc}=\Psi^m{}_{Jb,Jc}.
\eea
and (\ref{UV0}) as
\bea
\label{UV1}
U^m{}_{bc}={1\over 4}((d^c\o)^m{}_{bc}+(d^c\o)^m{}_{Jb,Jc}),\qquad
V^m{}_{bc}={1\over 4}((d^c\o)^m{}_{bc}-(d^c\o)^m{}_{Jb,Jc}).\label{vexp}
\eea
As another example, since $d^c\o$ has no $(0,3)+(3,0)$ component, it satisfies
\be
(d^c\o)_{Ji,j,k}+(d^c\o)_{i,Jj,k}+(d^c\o)_{i,j,Jk}=(d^c\o)_{Ji,Jj,Jk}.\label{30p}
\ee
When summing over repeated indices, $J$ can be raised or lowered as follows
\begin{equation}
X^k\alpha_{Jk}=X^kJ^l_{~k}\alpha_l=X^{Jk}\alpha_k.\nonumber
\end{equation}
Moreover, we can insert $J$ in summation indices at the cost of adding a minus sign:
\be
X^k\alpha_k=-X^{Jk}\alpha_{Jk}.\nonumber
\ee
In this short-hand notation, we have\footnote{Here  $(\o^{jk})$ denotes the inverse matrix of $(\o_{jk})$, $\o^{jk}\o_{kl}=\delta^j{}_l$. It is not the tensor one gets by raising indices using $g_{jk}$. In fact $\o^{jk}=-g^{ja}\o_{ab}g^{bk}$.}
\bea
\o_{jk}=g_{Jj,k}=-g_{j,Jk},&\quad&g_{jk}=\o_{j,Jk}=\o_{k,Jj},\nonumber\\
\o^{jk}=g^{Jj,k}=-g^{j,Jk},&\quad&g^{jk}=\o^{j,Jk}=\o^{k,Jj}.\nonumber
\eea

\subsubsection{The types of $TM$-valued $2$-forms}

All the tensors which we encountered above, namely the Nijenhuis tensor, the tensors $d^c\o$, $U$, and $V$, and the torsion tensors can be viewed as $TM$-valued $2$-forms, or equivalently $3$-tensors which are antisymmetric in the last two slots. Denote this space by $A^2(TM)$ for simplicity.
It is convenient to break up elements of $A^2(TM)$ into simpler components.

\medskip
Recall the involution ${\cal M}$ on $A^2(TM)$ defined by (\ref{M0}) or (\ref{M2}) in components.
Clearly $A^2(TM)$ splits into the direct sum of eigenspaces of $\cM$ with eigenvalue $\pm1$. We shall call the eigenvalue 1 subspace the space of $TM$-valued (1,1)-forms, denoted by $A^{1,1}(TM)$. That is, $\Psi\in A^{1,1}(TM)$ if and only if
\begin{equation}
\Psi^p_{~Jj,Jk}=\Psi^p_{~jk}\quad\textrm{or equivalently}\quad \Psi_{i,Jj,k}+\Psi_{i,j,Jk}=0.\nonumber
\end{equation}
The space of eigenvalue $-1$ can be decomposed further as follows.
We say a $TM$-valued 2-form $\Psi$ is of type (2,0) or (0,2) if
\footnote{To avoid confusion, we stress that this notion is specific to $A^2(TM)$ and is not the same as that of scalar-valued $(2,0)$ or $(0,2)$ forms.}
\begin{equation}
\Psi(JX,Y)=J\Psi(X,Y)\quad\textrm{or}\quad \Psi(JX,Y)=-J\Psi(X,Y).\nonumber
\end{equation}
In this way we have a direct sum decomposition
\begin{equation}
A^2(TM)=A^{1,1}(TM)\oplus A^{2,0}(TM)\oplus A^{0,2}(TM).\nonumber
\end{equation}
In terms of indices we see that
\begin{eqnarray}
\Psi\in A^{2,0}(TM)&\quad\textrm{if}\quad& \Psi^p_{~Jj,k}=\Psi^{Jp}_{~~jk}=\Psi^p_{~j,Jk} ~\textrm{ or }~ \Psi_{Ji,j,k}=-\Psi_{i,Jj,k}=-\Psi_{i,j,Jk},\nonumber\\
\Psi\in A^{0,2}(TM)&\quad\textrm{if}\quad& \Psi^p_{~Jj,k}=-\Psi^{Jp}_{~~jk}=\Psi^p_{~j,Jk} ~\textrm{ or }~ \Psi_{Ji,j,k}=\Psi_{i,Jj,k}=\Psi_{i,j,Jk}.\nonumber
\end{eqnarray}

Returning to the tensors which we have encountered, their types are as follows:

\medskip

$\bullet$ It is readily seen that the Nijenhuis tensor is of type $(0,2)$, therefore any contraction of $N$ using either $g$ or $\o$ yields 0.

\smallskip

$\bullet$ It is easy to see that the tensor $U^m{}_{jk}$ is of type $(1,1)$. As for $V$, we have by definition
\begin{eqnarray}
V_{Ji,j,k}&=&\frac{1}{4}((d^c\o)_{Ji,j,k}-(d^c\o)_{Ji,Jj,Jk})\overset{\textrm{(\ref{30p})}}{=}-\frac{1}{4}((d^c\o)_{i,Jj,k}+(d^c\o)_{i,j,Jk}) \nonumber\\
&=&V_{J(Ji),Jj,k}=-V_{i,Jj,k}.\nonumber
\end{eqnarray}
and thus $V$ is of type (2,0).
\smallskip

$\bullet$ Finally, tensors such as the torsions of connections on the Gauduchon line and consequently also differences of connections, correspond to forms of mixed types, whose decomposition in $\Lambda^{1,1}(TM)\oplus \Lambda^{2,0}(TM)\oplus \Lambda^{0,2}(TM)$ can be read off from formulas such as (\ref{Gamma}) and (\ref{T1}), since we know now the types of $N^m{}_{jk}$, $U^m{}_{jk}$, and $V^m{}_{jk}$.

\subsection{Proof of Theorem \ref{th:Laplacian}}

We can now give the proof of Theorem \ref{th:Laplacian}.
For simplicity, we shall denote in the subsequent calculations $g_\varphi$ and $J_\varphi$ by just $g$ and $J$. Recall that the operator $d^c$ is defined by $d^c=J^{-1}dJ$, and that, for a compatible structure $(\o,g,J)$, we have the identity
\bea\label{dagger-identity}
d^\Lambda:=d\Lambda-\Lambda d=(d^c)^\dagger.
\eea
This identity holds even if $J$ is not integrable. Since $\varphi$ is primitive, we can replace in (\ref{fl:IIA}) $\Lambda d$ by $-d^\Lambda$ and rewrite the equation as
\bea
\p_t\varphi
&=&
d\Lambda d(|\varphi|^2\star\!\varphi)-\rho_A
=
-d(d^c)^\dagger(|\varphi|^2\star\!\varphi)-\rho_A
=
-d(J^{-1}dJ)^\dagger(|\varphi|^2\star\!\varphi)-\rho_A
\nonumber\\
&=&
-dJ^{-1}d^\dagger J(|\varphi|^2\star\!\varphi)-\rho_A.
\eea
Here we have used the fact that the adjoint $J^\dagger$ of $J$ with respect to $g$ is $J^{-1}$, since $J$ is an isometry. We shall see later that $\star\varphi=J\varphi$, so that the above equation can be rewritten as
\bea
\p_t\varphi=dJd^\dagger(|\varphi|^2\varphi)-\rho_A.
\eea
Thus the theorem would be proved once we can establish the following lemma:

\begin{lemma}
\label{Jddagger}
Let $(M,g,J,\o)$ be a $6$-dimensional compact almost-Hermitian manifold with $d\o=0$. We have
\bea
Jd^\dagger\varphi=-d^\dagger\varphi+2N^\dagger\cdot\varphi
\eea
where $N^\dagger$ is the operator defined in (\ref{Ndagger}).
\end{lemma}

\medskip
To establish this, we begin by recalling that the adjoint of the operator $d$ on $3$-forms is given by
\bea
(d^\dagger\varphi)_{\alpha\beta}=-\na^\gamma\varphi_{\gamma\alpha\beta}
\eea
where $\na$ denotes the Levi-Civita connection of $g_{ij}$, which has no torsion. We need to apply the operator $J$ to both sides of this equation. For this, we need in turn the following lemma:

\begin{lemma}
\label{Jddagger1}
\bea
\nabla_k J^a{}_b = -2N_{Jk}{}^a{}_b.
\eea
\end{lemma}

The point of this lemma is that the Levi-Civita connection does not necessarily respect the almost-complex structure $J$. However, we can write it in terms of the projected Levi-Civita connection ${\frak D}$ which does, at the cost of having to handle in addition terms coming from difference of connections, which gives us the Nijenhuis tensor. This lemma follows directly from (\ref{Gamma}), where $U=V=0$ as $\omega$ is closed.

\if $A^\gamma{}_{jk}$. Explicitly, we write
\bea
J^\beta{}_j\na_k\varphi_{\gamma\alpha\beta}
=
{\frak D}_k(J^\beta{}_j\varphi_{\gamma\alpha\beta})
-
J^\beta{}_jA^\lambda{}_{k\gamma}\varphi_{\lambda\alpha\beta}
-
A^\lambda{}_{k\alpha}\varphi_{\gamma\lambda\beta}
-
A^\lambda{}_{k\beta}\varphi_{\gamma \alpha\lambda}
\nonumber
\eea
in view of the fact that $J$ is covariantly constant with respect to ${\frak D}$. This can be compared with
\bea
\na_k(J^\beta{}_j\varphi_{\gamma\alpha\beta})
=
{\frak D}_k(J^\beta{}_j\varphi_{\gamma\alpha\beta})
-
A^\lambda{}_{kj}(J^\beta{}_\lambda\varphi_{\gamma\alpha\beta})
-
A^\lambda{}_{k\gamma}(J^\beta{}_j\varphi_{\lambda\alpha\beta})
-
A^\lambda{}_{k\alpha}(J^\beta{}_j\varphi_{\gamma\lambda\beta})
\nonumber
\eea
and hence
\bea
J^\beta{}_j\na_k\varphi_{\gamma\alpha\beta}
&=&
\na_k(J^\beta{}_j\varphi_{\gamma\alpha\beta})
-A^\lambda{}_{k\gamma}\varphi_{\lambda\alpha, Jj}+A^\lambda{}_{kj}\varphi_{\gamma\alpha,J\lambda}
\nonumber\\
&=&
\na_k(J^\beta{}_j\varphi_{\gamma\alpha\beta})
-N^\lambda{}_{k\gamma}\varphi_{\lambda\alpha, Jj}+N^\lambda{}_{kj}\varphi_{\gamma\alpha,J\lambda}
\eea
since $d\o=d^c\o=0$, and consequently $A^\lambda{}_{k\gamma}=N^\lambda{}_{k\gamma}$. Contracting with $g^{k\gamma}$, we obtain the first formula in Lemma \ref{Jddagger1}. The second formula follows by iterating the first,
\bea
J^\alpha{}_kJ^\beta{}_j\na^\gamma\varphi_{\gamma\alpha\beta}
&=&
J^\alpha{}_k\na^\gamma(J^\beta{}_j\varphi_{\gamma\alpha \beta})+N^{\lambda\gamma}{}_j\varphi_{\gamma, Jk,J\lambda}
\nonumber\\
&=&
-\na^\gamma(J^\alpha{}_jJ^\beta{}_k\varphi_{\gamma jk})
-
N^{\lambda\mu}{}_k\varphi_{\mu,Jj,J\lambda}
+
N^{\lambda\mu}{}_j\varphi_{\mu, Jk, J\lambda}.
\eea
\fi

Returning to the proof of Lemma \ref{Jddagger}, since $\varphi_{\alpha,J\beta,J\gamma}=-\varphi_{\alpha\beta\gamma}$ by Lemma \ref{Jddagger1} and \ref{phihat}, we find
\bea
(Jd^\dagger\varphi)_{kj}
&=&
-J^\alpha{}_kJ^\beta{}_j\na^\gamma\varphi_{\gamma\alpha\beta}
\nonumber\\
&=&-\nabla^\gamma(J^\alpha{}_kJ^\beta{}_j\varphi_{\gamma\alpha\beta})
-2(N_{J\gamma}{}^\alpha{}_k J^\beta{}_j + N_{J\gamma}{}^\beta{}_j J^\alpha{}_k)\varphi^\gamma{}_{\alpha\beta}
\nonumber\\
&=&
-\na^\gamma\varphi_{\gamma, Jk,Jj}
-
2(N^{\lambda\mu}{}_k\varphi_{\mu j\lambda}-N^{\lambda\mu}{}_j\varphi_{\mu k\lambda})
\nonumber\\
&=&
-(d^\dagger\varphi)_{kj}+2(N^\dagger\cdot\varphi)_{kj}.
\eea
This completes the proof of Lemma \ref{Jddagger}. Replacing $\varphi$ by $|\varphi|^2\varphi$ in Lemma \ref{Jddagger}, we also obtain Theorem \ref{th:Laplacian}.

\section{The principal symbol of the Type IIA equation}
\setcounter{equation}{0}
\label{s:symbol}

Our next task is to identify the symbol and the eigenvalues of the Type IIA equation.

\subsection{Almost-complex structures and $3$-forms}

For this and for the rest of the paper, we make essential use of Hitchin's construction of an almost-complex structure $J_\varphi$ and a metric $g_\varphi$ from a $3$-form $\varphi$ on a $6$-dimensional manifold \cite{Hitchin}. We begin by recalling the results that we need.

\subsubsection{Hitchin's construction}\label{hitchin}

Let $V$ be a 6-dimensional oriented vector space over $\R$. Following Hitchin \cite{Hitchin}\footnote{Here we adopt the convention used in \cite{Fei, FeiD}.}, for any $3$-form $\varphi\in\Lambda^3V^*$, one can define a linear map $K_\varphi:V\to \Lambda^5V^*\cong V\otimes\Lambda^6V^*$ by
\be
K_\varphi(v)=-\iota_v\varphi\wedge\varphi=-e_i\otimes e^i\wedge\iota_v\varphi\wedge\varphi,
\ee
where $\{e_i\}$ is an arbitrary basis of $V$ and $\{e^i\}$ its dual basis in $V^*$. It follows that
\be
\lambda_\varphi:=\frac{1}{6}\tr_VK_\varphi^2\in(\Lambda^6V^*)^2
\ee
is well-defined and it makes sense to talk about the sign of $\lambda_\varphi$. In general $\lambda_\varphi$ is a homogeneous degree 4 polynomial in the components of $\varphi$. When $\lambda_\varphi<0$, as $V$ is oriented, one can take $\sqrt{-\lambda_\varphi}\in\Lambda^6V^*$ to be the positive square root of $-\lambda_\varphi$. It is proved in \cite{Hitchin}
\be
J_\varphi:=\frac{K_\varphi}{\sqrt{-\lambda_\varphi}}:V\to V
\ee
defines a complex structure on $V$. Note that $\lambda_\varphi<0$ is an open condition. In fact, the set $\{\varphi\in\Lambda^3V^*:\lambda_\varphi<0\}$ forms an open orbit in $\Lambda^3V^*$ of the natural $GL(V)$-action. Furthermore, there is a basis $\{e^i\}$ of $V^*$ where
$\varphi$ takes the following ``canonical form'',
\be \label{phi-can-form}
\varphi=\Re(e^1+ie^2)\wedge(e^3+ie^4)\wedge(e^5+ie^6)=e^{135}-e^{146}-e^{245}-e^{236}
\ee
and $e^{123456}$ defines a positive volume form. In this basis,
one can easily check that $J_\varphi e_{2k-1}=e_{2k}$ and $J_\varphi e_{2k}=-e_{2k-1}$ for $k=1,2,3$. Therefore $\varphi$ is the real part of a $(3,0)$-form with respect to the complex structure $J_\varphi$. It follows that the imaginary part
\be \label{phi-hat-can-form}
\hat{\varphi}=\Im(e^1+ie^2)\wedge(e^3+ie^4)\wedge(e^5+ie^6)=e^{136}+e^{145}+e^{235}-e^{246}
\ee
is also determined by $\varphi$ through $\hat{\varphi}=J_\varphi\varphi$, meaning that
\be
\hat\varphi(X,Y,Z)=\varphi(J_\varphi X,J_\varphi Y,J_\varphi Z)\nonumber
\ee
for any $X,Y,Z\in V$. Furthermore two forms $\varphi$ and $\tilde\varphi$ define the same complex structure
$J_\varphi$ if and only if they are related by $\C^*$-action
\be
\tilde\varphi=\rho\cdot\Re(e^{-i\theta}(\varphi+i\hat\varphi)).\nonumber
\ee

\subsubsection{$J_\varphi$ and symplectic structures}

Now let us assume that $V$ is equipped with a symplectic form $\omega\in\Lambda^2V^*$ so $V$ is canonically oriented by $\omega^3/3!$. A natural question is when the symplectic form $\omega$ is invariant under the induced complex structure $J_\varphi$. The answer is very simple:

\begin{lemma}
$\omega$ is $J_\varphi$-invariant if and only if $\varphi$ is primitive in the sense that $\omega\wedge\varphi=0$.
\end{lemma}
{\it Proof}. As stated above, we can choose a basis such that $\varphi=e^{135}-e^{146}-e^{245}-e^{236}$. In this coordinate, we may write $\omega=a_{ij}e^i\wedge e^j$ with $a_{ij}=-a_{ji}$. The condition that $\o(e_i,e_j)=\o(J_\varphi e_i,J_\varphi e_j)$ for any $i$ and $j$ is equivalent to the following system of linear equations
\bea
&&
a_{13}=a_{24},\qquad a_{14}=-a_{23},\qquad a_{15}=a_{26}, \nonumber\\
&&
a_{16}=-a_{25},\qquad a_{35}=a_{46},\qquad a_{36}=-a_{45}.
\eea
These are exactly the equations for $\varphi$ being primitive in the sense that $\omega\wedge\varphi=0$. Q.E.D.

\medskip
For primitive $\varphi$, we can consider then the Hermitian form $g_\varphi(X,Y)=\o(X,J_\varphi Y)$. We shall say that $\varphi$ is positive if $g_\varphi$ is positive, in which case $g_\varphi$ is a metric, and the triple $(\o,J_\varphi, g_\varphi)$ is compatible \footnote{This notion of positivity is different from the one defined in \cite{Hitchin}, which does not involve a symplectic form.}. The positivity of $\varphi$ is an open condition.
Once we have a Riemannian metric and an orientation, we have the associated Hodge star operator $\star$. It is straightforward to check that
\be
\star\!\varphi=\hat\varphi,\quad \star\hat\varphi=-\varphi.
\ee
In particular we know that $\hat\varphi$ is primitive if $\varphi$ is primitive.

\medskip
Altogether, assuming the presence of $\o$ and that $\varphi$ is primitive and positive, we may upgrade the previous choice of orthonormal basis $\{e_j\}_{j=1}^6$ to the following useful statement:

\begin{lemma}(Normal form of $\varphi$)
\label{up canonical form}~\\
There exists an orthonormal basis $\{e_j\}_{j=1}^6$ of $V$ (with respect to $g$) such that
\bea
\o&=&e^1\wedge e^2+e^3\wedge e^4+e^5\wedge e^6,\\
\varphi&=&M\Re(e^1+ie^2)\wedge(e^3+ie^4)\wedge(e^5+ie^6),
\eea
where $M=\dfrac{1}{2}|\varphi|>0$. It follows that
\be
\sqrt{-\lambda_\varphi}=\frac{1}{2}|\varphi|^2\frac{\omega^3}{3!}.
\ee
\end{lemma}

Using this upgraded version of canonical form of $\varphi$, one can check that the following key formula for the metric $g_\varphi$ holds:

\begin{lemma}
In any coordinate system, we can write
\be
\label{g}
(g_{\varphi})_{ij}=-|\varphi|^{-2} \varphi_{iab}\varphi_{jcd}\o^{ac}\o^{bd}=2|\varphi|^{-2}\frac{\iota_i\varphi\wedge\iota_j\varphi\wedge\o}{\o^3/3!}.
\ee
\end{lemma}

The metric $\tilde g_\varphi$ introduced in (\ref{tildeg0}) is then given by
\be
(\tilde g_\varphi)_{ij}=-\varphi_{iab}\varphi_{jcd}\o^{ac}\o^{bd}\label{tildemetric}.
\ee
Clearly $\tilde g_\varphi$ is conformal to $g_\varphi$ and its associated K\"ahler form is $\tilde\o_\varphi=|\varphi|^2\o$.
Since $|\varphi|^2$ is the square root of a complicated homogeneous degree 4 polynomial in components of $\varphi$, the metric $\tilde g_\varphi$ has the advantage that its expression is algebraic in $\varphi$, which makes it much easier to compute with. Also, the volume form of $g_\varphi$ is just $\o^3/3!$, but we can recapture $|\varphi|^2$ from the volume form of $\tilde g_\varphi$.

\subsubsection{Basic identities}

As always we fix a symplectic form $\omega$. Let $\varphi$ be a positive primitive 3-form so that it defines an almost complex structure $J_\varphi$ compatible with $\omega$. For simplicity, we shall denote in the subsequent calculations $g_\varphi$ and $J_\varphi$ by just $g$ and $J$. Furthermore we can define another primitive 3-form $\hat\varphi=J\varphi=\star\varphi$ such that $\Omega=\varphi+i\hat\varphi$ is a nowhere vanishing (3,0)-form with respect to $J$, i.e. a complex volume form that trivializes the canonical bundle of $(M,J)$.

\begin{lemma}\label{phihat}
The $3$-forms $\varphi$ and $\hat{\varphi}$ are related to each other by
\begin{eqnarray}
&&\varphi_{ijk}=\hat\varphi_{Ji,j,k}=\hat\varphi_{i,Jj,k}=\hat\varphi_{i,j,Jk}=-\varphi_{Ji,Jj,k}=-\varphi_{Ji,j,Jk}=-\varphi_{i,Jj,Jk} =-\hat\varphi_{Ji,Jj,Jk}\nonumber\\
&&\hat\varphi_{ijk}=-\varphi_{Ji,j,k}=-\varphi_{i,Jj,k}=-\varphi_{i,j,Jk}=-\hat\varphi_{Ji,Jj,k}=-\hat\varphi_{Ji,j,Jk}=-\hat\varphi_{i,Jj,Jk} =\varphi_{Ji,Jj,Jk}.\nonumber
\end{eqnarray}
\end{lemma}
{\it Proof}. Since $\varphi+i\hat{\varphi}$ is of type (3,0), we have
$\iota_{\p_k+iJ\p_k}(\varphi+i\hat\varphi)=0$.
By taking the real and imaginary parts of the above equation and its iterations gives the desired identities. Q.E.D.\\

Since $\varphi$ and $\hat\varphi$ are type $(3,0)+(0,3)$-forms, for any 1-form $\mu$, we know that both $\mu\wedge\varphi$ and $\mu\wedge\hat\varphi$ are of type $(3,1)+(1,3)$, so
\begin{equation}
\mu\wedge\varphi=-J(\mu\wedge\varphi)=-J\mu\wedge\hat\varphi.\label{rotate}
\end{equation}
It is not hard to verify that wedging with $\varphi$ or $\hat\varphi$ gives an isomorphism from the space of real $1$-forms to the space of real $(3,1)+(1,3)$-forms. Note that the primitiveness of $\varphi$ with respect to $\o$ implies the primitiveness of $\varphi$ with respect to $\tilde\o$, and hence $\tilde{\o}^{ji}\varphi_{ijk}=0$, or equivalently,
\begin{eqnarray}
&&\tilde{\o}_{ij}\varphi_{klm}-\tilde{\o}_{kj}\varphi_{ilm}-\tilde{\o}_{lj}\varphi_{kim}-\tilde{\o}_{mj}\varphi_{kli}-\tilde{\o}_{ik}\varphi_{jlm}\nonumber\\ &&-\tilde{\o}_{il}\varphi_{kjm}-\tilde{\o}_{im}\varphi_{klj}+\tilde{\o}_{kl}\varphi_{ijm}+\tilde{\o}_{ml}\varphi_{kji}+\tilde{\o}_{mk}\varphi_{jli}=0. \label{primitive}
\end{eqnarray}

We also have the following simple lemma:

\begin{lemma}
\label{lm:Omega}
The following are equivalent:

{\rm (a)}  $d\hat\varphi=0$;

{\rm (b)} The almost-complex structure $J$ is integrable.

\noindent
In both cases, the $(3,0)$-form $\Omega$ is holomorphic and the form $\varphi$ is harmonic.

\end{lemma}

\smallskip
\noindent{\it Proof.} If $J$ is integrable, then $\p\Omega=0$ since it is a $(4,0)$-form in a $3$-dimensional complex manifold. Similarly for $\bar\p\bar\Omega$. But $\varphi={1\over 2}(\Omega+\bar\Omega)$, so $d\varphi=0$ implies $0=\bar\p\Omega+\p\bar\Omega$, and hence $\Omega$ is holomorphic and $d\Omega=0$. This implies that $d\hat\varphi=0$.

Conversely, if $d\hat\varphi=0$, then we know that $\Omega=\varphi+\sqrt{-1}\hat\varphi$ is also closed. Therefore for any $(1,0)$-form $\lambda$, we know $\lambda\wedge\Omega=0$, hence
\bea
0=d(\lambda\wedge\Omega)=d\lambda\wedge\Omega,\nonumber
\eea
so we deduce that $d\lambda$ has no $(0,2)$-components, which implies that $J$ is integrable by Frobenius theorem. Q.E.D.

\medskip
The defining equation for $\tilde{g}$ tells us the effect of contracting twice with $\o^{ij}$ a quadratic polynomial in $\varphi$. It actually follows from a stronger identity with only one contraction, and which can be verified explicitly using the normal form of $\varphi$ in Lemma \ref{up canonical form}. Using the fact that $g^{ij}=\o^{i,Jj}$, we can readily deduce the effect of contracting with $g^{ij}$. We summarize these contractions in the following lemma:

\begin{lemma}\label{contract1}
The following quadratic identities hold:
\bea
\omega^{ij}\varphi_{iab}\varphi_{jcd}&=&\frac{1}{4}(\omega_{ac}\tilde{g}_{bd}-\omega_{bc}\tilde{g}_{ad}-\omega_{ad}\tilde{g}_{bc}+\omega_{bd}\tilde{g}_{ac})
\nonumber\\
g^{ij}\varphi_{iab}\varphi_{jcd}&=&\frac{1}{4}(g_{ac}\tilde g_{bd}-g_{bc}\tilde g_{ad}+\o_{ad}\tilde\o_{bc}-\o_{bd}\tilde\o_{ac})\\
g^{ac}g^{bd}\varphi_{iab}\varphi_{jcd}&=&\tilde{g}_{ij}.
\eea
\end{lemma}

\subsection{The variation $\delta(\hat\varphi)$}

The key variational formula for $\hat\varphi$ is given by the following lemma:

\begin{lemma}\label{variation}
\be
\delta(\hat\varphi)=-J_\varphi(\delta\varphi)+2\frac{\delta\varphi\wedge\varphi}{\varphi\wedge\hat\varphi}\varphi +2\frac{\delta\varphi\wedge\hat\varphi}{\varphi\wedge\hat\varphi}\hat\varphi.
\ee
When $\varphi$ is primitive and positive with respect to a symplectic form $\o$, the above formula can be written as
\bea
\delta(\hat\varphi)=-J_\varphi(\delta\varphi)-\frac{2(\delta\varphi,\hat\varphi)}{|\varphi|^2}\varphi+\frac{2(\delta\varphi,\varphi)}{|\varphi|^2}\hat\varphi.
\eea
\end{lemma}
{\it Proof.} This is a purely linear algebra problem. We break its proof into two steps.

\noindent $\bullet$ Step 1: The primitive case.

Choose a nondegenerate $(1,1)$-form $\o$ compatible with $J_\varphi$. We shall first prove Lemma \ref{variation} under the assumption that $\delta\varphi$ is primitive with respect to $\o$. By Lemma \ref{up canonical form}, one can find an orthonormal basis $\{e^i\}_{i=1}^6$ of $V^*$ such that
\bea
\o&=&e^1\wedge e^2+e^3\wedge e^4+e^5\wedge e^6,\nonumber\\
\varphi&=&M(e^1\wedge e^3\wedge e^5-e^1\wedge e^4\wedge e^6-e^2\wedge e^4\wedge e^5-e^2\wedge e^3\wedge e^6),\nonumber\\
\star\varphi&=&M(e^1\wedge e^3\wedge e^6+e^1\wedge e^4\wedge e^5+e^2\wedge e^3\wedge e^5-e^2\wedge e^4\wedge e^6)\nonumber,
\eea
where $M=\dfrac{1}{2}|\varphi|>0$. For simplicity, we denote $e^{123456}=\o^3/3!$ by $\epsilon$. Let $\delta\varphi=\mu=\dfrac{1}{3!}\mu_{ijk}e^i\wedge e^j\wedge e^k$. Straightforward computation gives us
\be
\delta(K_\varphi)=2M\epsilon
\begin{bmatrix}
A_1 & \mu_{135}-\mu_{146} & -\mu_{125} & \mu_{126} & \mu_{123} & -\mu_{124} \\ \mu_{236}+\mu_{245} & -A_1 & \mu_{126} & \mu_{125} & -\mu_{124} & -\mu_{123} \\ \mu_{345} & -\mu_{346} & A_2 & \mu_{135}-\mu_{236} & -\mu_{134} & \mu_{234} \\ -\mu_{346} & -\mu_{345} & \mu_{245}+\mu_{146} & -A_2 & \mu_{234} & \mu_{134} \\ -\mu_{356} & \mu_{456} & \mu_{156} & -\mu_{256} & A_3 & \mu_{135}-\mu_{245} \\ \mu_{456} & \mu_{356} & -\mu_{256} & -\mu_{156} & \mu_{236}+\mu_{146} & -A_3
\end{bmatrix},\nonumber
\ee
where
\bea
A_1&=&\frac{\mu_{246}+\mu_{136}+\mu_{145}-\mu_{235}}{2},\nonumber\\
A_2&=&\frac{\mu_{246}+\mu_{136}-\mu_{145}+\mu_{235}}{2},\nonumber\\
A_3&=&\frac{\mu_{246}-\mu_{136}+\mu_{145}+\mu_{235}}{2}.\nonumber
\eea
It follows that
\be
\delta(K_\varphi^2)=4M^3\epsilon^2(\mu_{236}+\mu_{245}-\mu_{135}+\mu_{146})\textrm{Id}_V=4M^2\epsilon^2\frac{\hat\varphi\wedge\mu}{\epsilon}\textrm{Id}_V.
\ee
Therefore
\bea
\delta(-\lambda_\varphi)&=&4M^2\epsilon^2\frac{\mu\wedge\hat\varphi}{\epsilon}=4M^2\epsilon^2(\mu,\varphi)\\
\delta\sqrt{-\lambda_\varphi}&=&\frac{1}{2}\delta(|\varphi|^2)\epsilon=\mu\wedge\hat\varphi=(\delta\varphi,\varphi)\epsilon,
\eea
which agrees with Hitchin's formula \cite[Proposition 4]{Hitchin}. As a consequence, we have
\be
\delta(J_\varphi)=\frac{\delta(K_\varphi)}{2M^2\epsilon}-\frac{\mu\wedge\hat\varphi}{2M^2\epsilon}J_\varphi.
\ee

For simplicity of notation, we introduce the following
\bea
F&:=&\frac{\delta(K_\varphi)}{2M\epsilon},\nonumber\\
\varphi_0&:=&e^1\wedge e^3\wedge e^5-e^1\wedge e^4\wedge e^6-e^2\wedge e^4\wedge e^5-e^2\wedge e^3\wedge e^6,\nonumber\\
\hat\varphi_0&:=&e^2\wedge e^3\wedge e^6+e^1\wedge e^4\wedge e^5+e^2\wedge e^3\wedge e^5-e^2\wedge e^4\wedge e^6,\nonumber\\
dz^k&:=&e^{2k-1}+ie^{2k},~k=1,2,3.\nonumber
\eea
Again by straightforward calculation, we get
\bea
F.dz^1&=&iB_0dz^1+(A_1+iB_1)d\bar z^1+C_3d\bar z^2+C_2d\bar z^3,\nonumber\\
F.dz^2&=&C_3d\bar z^1+iB_0dz^2+(A_2+iB_2)d\bar z^2+C_1d\bar z^3,\nonumber\\
F.dz^3&=&C_2d\bar z^1+C_1d\bar z^2+iB_0dz^3+(A_3+iB_3)d\bar z^3.\nonumber
\eea
where
\bea
B_0&=&\frac{1}{2}(\mu_{135}-\mu_{146}-\mu_{236}-\mu_{245})=\frac{\mu\wedge\hat{\varphi}_0}{2\epsilon},\nonumber\\
B_1&=&\frac{1}{2}(\mu_{135}-\mu_{146}+\mu_{236}+\mu_{245}),\nonumber\\
B_2&=&\frac{1}{2}(\mu_{135}+\mu_{146}-\mu_{236}+\mu_{245}),\nonumber\\
B_3&=&\frac{1}{2}(\mu_{135}+\mu_{146}+\mu_{236}-\mu_{245}),\nonumber\\
C_1&=&\mu_{156}-i\mu_{256}=-\mu_{134}+i\mu_{234},\nonumber\\
C_2&=&-\mu_{356}+i\mu_{456}=\mu_{123}-i\mu_{124},\nonumber\\
C_3&=&\mu_{345}-i\mu_{346}=-\mu_{125}+i\mu_{126}.\nonumber
\eea
Here to obtain expressions of $C_j$ one makes use of primitiveness of $\mu$. For completeness we also introduce
\be
A_0=\frac{\mu_{246}-\mu_{235}-\mu_{136}-\mu_{145}}{2}=\frac{\mu\wedge\varphi_0}{2\epsilon}.\nonumber
\ee
Collecting all these together, we get
\bea
F.(\varphi_0+i\hat\varphi_0)&=&F.(dz^1\wedge dz^2\wedge dz^3)\nonumber\\
&=&F.(dz^1)\wedge dz^2\wedge dz^3+dz^1\wedge F.(dz^2)\wedge dz^3+dz^1\wedge dz^2\wedge F.(dz^3)\nonumber\\
&=&3iB_0dz^1\wedge dz^2\wedge dz^3+(A_1+iB_1)d\bar z^1\wedge dz^2\wedge dz^3\nonumber\\ &&+(A_2+iB_2)dz^1\wedge d\bar z^2\wedge dz^3+(A_3+iB_3)dz^1\wedge dz^2\wedge d\bar z^3\nonumber\\ && +C_1dz^1\wedge(d\bar z^3\wedge dz^3+dz^2\wedge d\bar z^2)-C_2dz^2\wedge(d\bar z^3\wedge dz^3+dz^1\wedge d\bar z^1)\nonumber\\ &&+C_3dz^3\wedge(dz^1\wedge d\bar z^1+d\bar z^2\wedge dz^2).\label{faction}
\eea

Notice we can express $\mu$ as
\bea
\mu&=&\nonumber-\frac{i}{4}\bigg [(C_1dz^1+\overline{C_1}d\bar z^1)\wedge(d\bar z^3\wedge dz^3+dz^2\wedge d\bar z^2)\nonumber\\
&&-(C_2dz^2+\overline{C_2}d\bar z^2)\wedge(d\bar z^3\wedge dz^3+dz^1\wedge d\bar z^1)\nonumber\\
&&+(C_3dz^3+\overline{C_3}d\bar z^3)\wedge(dz^1\wedge d\bar z^1+d\bar z^2\wedge dz^2)\bigg ]+\frac{B_0}{2}\varphi_0-\frac{A_0}{2}\hat{\varphi}_0\nonumber\\
&&-\frac{i}{4}(A_1+iB_1)d\bar z^1\wedge dz^2\wedge dz^3+\frac{i}{4}(A_1-iB_1)d z^1\wedge d\bar z^2\wedge d\bar z^3\nonumber\\
&&-\frac{i}{4}(A_2+iB_2)dz^1\wedge d\bar z^2\wedge dz^3+\frac{i}{4}(A_2-iB_2)d\bar z^1\wedge d z^2\wedge d\bar z^3\nonumber\\
&&-\frac{i}{4}(A_3+iB_3)dz^1\wedge dz^2\wedge d\bar z^3+\frac{i}{4}(A_3-iB_3)d\bar z^1\wedge d\bar z^2\wedge d z^3.\nonumber
\eea
By taking real part of (\ref{faction}) we get
\bea
F.\varphi_0&=&2J_\varphi\mu-A_0\varphi_0-4B_0\hat\varphi_0.\nonumber
\eea
Since
\bea
J_\varphi.\varphi_0&=&\varphi_0(J_\varphi\cdot,\cdot,\cdot)+\varphi_0(\cdot,J_\varphi\cdot,\cdot)+\varphi_0(\cdot,\cdot,J_\varphi\cdot) =-3\hat\varphi_0,\nonumber\\
\delta(J_\varphi)&=&\frac{1}{M}(F-B_0J_\varphi),\nonumber
\eea
it follows that
\bea
\delta(J_\varphi).\varphi_0&=&\frac{1}{M}(2J_\varphi\mu-A_0\varphi_0-B_0\hat{\varphi}_0).
\eea
Consequently we obtain
\be
\delta\hat{\varphi}=J_\varphi\delta\varphi-\delta(J_\varphi).\varphi=-J_\varphi\mu+A_0\varphi_0+B_0\hat\varphi_0.\nonumber
\ee
Rewrite this equation in a coordinate-free manner, we obtain the desired formula.

\noindent$\bullet$ Step 2: The general case.

Choosing $\omega$ as before, a general variation $\delta\varphi$ takes the form
\bea
\delta\varphi=\mu+\omega\wedge\lambda,\nonumber
\eea
where $\mu$ is primitive with respect to $\o$ and $\lambda$ is a 1-form. By linearity, we only need to prove our formula for $\delta\varphi=\omega\wedge\lambda$. By symmetry, we may assume that $\lambda=Ne^1$ for some number $N$. Therefore $\o\wedge\lambda=Ne^1\wedge(e^{34}+e^{56})$ is a linear combination of $e^1\wedge(e^{34}-e^{56})$ and $e^1\wedge(Ke^{34}-e^{56})$ for some constant $K>1$. Notice that $e^1\wedge(e^{34}-e^{56})$ is primitive with respect to $\o$, and $e^1\wedge(Ke^{34}-e^{56})$ is primitive with respect to another $J_\varphi$-compatible $(1,1)$-form $\o'=e^{12}+Ke^{34}+e^{56}$. By linearity we reduce Step 2 to Step 1 with $\o$ replaced by $\o'$. Q.E.D.

\subsection{The eigenvalues of the principal symbol}

With all preparations from last section, we are now ready to compute the principal symbol of the Type IIA flow (\ref{fl:IIA}) without source. When $\varphi$ is primitive, the right hand side of the Type IIA flow is also primitive, so we only need to consider primitive variations in Lemma \ref{variation}, which takes the form
\bea
\label{linearize}
\delta(|\varphi|^2\hat\varphi)=-|\varphi|^2J(\delta\varphi)-2(\delta\varphi,\hat\varphi)\varphi+4(\delta\varphi,\varphi)\hat\varphi.\label{linearize}
\eea
Thus the symbol of the leading term in the Type IIA flow is given by
\be
\delta\varphi\mapsto\xi\wedge\Lambda\big\{\xi\wedge(-|\varphi|^2J(\delta\varphi)-2(\delta\varphi,\hat\varphi)\varphi+4(\delta\varphi,\varphi)\hat\varphi)\big\}.
\ee
and whether the flow is parabolic or not, depends on the eigenvalues of this operator.
Since by our assumption, the right hand side of the flow is primitive and admits an integrability operator $d$. Thus by the Hamilton-Nash-Moser theorem \cite{Ha1}, we can restrict $\delta\varphi$ to the space
\be
W=\{\delta\varphi\in\Lambda^3V^*:\xi\wedge\delta\varphi=0,\Lambda(\delta\varphi)=0\}.\nonumber
\ee
As before, we may choose an orthonormal basis $\{e^i\}_{i=1}^6$ of $V^*$ such that
\bea
\o&=&e^1\wedge e^2+e^3\wedge e^4+e^5\wedge e^6,\nonumber\\
\varphi&=&M(e^1\wedge e^3\wedge e^5-e^1\wedge e^4\wedge e^6-e^2\wedge e^4\wedge e^5-e^2\wedge e^3\wedge e^6),\nonumber\\
\star\varphi&=&M(e^1\wedge e^3\wedge e^6+e^1\wedge e^4\wedge e^5+e^2\wedge e^3\wedge e^5-e^2\wedge e^4\wedge e^6)\nonumber.
\eea
As we only care about the sign of the eigenvalues of the principal symbol, we may assume that $|\xi|=2M=1$. Moreover by rotational symmetry we may assume that $\xi=e^1$. Under such reduction, it is easy to see that
\be
W=\{e^1\wedge\gamma:\gamma\in\Lambda^2(V')^*,\Lambda'\gamma=0\},\nonumber
\ee
where $V'=\textrm{span}\{e_j\}_{j=3}^6$ equipped with the symplectic form $\o'=e^3\wedge e^4+e^5\wedge e^6$ and $\Lambda'$ is the contraction with respect to $\o'$. In this way we can also simplify the operator to
\be
e^1\wedge\gamma\mapsto e^1\wedge\left[J\gamma+\frac{1}{2}(\gamma,e^3\wedge e^6+e^4\wedge e^5)(e^3\wedge e^6+e^4\wedge e^5)+(\gamma,e^{35}-e^{46})(e^{35}-e^{46})\right],\nonumber
\ee
which is equivalent to
\be
\gamma\mapsto J\gamma+\frac{1}{2}(\gamma,e^3\wedge e^6+e^4\wedge e^5)(e^3\wedge e^6+e^4\wedge e^5)+(\gamma,e^{35}-e^{46})(e^{35}-e^{46}).
\ee
Then it is clear that the eigenvalues are $\lambda=1$ (multiplicity 4) with eigenvectors $\gamma=e^3\wedge e^4-e^5\wedge e^6, e^3\wedge e^5+e^4\wedge e^6, e^3\wedge e^6-e^4\wedge e^5$ and $e^3\wedge e^5-e^4\wedge e^6$, $\lambda=0$ with eigenvector $\gamma=e^3\wedge e^6+e^4\wedge e^5$. We summarize our findings in the following lemma:

\begin{lemma}
\label{lm:eigenvalues}
The leading symbol in the Type IIA flow, restricted to closed and primitive forms, is only weakly parabolic. More precisely, it has an eigenvalue $\lambda=1$ with multiplicity $4$, and an eigenvalue $\lambda=0$ with multiplicity $1$.
\end{lemma}



\section{Proof of Theorem \ref{th:short}: existence}
\setcounter{equation}{0}

In this section we establish the short-time existence of the Type IIA flow. As we saw in \S 4, the flow is not strictly parabolic, and the presence of the symplectic form prevents a direct application of either the reparametrization arguments of \cite{DeT} or the Hamilton-Nash-Moser theorem of \cite{Ha1}. Rather, we proceed as follows: first we do apply a reparametrization, but we have to accompany it at the same time with a flow of the symplectic form. This new coupled flow of $(\varphi,\o)$ is still not strictly parabolic, but one of its key properties is that it admits a strictly regularization with integrability condition, to which the Hamilton-Nash-Moser theorem can apply. Thus we obtain the short-time existence for a regularized version of the Type IIA flow. Next, we show that the regularized flow preserves the primitiveness of the data, and reduces to the Type IIA flow if the form $\varphi$ is known to be primitive. Altogether, we obtain the desired short-time existence of the Type IIA flow for primitive data. The uniqueness of the solution will be shown later in \S \ref{7.5}, as a consequence of the uniqueness of the flows of the corresponding metrics.

\subsection{A coupled flow for $(\varphi,\o)$}

More precisely, we consider a reparametrization of the Type IIA flow by the following time-dependent vector field
\bea
V^k=e^u\left(g^{pq}(\Gamma^k_{pq}-(\Gamma_0)^k_{pq})-g^{lk}u_l\right),\label{vf}
\eea
where $|\varphi|$, $u$ and $g$ are defined by
\bea
|\varphi|^2\frac{\o^3}{3!}=\varphi\wedge\hat\varphi,\quad
u=\log|\varphi|^2,\quad
g_{ij}=-|\varphi|^{-2}\varphi_{iab}\varphi_{jcd}\o^{ac}\o^{bd}\nonumber,
\eea
and $\Gamma$ and $\Gamma_0$ are Christoffel symbols associated to the evolving metric $g$ and the initial metric $g_0$.
\medskip
Under a reparametrization by the diffeomorphisms generated by the vector field $V^k$, the given symplectic form in the Type IIA flow would
become time-dependent and evolve by its Lie derivative. It is convenient to change notation slightly, and denote the given symplectic form by $\o_0$ while reserving $\o=\o(t)$ for the evolving symplectic form. This consideration inspires us to consider the following coupled flow for the pair $(\varphi,\o)$,
\bea
\begin{cases}
&\p_t\varphi=d\Lambda d(|\varphi|^2\hat\varphi)+d(\iota_V\varphi)\\
&\p_t\o=d(\iota_V\o)
\end{cases}
\label{aux}
\eea
with initial data $\varphi(0)=\varphi_0$, $\o(0)=\o_0$, where $\varphi_0$ would be a closed primitive positive $3$-form with respect to $\o_0$. Although the initial metric $g_0$ is almost K\"ahler, a priori we should not assume that $\varphi(t)$ is primitive with respect to $\o(t)$, hence $g(t)$ a priori may not even be almost Hermitian.

Our first task is to work out the eigenvalues of the principal symbol for this coupled flow. Note that because of the coupling, the principal symbol of (\ref{aux}) is now a linear operator acting on both $\delta\varphi$ and $\delta\o$, and not just on $\delta\varphi$, and we may no longer assume that $\delta\varphi$ is primitive. It is easy to see that the principal symbol of the linearization of (\ref{aux}) is determined by
\bea
(\delta\varphi,
\delta\o)\ \rightarrow\ (d\Lambda d\,\delta(|\varphi|^2\hat\varphi)+d(\iota_{\delta V}\varphi),\
d(\iota_{\delta V}\o)).
\eea
Now the leading order term in $\delta V$ is
\bea
e^ug^{pq}g^{kl}(\nabla_p(\delta g)_{lq}-\frac{1}{2}\nabla_l(\delta g)_{pq})-g^{lk}\nabla_l\delta(e^u).\nonumber
\eea
so if we define the vector field $W_\xi$ by
\bea
W_\xi^k=|\varphi|^2g^{pq}g^{kl}(\xi_p(\delta g)_{lq}-\frac{1}{2}\xi_l(\delta g)_{pq})-g^{lk}\xi_l\delta|\varphi|^2,\nonumber
\eea
it follows immediately that the principal symbol of the linearized operator is
\bea
(\delta\varphi,\delta \o)\mapsto (\xi\wedge\Lambda(\xi\wedge\delta(|\varphi|^2\hat\varphi))+\xi\wedge\iota_{W_\xi}\varphi,\xi\wedge\iota_{W_\xi}\o)
\eea
with integrability conditions $\xi\wedge\delta\varphi=\xi\wedge\delta\o=0$.

We work out more explicitly the symbol at a point $(\varphi,\o)$ where $\varphi$ is primitive with respect to $\o$. In this case, we may choose an orthonormal basis of $g$ such that
\bea
\o&=&e^{12}+e^{34}+e^{56},\nonumber\\
\varphi&=&\frac{|\varphi|}{2}(e^{135}-e^{146}-e^{245}-e^{236}),\quad
\hat\varphi=\frac{|\varphi|}{2}(e^{136}+e^{145}+e^{235}-e^{246}).\nonumber
\eea
Without loss of generality, we may further assume that $\xi=e^1$ and $|\varphi|=1$. In this case, we can write $\delta\varphi=e^1\wedge\gamma$ and $\delta\o=e^1\wedge\alpha$ for some 2-form $\gamma$ and 1-form $\alpha$ such that $\alpha,\gamma\in\bigwedge^*\{e^2,\dots,e^6\}$. It is straightforward to check that
\bea
\delta|\varphi|^2&=&2(\delta\varphi,\varphi)-|\varphi|^2(\delta\o,\o)=(\gamma,e^{35}-e^{46})-(\alpha,e^2),\nonumber\\
\delta(|\varphi|^2\hat\varphi)&=&-|\varphi|^2J(\delta\varphi)-2(\delta\varphi,\hat\varphi)\varphi+4(\delta\varphi,\varphi)\hat\varphi -|\varphi|^2(\delta\o,\o)\hat\varphi\nonumber\\
&=&e^2\wedge J\gamma-(\gamma,e^{36}+e^{45})\varphi+(2(\gamma,e^{35}-e^{46})-(\alpha,e^2))\hat\varphi\nonumber\\
W^k_\xi&=&(\delta g)_{k1}-\delta^k_1\left(\frac{1}{2}\tr_g\delta g+(\gamma,e^{35}-e^{46})-(\alpha,e^2)\right).
\eea
We see that the key is to compute $\delta g$, especially $(\delta g)_{k1}$. By definition of $\tilde g$ and straightforward calculation, we have
\be
\begin{split}
&(\delta\tilde g)_{11}=2(\gamma,e^{35}-e^{46}),\quad (\delta\tilde g)_{22}=0,\quad (\delta\tilde g)_{12}=-(\gamma,e^{36}+e^{45}),\\
&(\delta\tilde g)_{33}=(\delta\tilde g)_{55}=2(\gamma,e^{35})-(\alpha,e^2),\quad (\delta\tilde g)_{44}=(\delta\tilde g)_{66}=-2(\gamma,e^{46})-(\alpha,e^2),\\
&(\delta\tilde g)_{13}=(\gamma,e^{26})+\frac{1}{2}(\alpha,e^4),\quad (\delta\tilde g)_{14}=(\gamma,e^{25})-\frac{1}{2}(\alpha,e^3),\\
&(\delta\tilde g)_{15}=-(\gamma,e^{24})+\frac{1}{2}(\alpha,e^6),\quad (\delta\tilde g)_{16}=-(\gamma,e^{23})-\frac{1}{2}(\alpha,e^5),\nonumber\\
&(\delta g)_{ij}=(\delta\tilde g)_{ij}-((\gamma,e^{35}-e^{46})-(\alpha,e^2))\delta_{ij}.\nonumber
\end{split}
\ee
It follows that $\tr_g\delta g=2(\alpha,e^2)$, hence
\bea
&W_\xi^1=(\alpha,e^2),&
W_\xi^2=-(\gamma,e^{36}+e^{45})\nonumber\\
&W_\xi^3=(\gamma,e^{26})+\dfrac{1}{2}(\alpha,e^4),&
W_\xi^4=(\gamma,e^{25})-\dfrac{1}{2}(\alpha,e^3),\nonumber\\
&W_\xi^5=-(\gamma,e^{24})+\dfrac{1}{2}(\alpha,e^6),&
W_\xi^6=-(\gamma,e^{23})-\dfrac{1}{2}(\alpha,e^5).\nonumber
\eea
Consequently we find that
\bea
\xi\wedge\iota_{W_\xi}\varphi&=&\frac{1}{2}e^1\wedge\bigg[(\gamma_{36}+\gamma_{45})(e^{36}+e^{45})+(\gamma_{26}+\frac{\alpha_4}{2})e^{26}+(\gamma_{25} -\frac{\alpha_3}{2})e^{25}\nonumber\\
&&+(\gamma_{24}-\frac{\alpha_6}{2})e^{24}+(\gamma_{23}+\frac{\alpha_5}{2})e^{23}+\alpha_2(e^{35}-e^{46})\bigg],\nonumber\\
\xi\wedge\iota_{W_\xi}\o&=&\frac{1}{2}e^1\wedge(\alpha+2\alpha_2e^2 -2(\gamma_{25}e^3-\gamma_{26}e^4-\gamma_{23}e^5+\gamma_{24}e^6)).\nonumber
\eea
If we further write $\gamma=e^2\wedge\beta+\lambda$, $\alpha=\alpha_2e^2+\mu$, where $\beta,\mu,\lambda\in\bigwedge^*\{e^3,\dots,e^6\}$, we have
\bea
\xi\wedge\iota_{W_\xi}\varphi&=&\frac{1}{2}e^1\wedge\bigg[(\lambda_{36}+\lambda_{45})(e^{36}+e^{45})+e^2\wedge\beta -\frac{1}{2}e^2\wedge\iota_\mu(e^{35}-e^{46})\nonumber\\
&&+\alpha_2(e^{35}-e^{46})\bigg],\nonumber\\
\xi\wedge\iota_{W_\xi}\o&=&\frac{1}{2}e^1\wedge(\mu+2\alpha_2e^2+2\iota_\beta(e^{35}-e^{46})),\nonumber\\
\xi\wedge\Lambda(\xi\wedge\delta(|\varphi|^2\hat\varphi)&=&e^1\wedge\bigg[J\lambda+\frac{\lambda_{36}+\lambda_{45}}{2}(e^{36}+e^{45}) +(\lambda_{35}-\lambda_{46}-\frac{\alpha_2}{2})(e^{35}-e^{46})\bigg],\nonumber
\eea
It follows that the principal symbol is the linear map
\bea
(\beta,\lambda,\alpha_2,\mu)\mapsto\left(\frac{\beta}{2}-\frac{\iota_\mu}{4}(e^{35}-e^{46}),\lambda,\alpha_2,\frac{\mu}{2}+\iota_\beta(e^{35}-e^{46})\right)
\eea
This matrix is only positive semi-definite. The part $(\lambda,\alpha_2)\mapsto(\lambda,\alpha_2)$ is the identity map. However the other part
\bea
(\beta,\mu)\mapsto \left(\dfrac{\beta}{2}-\dfrac{\iota_\mu}{4}(e^{35}-e^{46}),\dfrac{\mu}{2}+\iota_\beta(e^{35}-e^{46})\right)\nonumber
\eea
has eigenvalues 0 and 1, both of multiplicities 4. So the coupled flow (\ref{aux}) for $(\varphi,\o)$ is still not strictly parabolic.

\subsection{A parabolic regularization of the coupled flow}

To solve this problem, we add an extra term on the right hand side of the evolution equation of $\varphi$ in (\ref{aux}). This term takes the form
\bea
-BdJd(|\varphi|^2\Lambda\hat\varphi),\nonumber
\eea
where $B$ is a constant to be determined. In fact, $\Lambda(\hat\varphi)$ is expected to be zero along the flow as $\varphi$ should always be primitive. Again let us consider the linearization of $|\varphi|^2\Lambda\hat\varphi$ at a primitive pair $(\o,\varphi)$ and we may assume that $|\varphi|^2=1$ at the point of linearization. It follows that
\bea
\delta(|\varphi|^2\Lambda\hat\varphi)=(\delta\Lambda)(\hat\varphi)-\Lambda(J\delta\varphi).
\eea
As before, we may assume that $\xi=e^1$ and $\delta\varphi=e^1\wedge(e^2\wedge\beta+\lambda)$ and $\delta\omega=e^1\wedge(\alpha_2e^2+\mu)$. The principal symbol for the extra term is
\bea
Be^{12}\wedge J\big[(\delta\Lambda)(\hat\varphi)-\Lambda(J\delta\varphi\big)].
\eea
The second term is easy to compute:
\bea
-Be^{12}\wedge J(\Lambda(J\delta\varphi))=Be^{12}\wedge\beta.
\eea
The first term is more complicated, notice that
\bea
B(\delta\Lambda)(\hat\varphi)_k=-\frac{B}{2}\o^{js}(\delta\o)_{st}\o^{ti}\hat\varphi_{ijk}=B\sum_{t=3}^6\mu_t\o^{ti}\hat\varphi_{2ik},\nonumber
\eea
therefore by straightforward calculation, this part of the principal symbol is equivalent to the linear map
\bea
\mu\mapsto \frac{B}{2}\iota_\mu(e^{35}-e^{46}).\nonumber
\eea
Therefore the principal symbol for the full evolution equation is equivalent to the linear map
\bea
(\beta,\lambda,\alpha_2,\mu)\mapsto\left(\frac{\beta}{2}(1+2B)+\frac{2B-1}{4}\iota_\mu(e^{35}-e^{46}), \lambda,\alpha_2,\frac{\mu}{2}+\iota_\beta(e^{35}-e^{46})\right).\nonumber
\eea
If $B>0$, then all the eigenvalues of the above matrix are positive. In this sense, the coupled flow with the additional $B$ term is parabolic.

\begin{lemma}
Consider the flow
\bea
\label{aux;B}
&&
\p_t\varphi=d\Lambda d(|\varphi|^2\hat\varphi)-BdJd(|\varphi|^2\Lambda\hat\varphi)+d(\iota_V\varphi),\nonumber\\
&&\p_t\o=d(\iota_V\o),
\eea
for any fixed, strictly positive constant $B$. Then for any initial value $\varphi_0$ which is a closed, positive, and primitive form with respect to the initial symplectic form $\o_0$, the flow exists and is smooth at least on some interval $[0,T)$ with $T>0$. Clearly the flow preserves both the closedness of both the form $\varphi$ and the symplectic form $\o$.
\end{lemma}

\noindent
{\it Proof.} Let $d$ be the exterior derivative. The preceding fact that the eigenvalues of the principal symbol of the flow (\ref{aux;B}) when restricted to closed and primitive forms are positive means that the flow (\ref{aux;B}) together with $d$ as the integrability operator satisfies the condition of the Hamilton-Nash-Moser theorem (\cite{Ha1}, Theorem 5.1). This theorem implies the short-time existence and uniqueness of the flow (\ref{aux;B}). Q.E.D.

\medskip

It should be noted that we need to treat $\varphi$ and $\o$ as tensors evolving independently at this moment, therefore we cannot assume that $\varphi$ is primitive with respect to $\o$ (though we shall prove it is indeed the case later). Consequently the metric $g$ defined above is not necessarily compatible with $J$ or $\omega$: we only know it is a Riemannian metric. As (\ref{aux;B}) preserves the closedness of $\varphi$ and $\o$, by performing the reverse reparametrization, we obtain immediately

\begin{lemma}
Fix any fixed positive constant $B$. Then the flow of $3$-forms $\varphi$
\bea
\label{aux;BB}
\p_t\varphi=d\Lambda d(|\varphi|^2\hat\varphi)-BdJd(|\varphi|^2\Lambda\hat\varphi)
\eea
admits a closed smooth solution $\varphi$ on some interval $[0,T)$ with $T>0$, for any initial value $\varphi_0$ which is a smooth closed, positive, and primitive form with respect to the symplectic form $\o_0$.
\end{lemma}

\bigskip

\subsection{Preservation of the primitiveness condition}

Next we shall show that, if the initial data $\varphi_0$ is primitive in the flow (\ref{aux;BB}), then $\varphi(t)$ remains primitive for all time. Since $\varphi$ primitive implies that $\hat\varphi$ is also primitive, it follows that the terms with coefficient $B$ in (\ref{aux;BB}) all drop out, and the flow reduces to the Type IIA flow, establishing Theorem \ref{th:short} in the case of no sources.

\medskip
From now on, we take $B=1$. Let $\varphi(t)$ be a solution to (\ref{aux;BB}) on $M\times[0,T)$ with $\varphi(0)$ being closed, positive, and primitive. Clearly for any $t$, $\varphi(t)$ stays closed. Let
\bea
\varphi=P+\beta\wedge\omega
\eea
be the primitive decomposition of $\varphi$, where $P$ is a primitive 3-form. It follows that
$$
\beta=\dfrac{\Lambda \varphi}{2}.
$$
We wish to show that $\beta=0$ by the maximum principle. To do so, we need to compute the evolution equation of $\beta$. We fix a background metric $\bar g=g(0)$ which is compatible with $\o$. We denote by $\bar\na$ the covariant derivatives with respect to $\bar g$. Since $\varphi$ is closed, $d\beta$ is primitive, and thus
\bea
\o^{jk}\bar\nabla_j\beta_k=0.
\eea
Furthermore,
\bea
(\Lambda\hat\varphi)_k=\frac{\o^{ji}}{2}\hat\varphi_{ijk}=-\frac{\o^{ji}}{2}\varphi_{i,j,Jk}=-(\Lambda\varphi)_{Jk},\nonumber
\eea
hence there exists a primitive 3-form $\hat P$ such that the primitive decomposition for $\hat\varphi$ is
\bea
\hat\varphi=\hat P-J\beta\wedge\o.
\eea
Using this decomposition for $\hat\varphi$,  we can derive the evolution equation for $\beta$,
\bea
\label{fl:beta}
\p_t\beta=-d\Lambda d(|\varphi|^2J\beta)+\Lambda(dJd(|\varphi|^2J\beta)).\label{beta}
\eea
We know that $\varphi$, $J$, and all their covariant derivatives are bounded in $M\times[0,\tau]$ for any $\tau<T$, therefore we can write (\ref{beta}) in the form
\bea
\p_t\beta=|\varphi|^2(-d\Lambda d(J\beta)+\Lambda(dJdJ\beta)+\bar\nabla\beta*S_1+\beta*S_2),\label{beta2}
\eea
where $S_1$ and $S_2$ are bounded tensors on $M\times[0,\tau]$ and $*$ represents certain contraction of indices. We need to compute the leading term of $\beta$ in (\ref{beta2}).
Notice that $(J\beta)_j=J^p{}_j\beta_p$, so
\bea
d(J\beta)_{jk}&=&J^p{}_k\bar\nabla_j\beta_p-J^p{}_j\bar\nabla_k\beta_p+O(\beta),\nonumber\\
\Lambda d(J\beta)&=&\o^{kj}J^p{}_k\bar\nabla_j\beta_p+O(\beta),\nonumber\\
(d\Lambda d(J\beta))_l&=&\o^{Jp,j}\bar\nabla_l\bar\nabla_j\beta_p+O(\beta,\bar\nabla\beta),\nonumber\\
(JdJ\beta)_{jk}&=&J^s{}_jJ^t{}_k(J^p{}_t\bar\nabla_s\beta_p-J^p{}_s\bar\nabla_t\beta_p+O(\beta))\nonumber\\
&=&J^t{}_k\bar\nabla_t\beta_j-J^t{}_j\bar\nabla_t\beta_k+O(\beta),\nonumber\\
(dJdJ\beta)_{jkl}&=&J^t{}_k(\bar\nabla_l\bar\nabla_t\beta_j-\bar\nabla_j\bar\nabla_t\beta_l)-J^t{}_j(\bar\nabla_l\bar\nabla_t\beta_k -\bar\nabla_k\bar\nabla_t\beta_l)\nonumber\\
&&+J^t{}_l(\bar\nabla_j\bar\nabla_t\beta_k-\bar\nabla_k\bar\nabla_t\beta_j)+O(\beta,\bar\nabla\beta),\nonumber\\
\Lambda(dJdJ\beta)_l&=&\o^{kj}J^t{}_k(\bar\nabla_l\bar\nabla_t\beta_j-\bar\nabla_j\bar\nabla_t\beta_l)+J^t{}_l\bar\nabla_t(\o^{kj}\bar\nabla_j\beta_k) +O(\beta,\bar\nabla\beta)\nonumber\\
&=&\o^{Jj,p}\bar\nabla_l\bar\nabla_j\beta_p+\o^{j,Jt}\bar\nabla_j\bar\nabla_t\beta_l+O(\beta,\bar\nabla\beta).\nonumber
\eea
It follows that
\bea
\p_t\beta=|\varphi|^2(L\beta+O(\beta,\bar\nabla\beta)),
\eea
where $L$ is defined by
\bea
(L\beta)_l=(\o^{Jj,p}+\o^{j,Jp})\bar\nabla_l\bar\nabla_j\beta_p+\o^{j,Jt}\bar\nabla_j\bar\nabla_t\beta_l.\nonumber
\eea
We further notice that
$\o^{Jj,p}+\o^{j,Jp}=O(\beta)$, $\o^{j,Jt}-\delta^{jt}=O(\beta)$, and as $\beta$ is a smooth function of $\varphi$, we also have $|\bar\nabla^2\beta|$ is uniformly bounded. Therefore one can also write
\bea
\p_t\beta=|\varphi|^2(\bar\Delta\beta+O(\beta,\bar\nabla\beta)).
\eea
It follows that
\bea
\p_t|\beta|^2_{\bar g}&=&2|\varphi|^2(\bar\Delta\beta,\beta)_{\bar g}+O(\beta,\bar\nabla\beta)*\beta*S\nonumber\\
&\leq& |\varphi|^2\bar\Delta(|\beta|^2_{\bar g})-2|\varphi|^2|\bar\nabla\beta|^2+O(\beta,\bar\nabla\beta)*\beta*S\nonumber\\
&\leq& |\varphi|^2\bar\Delta(|\beta|^2_{\bar g})+C|\beta|^2.
\eea
Since $\beta=0$ initially, by maximum principle we know that $\beta=0$ on $[0,\tau]$ for any $\tau<T$. Therefore $\varphi(t)$ is primitive for as long as the flow exists, and $J$ is always compatible with $\o$. The existence part of Theorem \ref{th:short} is proved in the case of no sources.

\section{Type IIA geometry: proof of Theorem \ref{th:holonomy}\label{ahg}}
\setcounter{equation}{0}
\label{s:geometry}

The goal of this section is to work out some properties specific to Type IIA geometry. What is crucial is that the almost complex structure $J_\varphi$ in Type IIA strings comes from a closed primitive positive 3-form $\varphi$ via Hitchin's construction. In fact, the closedness of $\varphi$ imposes subtle ``higher integrability'' conditions on $J_\varphi$ which in turn distinguish $J_\varphi$ from a generic almost complex structure. This feature gives rise to various identities that are not available in the more general almost-K\"ahler setting.
We begin with the curvature and Nijenhuis tensor on general almost-complex manifolds, and gradually specialize to almost-K\"ahler manifolds, and then to Type IIA geometry.

\subsection{Curvature tensors on general almost-complex manifolds}

For any affine connection $D$, we define its curvature tensor $R(D)$ and torsion $T(D)$ by
\bea
[D_i,D_j]X^m=R(D)_{ij}{}^m{}_lX^l-T(D)^l{}_{ij}D_lX^m.
\eea
The curvature tensor with four lower indices is defined in the usual way
\bea
R(D)_{ijkl}=R(D)_{ij}{}^p{}_lg_{pk}.
\eea
As in the case for Levi-Civita connection, we define the Ricci curvature of $D$, also denoted by $R(D)$, by
\bea
R(D)_{ik}=g^{jl}R(D)_{ijkl}.
\eea

Let now $J$ be an almost-complex structure on the Riemannian manifold $M$.
In subsequent developments, we shall need both the Levi-Civita connection $\nabla$ and the projected Levi-Civita connection $\mathfrak D=\mathfrak D^0$. Therefore we will reserve the Latin letter $R$ for various curvature tensors associated to $\nabla$ and the German letters $\mathfrak\Gamma$, $\mathfrak T$, and $\mathfrak R$ for Christoffel symbol, torsion, and curvature tensors associated to $\mathfrak D$. When we have other Hermitian metrics with decoration like $\tilde g$ or $\hat g$, we shall decorate the corresponding connections and curvature tensors with the same symbol. Identities for the curvature and torsion of the Chern connection have been worked out in the paper of Tosatti, Weinkove, and Yau \cite{TWY}. However, they are expressed there in complex frames, and it is difficult for us to apply their formulas, as we shall have to let the almost-complex structure evolve. Thus we develop here a formalism for curvature and torsion identities with the action of $J$ in real coordinate systems.

\smallskip

To pass back and forth from $\na$ to $\mathfrak\Gamma$, we note that from (\ref{Gamma}) and (\ref{T1}) that
\bea
\mathfrak \Gamma^m{}_{ij}&=&\Gamma^m{}_{ij}-N_{ij}{}^m-V_{ij}{}^m=:\Gamma^m{}_{ij}-A_{ij}{}^m,\label{christ}\\
\mathfrak T^m{}_{ij}&=&N^m{}_{ij}-U^m{}_{ij},\label{torsion}
\eea
where $A=V+N$ is of type $(2,0)+(0,2)$, and hence their curvature tensors are related by
\bea
R_{ijkl}&=&\mathfrak R_{ijkl}-(\mathfrak D_iA_{jkl}-\mathfrak D_jA_{ikl}+\mathfrak T^p{}_{ij}A_{pkl}+A_{ik}{}^pA_{jlp}-A_{jk}{}^pA_{ilp}).\label{curvrel}
\eea
As $R_{ijkl}$ is the curvature tensor of the Levi-Civita connection, it has various symmetries and satisfies the Bianchi identities. On the other hand, since $\mathfrak DJ=0$, its curvature $\mathfrak D$ satisfies
\bea
\mathfrak R_{ijkl}=\mathfrak R_{i,j,Jk,Jl}.
\eea
It is easy to deduce from the preceding relation between $R_{ijkl}$ and $\mathfrak R_{ijkl}$ how to modify the identity for each curvature if it is replaced by the other.

\medskip

The projected Levi-Civita connection $\mathfrak D$ induces a connection on the canonical bundle of $M$, whose curvature represents the first Chern class (up to a constant) of the almost complex manifold $(M,J)$. To be precise, if we use small Greek letter to denote the index for ``holomorphic'' tangent bundle $T^{1,0}M$, then
\bea
\frac{1}{4\pi}\mathfrak R_{ijkl}\o^{lk}=\frac{\sqrt{-1}}{2\pi}\mathfrak R_{ij}{}^\gamma{}_\gamma\in [c_1(M,J)],
\eea
Since $A$ is of type $(2,0)+(0,2)$, the contraction of its last two indices using $\o$ or $g$ vanishes, therefore by (\ref{curvrel}) we see that
\bea
\mathfrak R_{ijkl}\o^{lk}=R_{ijkl}\o^{lk}+2A_{ik}{}^pA_{jlp}\o^{lk}
\eea
is a closed 2-form. We end this subsection by deriving the following formula for Ricci curvature
\bea
\label{ricciahg}
R_{ij}=-2g^{kl}(\mathfrak D_iA_{kjl}-\mathfrak D_kA_{ijl}+\mathfrak T^p{}_{ik}A_{pjl})+\frac{1}{2}\o^{lk}R_{i,Jj,k,l}.
\eea
Indeed, by (\ref{curvrel}) we see that
\bea
R_{i,j,Jk,Jl}-R_{ijkl}&=&(\mathfrak D_iA_{jkl}-\mathfrak D_jA_{ikl}+\mathfrak T^p{}_{ij}A_{pkl}+A_{ik}{}^pA_{jlp}-A_{jk}{}^pA_{ilp})\nonumber\\
&&-(\mathfrak D_iA_{j,Jk,Jl}-\mathfrak D_jA_{i,Jk,Jl}+\mathfrak T^p{}_{ij}A_{p,Jk,Jl}+A_{i,Jk}{}^pA_{j,Jl,p}-A_{j,Jk}{}^pA_{i,Jl,p}).\nonumber
\eea
Recall that $A$ is of type $(2,0)+(0,2)$, so $A_{ijk}=-A_{i,Jj,Jk}$, therefore we get
\bea \label{Rij-JkJl}
R_{i,j,Jk,Jl}-R_{ijkl}=2(\mathfrak D_iA_{jkl}-\mathfrak D_jA_{ikl}+\mathfrak T^p{}_{ij}A_{pkl}).
\eea
Let us denote the right hand side of the above equation by $B_{ijkl}$. Then the above equation is equivalent to
\bea
-R_{i,j,Jk,l}-R_{i,j,k,Jl}&=&2(\mathfrak D_iA_{j,k,Jl}-\mathfrak D_jA_{i,k,Jl}+\mathfrak T^p{}_{ij}A_{p,k,Jl})\nonumber\\
&=&B_{i,j,Jk,l}=B_{i,j,k,Jl}.\label{rel2}
\eea
Let $\{e_a\}$ be an orthonormal frame for the given Riemannian metric, and so is the frame $\{Je_a\}$. By definition of Ricci curvature, we have
\bea
R_{i,Jj}&=&\sum_a R(i,e_a,Jj,e_a)=\sum_a R(i,Je_a,Jj,Je_a)\nonumber\\
&=&\sum_a (R(i,Je_a,j,e_a)+B(i,Je_a,j,e_a))\nonumber\\
&=&\sum_a(-R(j,e_a,J(Ji),Je_a)-B(i,Je_a,Jj,Je_a))\nonumber\\
&=&\sum_a(-R(j,e_a,Ji,e_a)-B(j,e_a,Ji,e_a)-B(i,e_a,Jj,e_a))\nonumber\\
&=&-R_{j,Ji}-g^{kl}(B_{i,k,Jj,l}+B_{j,k,Ji,l}).\label{imp1}
\eea
On the other hand, by taking trace of (\ref{rel2}) and using Bianchi identity of $R$, we have
\bea
g^{jl}B_{i,j,Jk,l}&=&-g^{jl}(R_{i,j,Jk,l}+R_{i,j,k,Jl})\nonumber\\
&=&-R_{i,Jk}+g^{jl}(R_{j,k,i,Jl}+R_{k,i,j,Jl})\nonumber\\
&=&-R_{i,Jk}+g^{jl}(R_{k,j,Ji,l}+B_{k,j,Ji,l})+R_{kijl}\o^{lj}\nonumber\\
&=&-R_{i,Jk}+R_{k,Ji}+g^{jl}B_{k,j,Ji,l}+R_{kijl}\o^{lj}.\label{imp2}
\eea
(\ref{imp1}) and (\ref{imp2}) can be rewritten as
\bea
R_{i,Jj}+R_{j,Ji}&=&-g^{kl}(B_{i,k,Jj,l}+B_{j,k,Ji,l}),\nonumber\\
R_{i,Jj}-R_{j,Ji}&=&-g^{kl}(B_{i,k,Jj,l}-B_{j,k,Ji,l})-R_{ijkl}\o^{lk}.\nonumber
\eea
Adding these two equations up we get
\bea
R_{i,Jj}=-g^{kl}B_{i,k,Jj,l}-\frac{1}{2}R_{ijkl}\o^{lk}\nonumber
\eea
which is equivalent to
\bea
R_{ij}&=&-g^{kl}B_{ikjl}+\frac{1}{2}\o^{lk}R_{i,Jj,k,l}\nonumber\\
&=&-2g^{kl}(\mathfrak D_iA_{kjl}-\mathfrak D_kA_{ijl}+\mathfrak T^p{}_{ik}A_{pjl})+\frac{1}{2}\o^{lk}R_{i,Jj,k,l}.\nonumber
\eea
This gives the desired formula.

\subsection{Quadratic expressions in the Nijenhuis tensor}

We shall encounter frequently later quadratic expressions of the Nijenhuis tensor. It is convenient to introduce the following two symmetric tensors quadratic in $N$:
\bea
(N^2_+)_{ij}&:=&N^{pk}{}_iN_{pkj}\geq0,\nonumber\\
(N^2_-)_{ij}&:=&N^{kp}{}_iN_{pkj}.\nonumber
\eea
Since $N$ is skew-symmetric in the last two slots and it satisfies the Bianchi identity (\ref{nbianchi}), all the other similar tensors can be expressed as a linear combination of $N^2_+$ and $N^2_-$. For example
\be
0\leq N_{ipk}N_j{}^{pk}=(N_{pki}-N_{kpi})(N^{pk}{}_j-N^{kp}{}_j)=2(N^2_+)_{ij}-2(N^2_-)_{ij}.\nonumber
\ee
Obviously $g^{ij}N_{ipk}N_j{}^{pk}=|N|^2=g^{ij}N^{pk}{}_iN_{pkj}$, so we find that
\be
|N|^2=\tr N^2_+=g^{ij}(N^2_+)_{ij}=2g^{ij}(N^2_-)_{ij}=2\tr N^2_-.
\ee
Also we observe that both $N^2_+$ and $N^2_-$ are $J$-invariant in the sense that $(N^2_{\pm})_{ij}=(N^2_{\pm})_{Ji,Jj}$. In general, for any symmetric 2-tensor $A=A_{ij}$, we define its $J$-invariant and $J$-anti-invariant parts respectively by
\bea
(A^J)_{ij}:=\frac{1}{2}(A_{ij}+A_{Ji,Jj}),\qquad
(A^{-J})_{ij}:=\frac{1}{2}(A_{ij}-A_{Ji,Jj}).\nonumber
\eea
In this notation, we have $N^2_\pm=(N^2_\pm)^J$.

Clearly we have $A=A^J+A^{-J}$ and this decomposition is orthogonal with respect to the inner product induced by the metric $g$. Later such a decomposition will play an important role in our calculations.

\subsection{Curvature tensors in almost-K\"ahler geometry\label{akg}}

In this subsection we restrict ourselves to the case $d\o=0$, namely the case $(M,J,g)$ is an almost-K\"ahler manifold. Since $\o$ is a symplectic form, we know $d^c\o=0$, hence both $U$ and $V$ defined in (\ref{UV}) are zero. Therefore (\ref{christ}) and (\ref{torsion}) specialize to $A=N=\mathfrak T$. Therefore the previously deduced formula (\ref{curvrel}) becomes
\bea
R_{ijkl}=\mathfrak R_{ijkl}-(\mathfrak D_iN_{jkl}-\mathfrak D_jN_{ikl}+N^p{}_{ij}N_{pkl}+N_{ik}{}^pN_{jlp}-N_{jk}{}^pN_{ilp}),\label{curvrel2}
\eea
thus we have
\bea
R_{ij}=\mathfrak R_{ij}+\mathfrak D^kN_{ijk}-(N^2_-)_{ij}.
\eea
Combining (\ref{ricciahg}) with (\ref{curvrel2}), we also obtain
\bea
R_{ij}&=&2\mathfrak D^kN_{ijk}-2(N^2_+)_{ij}+\frac{1}{2}\o^{lk}R_{i,Jj,k,l}\nonumber\\
&=&2\mathfrak D^kN_{ijk}-2(N^2_+)_{ij}\nonumber\\
&&+\frac{1}{2}\o^{lk}(\mathfrak R_{i,Jj,k,l}-(\mathfrak D_iN_{Jj,k,l}-\mathfrak D_{Jj}N_{ikl}+N^p{}_{i,Jj}N_{pkl}+N_{ik}{}^pN_{Jj,l,p}-N_{Jj,k}{}^pN_{ilp}))\nonumber\\
&=&2\mathfrak D^kN_{ijk}-2(N^2_-)_{ij}+\frac{1}{2}\o^{lk}\mathfrak R_{i,Jj,k,l}.\label{ricciakg1}
\eea
Alternatively
\bea
R_{ij}&=&2\mathfrak D^kN_{ijk}-2(N^2_+)_{ij}+\frac{1}{2}\o^{lk}R_{k,l,i,Jj}\nonumber\\
&=&2\mathfrak D^kN_{ijk}-2(N^2_+)_{ij}\nonumber\\
&&+\frac{1}{2}\o^{lk}(\mathfrak R_{k,l,i,Jj}-(\mathfrak D_kN_{l,i,Jj}+N^p{}_{kl}N_{p,i,Jj}-\mathfrak D_lN_{k,i,Jj}+N_{ki}{}^pN_{l,Jj,p}-N_{li}{}^pN_{k,Jj,p}))\nonumber\\
&=&\mathfrak D^k(N_{ijk}+N_{jik})-(N^2_+)_{ij}+\frac{1}{2}\o^{lk}\mathfrak R_{k,l,i,Jj}.\label{ricciakg2}
\eea
From (\ref{ricciakg2}) we can immediately read off that
\bea
(R^J)_{ij}&=&-(N^2_+)_{ij}+\frac{1}{2}\o^{lk}\mathfrak R_{k,l,i,Jj},\label{jinvricci}\\
(R^{-J})_{ij}&=&\mathfrak D^k(N_{ijk}+N_{jik}).\label{-jinvricci}
\eea
Taking the trace of (\ref{ricciakg1}) and plugging in (\ref{curvrel2}), we see that
\bea
R=\frac{1}{2}\o^{ji}\o^{lk}\mathfrak R_{ijkl}-|N|^2=\frac{1}{2}\o^{ji}\o^{lk}R_{ijkl}-2|N|^2.\label{scalarcurv}
\eea
In the literature, the expression $\dfrac{1}{2}\o^{ji}\o^{lk}R_{ijkl}$ is sometimes known as the $\star$-scalar curvature. This relation (\ref{scalarcurv}) was first discovered by Blair-Ianus \cite{BI}, and Blair \cite{Blair} together with Oproiu \cite{O}.

Combining (\ref{curvrel2}) with the symmetry of $R$, we can derive that
\bea
\mathfrak R_{ijkl}-\mathfrak R_{klij}&=&{\mathfrak D}_iN_{jkl}-{\mathfrak D}_jN_{ikl}-{\mathfrak D}_kN_{lij}+{\mathfrak D}_lN_{kij}\nonumber\\ &&+N_{ik}^{~~p}N_{jlp}-N_{jk}^{~~p}N_{ilp}-N_{ki}^{~~p}N_{ljp}+N_{li}^{~~p}N_{kjp}.\label{frakrcomm}
\eea
As $\mathfrak R_{ijkl}=\mathfrak R_{i,j,Jk,Jl}$, by making use of (\ref{frakrcomm}), we get
\bea
&&\mathfrak R_{Ji,Jj,k,l}-\mathfrak R_{ijkl}=\mathfrak R_{Ji,Jj,k,l}-\mathfrak R_{k,l,Ji,Jj}+\mathfrak R_{klij}-\mathfrak R_{ijkl}\nonumber\\
&=&{\mathfrak D}_{Ji}N_{Jj,k,l}-{\mathfrak D}_{Jj}N_{Ji,k,l}-{\mathfrak D}_iN_{jkl}+{\mathfrak D}_jN_{ikl}+2{\mathfrak D}_kN_{lij}-2{\mathfrak D}_lN_{kij}.\label{prep}
\eea
Notice that the LHS of (\ref{prep}) does not change if one replace $k$ and $l$ by $Jk$ and $Jl$ respectively, so we get an interesting identity satisfied by $\mathfrak DN$
\bea
&&{\mathfrak D}_{Ji}N_{Jj,k,l}-{\mathfrak D}_{Jj}N_{Ji,k,l}-{\mathfrak D}_{Jk}N_{Jl,i,j}+{\mathfrak D}_{Jl}N_{Jk,i,j}\nonumber\\
&=&{\mathfrak D}_iN_{jkl}-{\mathfrak D}_jN_{ikl}-{\mathfrak D}_kN_{lij}+{\mathfrak D}_lN_{kij},\label{sym}
\eea
which allows us to rewrite one partial derivative of $N$ in terms of some other combination of partial derivatives.

In the same vein we can derive the Bianchi-type identity for $\mathfrak R$
\bea
&&{\mathfrak R}_{ijkl}+{\mathfrak R}_{jkil}+{\mathfrak R}_{kijl}\nonumber\\
&=&-{\mathfrak D}_iN_{ljk}-{\mathfrak D}_jN_{lki}-{\mathfrak D}_kN_{lij}+N^p{}_{ij}N_{lkp}+N^p{}_{jk}N_{lip}+N^p{}_{ki}N_{ljp}.\label{bianchi}
\eea

Equation (\ref{prep}) accounts for the (2,0)+(0,2)-part of the curvature tensor $\mathfrak R$. If we use Greek letters for barred and unbarred directions, then (\ref{prep}) can be translated into
\bea
{\mathfrak R}_{\alpha\beta\bar\gamma\delta}&=&{\mathfrak D}_\alpha N_{\beta\bar\gamma\delta}-{\mathfrak D}_\beta N_{\alpha\bar\gamma\delta}-{\mathfrak D}_{\bar\gamma} N_{\delta\alpha\beta}+{\mathfrak D}_{\delta} N_{\bar\gamma\alpha\beta}\nonumber\\
&=&-{\mathfrak D}_{\bar\gamma}N_{\delta\alpha\beta},\nonumber\\
{\mathfrak R}_{**\gamma\delta}&=&{\mathfrak R}_{**\bar\gamma\bar\delta}=0,\nonumber
\eea
which is the content of (2.17) in \cite{TWY}. Replace $i,j,k,l$ by $\bar\alpha,\beta,\bar\gamma,\delta$ respectively in the Bianchi-type identity (\ref{bianchi}), we get
\bea
{\mathfrak R}_{\bar\alpha\beta\bar\gamma\delta}-{\mathfrak R}_{\bar\gamma\beta\bar\alpha\delta}={\mathfrak R}_{\bar\alpha\beta\bar\gamma\delta}+{\mathfrak R}_{\beta\bar\gamma\bar\alpha\delta}+{\mathfrak R}_{\bar\gamma\bar\alpha\beta\delta}
=N^p{}_{\bar\gamma\bar\alpha}N_{\delta\beta p}\nonumber=N^\lambda{}_{\bar\gamma\bar\alpha}N_{\delta\beta\lambda}.\nonumber
\eea
This is the content of (2.16) in \cite{TWY}.

\subsection{The holonomy of Type IIA geometry}

We now restrict ourselves further to Type IIA geometry, namely a triple $(M,\o,\varphi)$ where $(M,\o)$ is a symplectic 6-manifold and $\varphi$ is a closed positive $\o$-primitive 3-form.

\smallskip
Our first task is to prove Theorem \ref{th:holonomy}(a). Recall that $|\varphi|$ is the norm of $\varphi$ with respect to the metric $g_\varphi$, and that we have defined the metric $\tilde g_\varphi$ by $\tilde g_\varphi=|\varphi|^2g_\varphi$.
It is not hard to see that
\bea
|\varphi|_{\tilde g}=|\varphi|^{-2}.\nonumber
\eea
Henceforth we shall denote $g_\varphi$, $\tilde g_\varphi$, and $J_\varphi$ by just $g$, $\tilde g$, and $J$ for simplicity. It is clear that $J$ is compatible with $\tilde g$ and the corresponding K\"ahler form $\tilde\o=|\varphi|^2\o$ satisfies
\bea
d\tilde\o=-\alpha\wedge\tilde\o,\qquad
d^c\tilde\o=J\alpha\wedge\tilde\o,\nonumber
\eea
where
\be
\alpha=d\log|\varphi|_{\tilde g}=-d\log|\varphi|^2.\label{alpha}
\ee
It follows from (\ref{vexp}) that
\bea
U_{ijk}&=&\frac{1}{4}(2\alpha_{Ji}\tilde\o_{jk}+\alpha_{Jj}\tilde\o_{ki}+\alpha_{Jk}\tilde\o_{ij}+\alpha_j\tilde g_{ki}-\alpha_k\tilde g_{ij}),\label{uexp2}\\
V_{ijk}&=&\frac{1}{4}(\alpha_{Jj}\tilde\o_{ki}+\alpha_{Jk}\tilde\o_{ij}-\alpha_j\tilde g_{ki}+\alpha_k\tilde g_{ij}).\label{vexp2}
\eea


We need now the following lemmas for computational purposes.

\begin{lemma}\label{torsiond}~\\
Let $\mu$ be any differential form, $D$ any affine connection, $T=T(D)$ the torsion tensor associated to $D$. Then we have the following formula
\begin{equation}
d\mu=d x^j\wedge D_j\mu+T\boxtimes\mu,\label{dconnection}
\end{equation}
where $\boxtimes$ is a multiplication operation linear in both factors. We only need the explicit expression of $\boxtimes$ when $\mu$ is a $3$-form, in which case $T\boxtimes\mu$ is a 4-form given by
\begin{equation} \label{torsion-box1}
(T\boxtimes\mu)_{ijkl}=T^p{}_{ij}\mu_{pkl}+T^p{}_{kl}\mu_{pij}-T^p{}_{ik}\mu_{pjl}-T^p{}_{jl}\mu_{pik}+T^p{}_{il}\mu_{pjk}+T^p{}_{jk}\mu_{pil},
\end{equation}
as well as the case $\mu$ is a 2-form, where $T\boxtimes\mu$ is a 3-form of the form
\begin{equation} \label{torsion-box2}
(T\boxtimes\mu)_{ijk}=T^p{}_{ij}\mu_{pk}+T^p{}_{jk}\mu_{pi}+T^p{}_{ki}\mu_{pj}.
\end{equation}
\end{lemma}

\noindent{\it Proof:} We give the proof of (\ref{torsion-box2}) and leave (\ref{torsion-box1}) to the reader. For $\mu = \dfrac{1}{2} \mu_{ij} dx^i \wedge dx^j$, we have
\be
d \mu = {1 \over 2} \partial_\alpha \mu_{i j} \, dx^\alpha \wedge d x^i \wedge d x^j.
\ee
We write $D_i W_j = \partial_i W_j - \Gamma(D)^k{}_{ij} W_k$, and obtain
\be
d \mu = {1 \over 2} (D_k \mu_{i j} + \Gamma(D)^\beta{}_{k i} \mu_{\beta j} + \Gamma(D)^\beta{}_{k j} \mu_{i \beta}) \, dx^k \wedge d x^i \wedge d x^j.
\ee
This becomes
\bea
d \mu &=& dx^k \wedge D_k \mu + {1 \over 3!} \left( \Gamma(D)^\beta{}_{k i} \mu_{\beta j} + \Gamma(D)^\beta{}_{j k} \mu_{\beta i} + \Gamma(D)^\beta{}_{i j} \mu_{\beta k}  \right) dx^k \wedge dx^i \wedge d x^j \nonumber\\
&&+ {1 \over 3!} \left( \Gamma(D)^\beta{}_{k j} \mu_{i \beta} + \Gamma(D)^\beta{}_{i k} \mu_{j \beta} + \Gamma(D)^\beta{}_{j i} \mu_{k \beta} \right) dx^k \wedge dx^i \wedge d x^j
\eea
which leads to
\be
d \mu = dx^k \wedge D_k \mu + {1 \over 3!} (T^\beta{}_{k i} \mu_{\beta j} + T^\beta{}_{j k} \mu_{\beta i} + T^\beta{}_{i j} \mu_{\beta k}) dx^k \wedge dx^i \wedge dx^j.
\ee
Q.E.D.


\begin{lemma}\label{22}~\\
In the notation in Lemma \ref{torsiond}, the 4-form $N\boxtimes\varphi$ is of type $(2,2)$.
\end{lemma}
{\it Proof of the Lemma}: By Lemma \ref{torsiond}, we know that
\bea
(N\boxtimes\varphi)_{ijkl}=N^p{}_{ij}\varphi_{pkl}+N^p{}_{kl}\varphi_{pij}-N^p{}_{ik}\varphi_{pjl}-N^p{}_{jl}\varphi_{pik}+N^p{}_{il}\varphi_{pjk} +N^p{}_{jk}\varphi_{pil}.\nonumber
\eea
Since $N\in A^{0,2}(TM)$ and $\varphi$ satisfies Lemma \ref{phihat}, we find
\bea
(J(N\boxtimes\varphi))_{ijkl}&=&N^p{}_{Ji,Jj}\varphi_{p,Jk,Jl}+N^p{}_{Jk,Jl}\varphi_{p,Ji,Jj}-N^p{}_{Ji,Jk}\varphi_{p,Jj,Jl}\nonumber\\
&&-N^p{}_{Jj,Jl}\varphi_{p,Ji,Jk}+N^p{}_{Ji,Jl}\varphi_{p,Jj,Jk} +N^p{}_{Jj,Jk}\varphi_{p,Ji,Jl}\nonumber\\
&=&N^p{}_{ij}\varphi_{pkl}+N^p{}_{kl}\varphi_{pij}-N^p{}_{ik}\varphi_{pjl}-N^p{}_{jl}\varphi_{pik}+N^p{}_{il}\varphi_{pjk} +N^p{}_{jk}\varphi_{pil}\nonumber\\
&=&(N\boxtimes\varphi)_{ijkl}.\nonumber
\eea
As $J$ acts on $(3,1)+(1,3)$-forms as $-1$ and acts on $(2,2)$-forms as $1$, we deduce that $N\boxtimes\varphi$ is a $(2,2)$-form. Q.E.D.

\begin{lemma}\label{uvarphi}~\\
Using the notation in Lemma \ref{torsiond}, we have
\bea
d^c\tilde\o\boxtimes\varphi=\mathcal M(d^c\tilde\o)\boxtimes\varphi=2\alpha\wedge\varphi,\\
d^c\tilde\o\boxtimes\hat\varphi=\mathcal M(d^c\tilde\o)\boxtimes\hat\varphi=2\alpha\wedge\hat\varphi.
\eea
\end{lemma}
{\it Proof}. As we have seen $d^c\tilde\o=J\alpha\wedge\tilde\o$, so the first term in $(d^c\tilde\o)\boxtimes\varphi$ is
\begin{eqnarray}
(d^c\tilde\o)^p_{~ij}\varphi_{pkl}&=&\tilde{g}^{pq}((J\alpha)_q\tilde\o_{ij}+(J\alpha)_i\tilde\o_{jq}+(J\alpha)_j\tilde\o_{qi})\varphi_{pkl}\nonumber\\
&=&\tilde{g}^{pq}\alpha_{Jq}\tilde\o_{ij}\varphi_{pkl}+\alpha_{Ji}\varphi_{Jj,k,l}-\alpha_{Jj}\varphi_{Ji,k,l}\nonumber\\
&=&\tilde{g}^{pq}\alpha_{Jq}\tilde\o_{ij}\varphi_{pkl}-\alpha_{Ji}\hat\varphi_{jkl}+\alpha_{Jj}\hat\varphi_{ikl}.\nonumber
\end{eqnarray}
Hence
\begin{eqnarray}
(d^c\tilde\o\boxtimes\varphi)_{ijkl}&=&\tilde{g}^{pq}\alpha_{Jq}(\tilde\o_{ij}\varphi_{pkl}+\tilde\o_{kl}\varphi_{pij}-\tilde\o_{ik}\varphi_{pjl} -\tilde\o_{jl}\varphi_{pik}+\tilde\o_{il}\varphi_{pjk}+\tilde\o_{jk}\varphi_{pil})\nonumber\\
&&-\alpha_{Ji}\hat\varphi_{jkl}+\alpha_{Jj}\hat\varphi_{ikl}-\alpha_{Jk}\hat\varphi_{lij}+\alpha_{Jl}\hat\varphi_{kij}+\alpha_{Ji}\hat\varphi_{kjl} -\alpha_{Jk}\hat\varphi_{ijl}\nonumber\\
&&+\alpha_{Jj}\hat\varphi_{lik}-\alpha_{Jl}\hat\varphi_{jik}-\alpha_{Ji}\hat\varphi_{ljk}+\alpha_{Jl}\hat\varphi_{ijk}-\alpha_{Jj}\hat\varphi_{ilk} +\alpha_{Ji}\hat\varphi_{jlk}\nonumber\\
&=&\tilde{g}^{pq}\alpha_{Jq}(\tilde\o_{ip}\varphi_{jkl}-\tilde\o_{jp}\varphi_{ikl}+\tilde\o_{kp}\varphi_{ijl}-\tilde\o_{lp}\varphi_{ijk})-3(J\alpha\wedge\hat\varphi)_{ijkl}\nonumber\\
&=&\alpha_{Jq}(J^q_{~i}\varphi_{jkl}-J^q_{~j}\varphi_{ikl}+J^q_{~k}\varphi_{ijl}-J^q_{~l}\varphi_{ijk})+3(\alpha\wedge\varphi)_{ijkl}\nonumber\\
&=&2(\alpha\wedge\varphi)_{ijkl}.\nonumber
\end{eqnarray}

In this proof we only used the fact that $\varphi$ is primitive (\ref{primitive}) so in the same manner we have
\begin{equation}
d^c\tilde\o\boxtimes\hat\varphi=2\alpha\wedge\hat\varphi.\nonumber
\end{equation}
The other identities can be proved similarly. Q.E.D.\\

\noindent Now we are ready to prove Theorem \ref{th:holonomy} (a).

Since $\tilde{\mathfrak D}J=0$, there exists a complex-valued 1-form $\theta=\alpha+\sqrt{-1}\beta$ such that
\bea
\tilde{\mathfrak D}\Omega=\theta\otimes\Omega.\nonumber
\eea
Taking its real and imaginary parts, we get
\bea
\tilde{\mathfrak D}\varphi&=&\alpha\otimes\varphi-\beta\otimes\hat\varphi,\label{dvarphi}\\
\tilde{\mathfrak D}\hat\varphi&=&\beta\otimes\varphi+\alpha\otimes\hat\varphi.\label{dstarvarphi}
\eea
The 1-form $\alpha$ is very easy to find: as $\tilde{\mathfrak D}\tilde g=0$, we know that
\begin{equation}
d|\varphi|^2_{\tilde{g}}=\tilde{\mathfrak D}\tilde{g}(\varphi,\varphi)=2\tilde{g}(\tilde{\mathfrak D}\varphi,\varphi)=2|\varphi|^2_{\tilde{g}}\alpha,\nonumber
\end{equation}
hence we conclude that
\begin{equation}
\alpha=\frac{1}{2}d\log|\varphi|^2_{\tilde{g}}=d\log|\varphi|_{\tilde{g}}=-d\log|\varphi|^2,
\end{equation}
which is the exactly same expression we assigned to $\alpha$ in (\ref{alpha}). To find $\beta$, we plug (\ref{dvarphi}) in (\ref{dconnection}) to get
\bea
0=d\varphi=\alpha\wedge\varphi-\beta\wedge\hat\varphi+\tilde{\mathfrak T}\boxtimes\varphi.\label{varphi2}
\eea
Apply (\ref{T1}) to the Hermitian metric $\tilde g$ with $t=0$, we get
\bea
\tilde{\mathfrak T}=N-U,\label{torsiontilde}
\eea
where $U=\dfrac{1}{4}(d^c\tilde\o+\mathcal M(d^c\tilde\o))$. According to Lemma \ref{uvarphi} we have
\bea
\tilde{\mathfrak T}\boxtimes\varphi&=&N\boxtimes\varphi-\frac{1}{4}d^c\tilde\o\boxtimes\varphi-\frac{1}{4}\mathcal M(d^c\tilde\o)\boxtimes\varphi\nonumber\\
&=&N\boxtimes\varphi-\alpha\wedge\varphi.\nonumber
\eea
Consequently (\ref{varphi2}) can be simplified to
\bea
N\boxtimes\varphi=\beta\wedge\hat\varphi.\nonumber
\eea
By Lemma \ref{22}, the LHS of the above equation is a $(2,2)$-form while the RHS is a $(3,1)+(1,3)$-form. Therefore we conclude that $N\boxtimes\varphi=0$ and $\beta=0$. As a result,
\bea
\tilde{\mathfrak D}\Omega=\alpha\otimes\Omega,\label{tildedomega}
\eea
which implies immediately that
$\tilde{\mathfrak D}\left(\frac{\Omega}{|\Omega|_{\tilde g}}\right)=0$. Q.E.D.

\medskip
\noindent{\it Remark}: Heuristically we can argue as follows. Since $N$ accounts for the non-integrability of $J$, the form $N\boxtimes\varphi$ is responsible for the ``exotic'' component of $d\varphi$ which vanishes automatically in the integrable case. Because $\varphi$ is a $(3,0)+(0,3)$-form, $d\varphi$ would be a $(3,1)+(1,3)$-form if $J$ is integrable. As a result
\be
N\boxtimes\varphi=\textrm{(2,2) component of }d\varphi=0.\nonumber
\ee

\begin{corollary}~\\
The pair $(N,\varphi)$ satisfies
\bea
N^p{}_{ij}\varphi_{pkl}+N^p{}_{kl}\varphi_{pij}&=&0,\label{switch}\\
N^p{}_{ij}\hat\varphi_{pkl}-N^p{}_{kl}\hat\varphi_{pij}&=&0.\label{switch2}
\eea
\end{corollary}
{\it Proof}: In the proof of Theorem (\ref{th:holonomy}) Part (a), we showed that $N\boxtimes\varphi=0$, namely
\bea
N^p{}_{ij}\varphi_{pkl}+N^p{}_{kl}\varphi_{pij}-N^p{}_{ik}\varphi_{pjl}-N^p{}_{jl}\varphi_{pik}+N^p{}_{il}\varphi_{pjk} +N^p{}_{jk}\varphi_{pil}=0.\label{nvarphi}
\eea
Replace $i$ and $j$ in (\ref{nvarphi}) by $Ji$ and $Jj$, by using symmetry of $N$ and $\varphi$, we get instead
\bea
-N^p{}_{ij}\varphi_{pkl}-N^p{}_{kl}\varphi_{pij}-N^p{}_{ik}\varphi_{pjl}-N^p{}_{jl}\varphi_{pik}+N^p{}_{il}\varphi_{pjk} +N^p{}_{jk}\varphi_{pil}=0.\label{nvarphi2}
\eea
By combining (\ref{nvarphi}) and (\ref{nvarphi2}) we prove the corollary. Equation (\ref{switch2}) follows from (\ref{switch}) and Lemma \ref{phihat}. Q.E.D.

\subsection{The curvature in Type IIA geometry}

Next, we prove Theorem \ref{th:holonomy} (b).

\smallskip
As we have seen in Theorem \ref{th:holonomy} (a), in Type IIA geometry, the nowhere vanishing $(3,0)$-form $\Omega/|\Omega|_{\tilde g}$ is parallel under the connection $\tilde{\mathfrak D}$. A direct consequence is that the first Chern form associated to $\tilde{\mathfrak D}$ is identically zero, that is,
\bea
-\frac{\sqrt{-1}}{2}\tilde{\mathfrak R}_{ijkl}\tilde\o^{lk}=\tilde{\mathfrak R}_{ij}{}^\beta{}_\beta=0.
\eea
As $\tilde g=|\varphi|^2g$, one can relate $\tilde{\mathfrak D}$ with $\mathfrak D$ by the conformal change formula. Combining it with (\ref{tildedomega}), it is not hard to see that
\bea
\mathfrak D\Omega=-\frac{1}{2}(\alpha-\sqrt{-1}J\alpha)\otimes\Omega.\label{chern}
\eea
As a consequence, the curvature tensor $\mathfrak R$ satisfies
\bea
-\frac{\sqrt{-1}}{2}\mathfrak R_{ijkl}\o^{lk}=\mathfrak R_{ij}{}^\beta{}_\beta=\frac{1}{2}d(\alpha-\sqrt{-1}J\alpha)_{ij}=-\sqrt{-1}(dd^c\log|\varphi|)_{ij},\label{complexricci}
\eea
and we recover the well-known formula for Ricci curvature in the K\"ahler case. In fact, (\ref{chern}) implies that $\mathfrak D^{0,1}\Omega=0$. For an almost K\"ahler manifold, the Gauduchon line of connections \cite{G} collapses to a point, so $\mathfrak D$ is also the Chern connection (in the almost complex setting), hence $\mathfrak D^{0,1}=\bar\p$, and we conclude that $\Omega$ is a holomorphic section of the canonical bundle associated to $(M,J)$. Theorem \ref{th:holonomy} (b) is proved.

\bigskip

We complete this section with some identities linking the curvature and Nijenhuis tensor.
Recall the globally defined function $u=\log|\varphi|^2$. In this notation we have $\tilde g=e^ug$ and $\alpha=-du$. Furthermore (\ref{complexricci}) can be rewritten as
\bea
\mathfrak R_{ijkl}\o^{lk}&=&(dd^cu)_{ij}=-(\mathfrak D_i(du)_{Jj}-\mathfrak D_j(du)_{Ji}+N^k{}_{ij}u_{Jk})\nonumber\\
&=&-(\nabla^2u)_{i,Jj}+(\nabla^2u)_{j,Ji}-2N^k{}_{ij}u_{Jk}.\label{complexricci2}
\eea
Substitute (\ref{complexricci2}) back to (\ref{ricciakg1}), we get
\bea \label{ricciakg2}
R_{ij}&=&2\mathfrak D^kN_{ijk}-2(N^2_-)_{ij}+\frac{1}{2}(\nabla^2u)_{ij}+\frac{1}{2}(\nabla^2u)_{Ji,Jj}-u_kN^k{}_{ij}.\nonumber
\eea
Since $R_{ij}$ is symmetric, we conclude
\bea
R_{ij}&=&\mathfrak D^k(N_{ijk}+N_{jik})-2(N^2_-)_{ij}+\frac{1}{2}(\nabla^2u)_{ij}+\frac{1}{2}(\nabla^2u)_{Ji,Jj},\label{ricciakg}\\
&=&\nabla^k(N_{ijk}+N_{jik})+2(N^2_-)_{ij}-2(N^2_+)_{ij}+\frac{1}{2}(\nabla^2u)_{ij}+\frac{1}{2}(\nabla^2u)_{Ji,Jj}.\label{ricciakgnabla}
\eea
and that $N$ satisfies
\bea
\mathfrak D^kN_{kij}=\nabla^kN_{kij}=-u^kN_{kij}.\label{divn}
\eea
Therefore the $J$-invariant and $J$-anti-invariant components of the Ricci curvature are given by (\ref{-jinvricci}) and the following refinement of (\ref{jinvricci})
\bea \label{jinvricci2}
(R^J)_{ij}=-2(N^2_-)_{ij}+\left((\nabla^2u)^J\right)_{ij},\qquad
(R^{-J})_{ij}=\mathfrak D^k(N_{ijk}+N_{jik}).
\eea
The scalar curvature is
\bea
R=\Delta u-|N|^2.
\eea
(\ref{complexricci2}) also implies that
\bea
\frac{1}{2}\mathfrak R_{ijkl}\o^{ji}\o^{lk}=\Delta u.\nonumber
\eea
It follows from (\ref{scalarcurv}) that the $\star$-scalar curvature is given by
\bea
\frac{1}{2}R_{ijkl}\o^{ji}\o^{lk}=\Delta u+|N|^2.
\eea

Similarly we can derive the formulae for $\tilde R$, the curvature tensor associated to the conformal metric $\tilde g$:
\bea
\tilde R_{ij}&=&-(\tilde{\mathfrak D}^s-\frac{1}{2}u^s)(N_{isj}+N_{jsi})+\frac{1}{2}((\tilde\nabla^2u)_{Ji,Jj}-3(\tilde\nabla^2u)_{ij}-\tilde\Delta u\tilde g_{ij})\nonumber\\
&&-2(N^2_-)_{ij}-\frac{1}{2}u_iu_j+\frac{1}{2}u_{Ji}u_{Jj}+\frac{1}{2}|du|^2_{\tilde g}\tilde g_{ij},\label{tildericci}\\
\tilde R&=&-4\tilde\Delta u+3|du|^2_{\tilde g}-|N|^2_{\tilde g},\label{tildescalar}
\eea
with (\ref{divn}) becoming
\bea
2\tilde{\mathfrak D}^kN_{kij}=u^kN_{kij},\qquad
2\tilde\nabla^kN_{kij}=3u^kN_{kij}.\nonumber
\eea

\subsection{The Nijenhuis tensor in Type IIA geometry}

As we have seen in previous sections, on an almost K\"ahler manifold, the Nijenhuis tensor $N$ is a $(0,2)$-type $TM$-valued 2-form satisfying the Bianchi identity (\ref{nbianchi}). Moreover, one can define two $J$-invariant symmetric tensors $N^2_+$ and $N^2_-$ satisfying
\bea
\tr N^2_+=|N|^2=2\tr N^2_-.\nonumber
\eea

When an almost K\"ahler structure is enhanced to a Type IIA structure, the integrability of $J$ is improved, hence one should expect more identities satisfied by $N$. For example, we have already seen that certain divergences of $N$ are actually terms of lower order term (\ref{divn}).
In this subsection, we shall derive more identities and differential equations satisfied by $N$, showing that an almost-complex structure coming from a Type IIA geometry is more ``integrable'' than a generic one. We shall also complete the proof of Theorem \ref{th:holonomy} by proving Part (c).

\medskip

First, we show that $N^2_+$ and $N^2_-$ are related to each other:
\begin{proposition}\label{quadraticn}~\\
For any Type IIA structure $(M,\o,\varphi)$, the Nijenhuis tensor $N$ satisfies
\bea
N^2_-=2N^2_+-\frac{1}{4}|N|^2g.
\eea
\end{proposition}
{\it Proof:} In view of (\ref{switch}), we notice that
\be
\varphi_{iap}N^{st}_{~~b}N^p_{~st}=-\varphi_{pst}N^{st}{}_bN^p{}_{ia}=\varphi_{stb}N^{st}{}_pN^p{}_{ia},\nonumber
\ee
so
\be
|\varphi|^{-2}N^{st}_{~~b}N^p_{~st}\varphi_{iap}\varphi_{jcd}\o^{ac}\o^{bd}=|\varphi|^{-2}N^{st}{}_pN^p{}_{ia}\varphi_{stb}\varphi_{jcd}\o^{ac}\o^{bd}.\nonumber
\ee
Apply Lemma \ref{contract1} to both sides, we get
\be
2N^{st}{}_iN_{jst}-N^{stk}N_{kst}g_{ij}=2N_j{}^{st}N_{tsi}+2N^{st}{}_jN_{tsi}.\nonumber
\ee
Converting everything into $N^2_+$ and $N^2_-$, we get
\be
-2N^2_++2N^2_--N^{stk}N_{kst}g=2N^2_+.\nonumber
\ee
By Bianchi identity, we know that
\be
-N^{stk}N_{kst}=N^{stk}(N_{stk}-N_{tsk})=\tr N^2_+-\tr N^2_-=\frac{1}{2}|N|^2.\nonumber
\ee
Consequently
\be
N^2_-=2N^2_+-\frac{1}{4}|N|^2g.\nonumber
\ee
As a corollary, we obtain the inequality
\be
0\leq N^2_+\leq\frac{1}{4}|N|^2g
\ee
since $0\leq N^2_+-N^2_-$. Q.E.D.

\medskip
We now come to the proof of Theorem \ref{th:holonomy} (c), which is a very powerful tool for proving identities involving $N$:

\medskip

Let us choose a frame at a given point as in Lemma \ref{up canonical form}. Since $N$ is a $(0,2)$-type $TM$-valued 2-form satisfying the Bianchi identity, we get the following relations:
\bea
&&0=N_{*jj}=N_{*12}=N_{*21}=N_{*34}=N_{*43}=N_{*56}=N_{*65},\nonumber\\
&&N_{135}=-N_{153}=-N_{146}=N_{164}=-N_{236}=N_{263}=-N_{245}=N_{254},\nonumber\\
&&N_{136}=-N_{163}=N_{145}=-N_{154}=N_{235}=-N_{253}=-N_{246}=N_{264},\nonumber\\
&&N_{315}=-N_{351}=-N_{326}=N_{362}=-N_{416}=N_{461}=-N_{425}=N_{452},\nonumber\\
&&N_{316}=-N_{361}=N_{325}=-N_{352}=N_{415}=-N_{451}=-N_{426}=N_{462},\nonumber\\
&&N_{513}=-N_{531}=-N_{524}=N_{542}=-N_{614}=N_{641}=-N_{623}=N_{632},\nonumber\\
&&N_{514}=-N_{541}=N_{523}=-N_{532}=N_{613}=-N_{631}=-N_{624}=N_{642},\nonumber
\eea
and
\bea
&&N_{113}=-N_{131}=-N_{124}=N_{142}=-N_{214}=N_{241}=-N_{223}=N_{232},\nonumber\\
&&N_{114}=-N_{141}=N_{123}=-N_{132}=N_{213}=-N_{231}=-N_{224}=N_{242},\nonumber\\
&&N_{115}=-N_{151}=-N_{126}=N_{162}=-N_{216}=N_{261}=-N_{225}=N_{252},\nonumber\\
&&N_{116}=-N_{161}=N_{125}=-N_{152}=N_{215}=-N_{251}=-N_{226}=N_{262},\nonumber\\
&&N_{331}=-N_{313}=-N_{342}=N_{324}=-N_{432}=N_{423}=-N_{441}=N_{414},\nonumber\\
&&N_{332}=-N_{323}=N_{341}=-N_{314}=N_{431}=-N_{413}=-N_{442}=N_{424},\nonumber\\
&&N_{335}=-N_{353}=-N_{346}=N_{364}=-N_{436}=N_{463}=-N_{445}=N_{454},\nonumber\\
&&N_{336}=-N_{363}=N_{345}=-N_{354}=N_{435}=-N_{453}=-N_{446}=N_{464},\nonumber\\
&&N_{551}=-N_{515}=-N_{562}=N_{526}=-N_{652}=N_{625}=-N_{661}=N_{616},\nonumber\\
&&N_{552}=-N_{525}=N_{561}=-N_{516}=N_{651}=-N_{615}=-N_{662}=N_{626},\nonumber\\
&&N_{553}=-N_{535}=-N_{564}=N_{546}=-N_{654}=N_{645}=-N_{663}=N_{636},\nonumber\\
&&N_{554}=-N_{545}=N_{563}=-N_{536}=N_{653}=-N_{635}=-N_{664}=N_{646},\nonumber
\eea
with constraints
\bea
N_{135}+N_{351}+N_{513}&=&0,\nonumber\\
N_{136}+N_{361}+N_{613}&=&0.\nonumber
\eea
Furthermore, by evaluating (\ref{switch}) at the given point, we get
\bea
0&=&N_{135}=N_{136}=N_{315}=N_{316}=N_{513}=N_{514},\nonumber\\
0&=&N_{331}-N_{551}=N_{113}-N_{553}=N_{115}-N_{335},\nonumber\\
0&=&N_{114}+N_{554}=N_{116}+N_{336}=N_{332}+N_{552},\nonumber
\eea
Therefore $N$ has only 6 independent components at the given point. Q.E.D.\\

We can choose and name such independent components as
\be
a:=N_{331},\quad b:=N_{332},\quad c:=N_{113},\quad d:=N_{114},\quad e:=N_{115},\quad f:=N_{116}.\nonumber
\ee
It follows that
\be
|N|^2=16(a^2+b^2+c^2+d^2+e^2+f^2).\nonumber
\ee
We can further express $N_+^2$ and $N^2_-$ in terms of these components. For instance, it is straightforward to verify that
\bea
N^2_+=2
\begin{bmatrix}
r^2+a^2+b^2 & 0 & ac+bd & -ad+bc & ae-bf & -af-be \\
0 & r^2+a^2+b^2 & ad-bc & ac+bd & af+be & ae-bf \\
ac+bd & ad-bc & r^2+c^2+d^2 & 0 & ce+df & cf-de \\
-ad+bc & ac+bd & 0 & r^2+c^2+d^2 & -cf+de & ce+df\\
ae-bf & af+be & ce+df & -cf+de & r^2+e^2+f^2 & 0 \\
-af-be & ae-bf & cf-de & ce+df & 0 & r^2+e^2+f^2
\end{bmatrix},\nonumber
\eea
where $r^2=a^2+b^2+c^2+d^2+e^2+f^2=\dfrac{1}{16}|N|^2$. Similarly we can find $N^2_-$ as well. This normal form also allows us to quickly prove that
\be
|N^2_+|^2=48r^4=\frac{3}{16}|N|^4.\label{n2+norm}
\ee

As an application of Theorem \ref{th:holonomy} (c), we prove that $N$ satisfies the following differential equation:
\begin{lemma}~\\
Given a Type IIA structure $(M,\o,\varphi)$, the Nijenhuis tensor $N$ satisfies
\bea
8N^{sti}\nabla_iN_{stj}=8N^{sti}{\mathfrak D}_iN_{stj}={\mathfrak D}_j|N|^2+u_j|N|^2.\label{imp}
\eea
\end{lemma}
{\it Proof}: By (\ref{sym}), we have
\bea
&&{\mathfrak D}_iN_{sjt}-{\mathfrak D}_sN_{ijt}-{\mathfrak D}_jN_{tis}+{\mathfrak D}_tN_{jis}\nonumber\\
=&&{\mathfrak D}_{Ji}N_{Js,j,t}-{\mathfrak D}_{Js}N_{Ji,j,t}-{\mathfrak D}_{Jj}N_{Jt,i,s}+{\mathfrak D}_{Jt}N_{Jj,i,s}.\nonumber
\eea
Contracting this equation with $N^{sti}$, we get
\be
2N^{sti}({\mathfrak D}_iN_{sjt}-{\mathfrak D}_sN_{ijt}+{\mathfrak D}_tN_{jis})=N^{sti}({\mathfrak D}_jN_{tis}-{\mathfrak D}_{Jj}N_{Jt,i,s}).\nonumber
\ee
The LHS can be simplified as follows
\bea
\textrm{LHS}&=&-2N^{sti}{\mathfrak D}_iN_{stj}+2N^{its}{\mathfrak D}_iN_{stj}+2N^{sit}{\mathfrak D}_iN_{jts}\nonumber\\
&=&2N^{tis}({\mathfrak D}_iN_{stj}+{\mathfrak D}_iN_{jst})\nonumber\\
&=&-2N^{sti}{\mathfrak D}_iN_{stj}.\nonumber
\eea
On the other hand, we see that
\bea
N^{sti}{\mathfrak D}_jN_{tis}=-N^{sti}{\mathfrak D}_jN_{tsi}=-\frac{1}{2}{\mathfrak D}_j(N^{sti}N_{tsi})=-\frac{1}{4}{\mathfrak D}_j|N|^2.\nonumber
\eea
Therefore to prove the lemma, we only need to show that
\be
N^{sti}{\mathfrak D}_{Jj}N_{Jt,i,s}=\frac{1}{4}u_j|N|^2,\nonumber
\ee
or equivalently
\be
N^{sti}{\mathfrak D}_jN_{t,s,Ji}=\frac{1}{4}u_{Jj}|N|^2.\label{key}
\ee

We only need to verify (\ref{key}) pointwise. To do so, at any given point, we expand the LHS of (\ref{key}) using the normal form of $\varphi$ in Lemma \ref{up canonical form}. For simplicity of notation, let us write $B={\mathfrak D}_jN$. Clearly $B$ has the same symmetry as $N$, namely it is a $TM$-valued type $(0,2)$-form and it satisfies the Bianchi identity. By Lemma \ref{up canonical form}, we get
\bea
\textrm{LHS}&=& \sum_{s,t,i}N_{sti}(\mathfrak D_jN)_{t,s,Ji}\nonumber\\
&=&N_{331}(B_{332}-B_{134}+B_{431}-B_{233}+B_{341}-B_{244}-B_{442}-B_{143})\nonumber\\
&&+N_{331}(B_{552}-B_{156}+B_{651}-B_{255}+B_{561}+B_{266}-B_{662}-B_{165})\nonumber\\
&&+N_{332}(-B_{331}-B_{234}+B_{432}+B_{133}+B_{342}-B_{144}+B_{441}-B_{243})\nonumber\\
&&-N_{332}(-B_{551}-B_{256}+B_{652}+B_{155}+B_{562}-B_{166}+B_{661}-B_{265})\nonumber\\
&&+\dots\nonumber\\
&=&4N_{331}(B_{332}+B_{552})+4N_{113}(B_{114}+B_{554})+4N_{115}(B_{116}+B_{336})\nonumber\\
&&+4N_{332}(B_{551}-B_{331})+4N_{114}(B_{553}-B_{113})+4N_{116}(B_{335}-B_{115}).\nonumber
\eea
Since
\be
\mathfrak D_j\varphi=\frac{1}{2}(u_j\varphi+u_{Jj}\hat\varphi),\nonumber
\ee
by taking derivative of Lemma \ref{switch}, we get
\be
\mathfrak D_jN^p{}_{ab}\varphi_{pcd}+\mathfrak D_jN^p{}_{cd}\varphi_{pab}=-u_{Jj}N^p{}_{ab}\hat\varphi_{pcd},\nonumber
\ee
or equivalently
\be
B^p{}_{ab}\varphi_{pcd}+B^p{}_{cd}\varphi_{pab}=-u_{Jj}N^p{}_{ab}\hat\varphi_{pcd}.\nonumber
\ee
Evaluating the above equation at the given point using Lemma \ref{up canonical form}, we get the following relations
\bea
\begin{split}B_{331}-B_{551}=-u_{Jj}N_{332},\quad &B_{332}+B_{552}=u_{Jj}N_{331},\\
B_{113}-B_{553}=-u_{Jj}N_{114},\quad &B_{114}+B_{554}=u_{Jj}N_{113},\\
B_{115}-B_{335}=-u_{Jj}N_{116},\quad &B_{116}+B_{336}=u_{Jj}N_{115}.
\end{split}\nonumber
\eea
It follows that
\bea
\textrm{LHS of (\ref{key})}&=&4u_{Jj}(N_{331}^2+N_{332}^2+N_{113}^2+N_{114}^2+N_{115}^2+N_{116}^2)\nonumber\\
&=&\frac{1}{4}u_{Jj}|N|^2.\nonumber
\eea
Q.E.D.\\

\if To end this subsection, we prove an identity involving both the $\mathfrak D$-curvature tensor $\mathfrak R$ and the Nijenhuis tensor $N$.
\begin{lemma}
\bea
\mathfrak R_{tklj}N_i{}^{lt}N^{ijk}=-\frac{1}{16}|N|^4.\label{contractcurv}
\eea
\end{lemma}
{\it Proof of the Lemma}: By Bianchi identity of $N$ we notice that
\bea
\mathfrak R_{tklj}N_i{}^{lt}N^{ijk}&=&\mathfrak R_{tklj}(N^{lt}{}_i-N^{tl}{}_i)(N^{jki}-N^{kji})\nonumber\\
&=&(\mathfrak R_{tklj}+\mathfrak R_{ljtk}-\mathfrak R_{tjlk}-\mathfrak R_{lktj})N^{lt}{}_iN^{jki}\nonumber\\
&=&(\mathfrak R_{tklj}+\mathfrak R_{kltj}+\mathfrak R_{ljtk}+\mathfrak R_{jtlk})N^{lt}{}_iN^{jki}.\nonumber
\eea
Invoking (\ref{bianchi}) we see that
\bea
(\mathfrak R_{tklj}+\mathfrak R_{kltj})N^{lt}{}_iN^{jki}&=&(\mathfrak R_{tklj}+\mathfrak R_{kltj}+\mathfrak R_{ltkj})N^{lt}{}_iN^{jki}\nonumber\\
&=&-(\mathfrak D_tN_{jkl}+\mathfrak D_kN_{jlt}+\mathfrak D_lN_{jtk}-N^p{}_{lt}N_{jkp})N^{lt}{}_iN^{jki}\nonumber\\
&=&\mathfrak D_tN_{kjl}(N^{tl}{}_iN^{kji}-N^{lt}{}_iN^{jki})-\frac{1}{16}|N|^4.\nonumber\\
(\mathfrak R_{ljtk}+\mathfrak R_{jtlk})N^{lt}{}_iN^{jki}&=&(\mathfrak R_{ljtk}+\mathfrak R_{jtlk}+\mathfrak R_{tljk})N^{lt}{}_iN^{jki}\nonumber\\
&=&-(\mathfrak D_lN_{kjt}+\mathfrak D_jN_{ktl}+\mathfrak D_tN_{klj}-N^p{}_{tl}N_{kjp})N^{lt}{}_iN^{jki}\nonumber\\
&=&-\mathfrak D_tN_{kjl}(N^{tl}{}_iN^{kji}-N^{lt}{}_iN^{jki}).\nonumber
\eea
Add these two equations up we prove the desired identity. Q.E.D.
\fi

\section{The flow of the metric in the Type IIA flow}
\setcounter{equation}{0}

The main task of this section is to prove Theorem \ref{th:g-varphi}, which gives explicit formulas for the flows of $\varphi$ and $\tilde g_\varphi$ in terms of the curvature and Nijenhuis tensors.

\subsection{A tensor coefficients ODE for $\varphi$: proof of Theorem \ref{th:g-varphi}(a)\label{step1}}

We begin with the proof of Theorem \ref{th:g-varphi} (a), which gives the flow of $\varphi$.
Since we will be mainly working with the metric $\tilde g$, we shall use $\tilde g$ to raise or lower indices in this subsection.


It is clear that $\varphi(t)$ is closed and primitive for any $t$. We can also assume that $\varphi$ is positive, since this is an open condition and later estimates (\ref{minimum}) will show that this property is preserved along the flow.
Therefore we get a family of Type IIA structures
$(M,\o,\varphi(t))$. So we can apply formulae in Type IIA geometry to
expand the right hand side of the flow equation.
Now the right hand side is given by $d\Lambda
d(|\varphi|^2\hat\varphi)=d\Lambda d(e^u\hat\varphi)$. Combining
(\ref{dstarvarphi}), (\ref{dconnection}), and $\beta=0$, we obtain
\be
d \hat{\varphi} = \alpha \wedge \hat{\varphi} + \tilde{{\frak T}}
\boxtimes \hat{\varphi}
\ee
where $\tilde{{\frak T}} = N - {1 \over 4} (d^c \tilde{\omega} + {\cal
  M}(d^c \tilde{\omega}))$ by (\ref{torsiontilde}). Applying Lemma \ref{uvarphi},
we obtain $d\hat\varphi=N\boxtimes\hat\varphi$ and
\bea
d(e^u\hat\varphi)=e^u(du\wedge\hat\varphi+N\boxtimes\hat\varphi).\nonumber
\eea
To proceed, we need to compute $\Lambda(du\wedge\hat\varphi)$ and $\Lambda(N\boxtimes\hat\varphi)$, which are 2-forms of type $(2,0)+(0,2)$ and of type $(1,1)$ respectively. We have the following lemmas:

\begin{lemma}\label{1term}
\be
(\Lambda(du\wedge\hat\varphi))_{kl}=g^{ji}u_i\varphi_{jkl}.
\ee
\end{lemma}
{\it Proof}: Since $(du\wedge\hat\varphi)_{ijkl}=u_i\hat\varphi_{jkl}-u_j\hat\varphi_{ikl}+u_k\hat\varphi_{ijl}-u_l\varphi_{ijk}$, by definition of $\Lambda$,  we have
\be
(\Lambda(du\wedge\hat\varphi))_{kl}= \frac{1}{2}\omega^{ji}(u_i\hat\varphi_{jkl}-u_j\hat\varphi_{ikl}+u_k\hat\varphi_{ijl}-u_l\hat\varphi_{ijk}).\nonumber
\ee
The last two terms in the above expression are zero since $\hat\varphi$ is primitive. In addition, the first two terms are identical due to the symmetry in switching $i$ and $j$, so it follows that
\bea
(\Lambda(du\wedge\hat\varphi))_{kl}&=&\omega^{ji}u_i\hat\varphi_{jkl}=g^{Jj,i}u_i\hat\varphi_{jkl}=g^{ji}u_i\hat\varphi_{Jj,k,l} \nonumber\\
&=&g^{ji}u_i\varphi_{jkl},\nonumber
\eea
where in the last step we make use of Lemma \ref{phihat}. Q.E.D.

\begin{lemma}\label{2term}
\be
\Lambda(N\boxtimes\hat\varphi)_{kl}=2g^{ji}N^p_{~il}\varphi_{pjk}=-2g^{ji}N^p_{~ik}\varphi_{pjl}.
\ee
\end{lemma}
{\it Proof}: By (\ref{torsion-box1}) and (\ref{switch2}) we see that
\be
\Lambda(N\boxtimes\hat\varphi)_{kl}=\omega^{ji}(N^p_{~kl}\hat\varphi_{pij}-N^p_{~ik}\hat\varphi_{pjl}+N^p_{~il}\hat\varphi_{pjk}).\nonumber
\ee
Notice that the first term above vanishes due to the primitiveness of $\hat\varphi$, and again, the last two terms are identical because of the symmetry of switching $i$ and $j$, so we conclude that
\be
\Lambda(N\boxtimes\hat\varphi)_{kl}=2\omega^{ji}N^p_{~il}\hat\varphi_{pjk}=-2\omega^{ji}N^p_{~il}\varphi_{p,Jj,k}=2g^{ji}N^p_{~il}\varphi_{pjk}.\nonumber
\ee
Here we again used Lemma \ref{phihat} to simplify our expression. Q.E.D.

Combining Lemma \ref{1term} and Lemma \ref{2term}, we see immediately that
\bea
\mu_{kl}&:=&(\Lambda d(e^u\hat\varphi))_{kl}=e^ug^{ji}(u_i\varphi_{jkl}+2N^p_{~il}\varphi_{pjk})\nonumber\\
&=&e^{2u}(u^s\varphi_{skl}+2N^{st}{}_l\varphi_{stk})=e^{2u}(u^s\varphi_{skl}-2N^{st}{}_k\varphi_{stl}).\label{mu}
\eea
To compute $d\Lambda d(e^u\hat\varphi)=d\mu$, we make use of Lemma \ref{torsiond} to get
\bea
(d\mu)_{iab}&=&(\tilde{\mathfrak T}\boxtimes\mu)_{iab}+\sum_{\mathrm{cyc~} i,a,b}\tilde{\frak D}_i\mu_{ab}.\label{2terms}
\eea
The first term in (\ref{2terms}) is already in good shape, since by (\ref{torsion-box2}) we get
\bea
(\tilde{\mathfrak T}\boxtimes\mu)_{iab}&=&\sum_{\mathrm{cyc~} i,a,b}{\tilde{\mathfrak T}}^p{}_{ia}\mu_{pb}=e^{2u}\sum_{\mathrm{cyc~} i,a,b}{\tilde{\mathfrak T}}^p{}_{ia}(u^s\varphi_{spb}-2N^{st}{}_p\varphi_{stb})\nonumber\\
&=&e^{2u}\sum_{\mathrm{cyc~} i,a,b}\varphi_{sta}(2N^{st}{}_p\tilde{\mathfrak T}^p{}_{ib}-u^s\tilde{\mathfrak T}^t{}_{ib}),\label{term1}
\eea
which is linear in $\varphi$. For the second term in (\ref{2terms}), we need
\begin{lemma}
\be
\tilde{\frak D}_i\mu_{ab}=e^{2u}(\varphi_{sab}(\tilde{\frak D}_i+u_i)u^s+2\varphi_{sta}(\tilde{\frak D}_i+u_i)N^{st}{}_b).\label{term3}
\ee
\end{lemma}
{\it Proof}: Plugging in (\ref{mu}), we see that
\bea
\tilde{\frak D}_i\mu_{ab}&=&\tilde{\frak D}_i(e^{2u}(u^s\varphi_{sab}+2N^{st}{}_b\varphi_{sta}))\nonumber\\
&=&e^{2u}(2u_i(u^s\varphi_{sab}+2N^{st}{}_b\varphi_{sta})+u^s\tilde{\frak D}_i\varphi_{sab}+2N^{st}{}_b\tilde{\frak D}_i\varphi_{sta}\nonumber\\
&&+\varphi_{sab}\tilde{\frak D}_iu^s+2\varphi_{sta}\tilde{\frak D}_iN^{st}{}_b)\nonumber\\
&\overset{\textrm{(\ref{tildedomega})}}{=}&e^{2u}(\varphi_{sab}(\tilde{\frak D}_i+u_i)u^s+2\varphi_{sta}(\tilde{\frak D}_i+u_i)N^{st}{}_b).\nonumber
\eea
Q.E.D.\\

Combining (\ref{term1}) and (\ref{term3}) we obtain the evolution equation for $\varphi$ stated in Theorem \ref{th:g-varphi} (a).

\medskip
Next we justify the remark made after Theorem \ref{th:g-varphi}, to the effect that the function $u$ is determined by $\tilde g$. Indeed $\tilde g=e^ug$, and thus
To prove the second part of the statement, we notice that $\tilde g=e^ug$, therefore the volume element associated to $\tilde g$ satisfies
\bea
\mathrm{dvol}_{\tilde g}=e^{3u}\mathrm{dvol}_g=e^{3u}\frac{\o^3}{3!}.\nonumber
\eea
Therefore (in Darboux coordinate) we may write the global function $u$ as
\bea
u=\frac{1}{6}\log\det~\tilde g,\label{u}
\eea
which is entirely determined by $\tilde g$. Hence the metric $g$ is also determined by $\tilde g$, and so is the almost complex structurer $J$ since $\o$ is fixed. It follows that the Nijenhuis tensor $N$, the projected Levi-Civita connection $\tilde{\frak D}$ and its torsion $\tilde{\frak T}$ are also determined by $\tilde g$.

\medskip

For the convenience of later calculations, we derive a more explicit evolution equation for $\varphi$ than what we have in Theorem \ref{th:g-varphi} Part (a).
The starting point is (\ref{2terms}), which can be expanded as
\bea
(d\mu)_{iab}=(N\boxtimes\mu)_{iab}-(U\boxtimes\mu)_{iab}+\sum_{\mathrm{cyc~} i,a,b}\tilde{\frak D}_i\mu_{ab}\label{3terms}
\eea
by using (\ref{torsiontilde}). The first two terms in (\ref{3terms}) can be expressed as follows:

\begin{lemma}\label{nmu}
\be
(N\boxtimes\mu)_{iab}=e^{2u}\sum_{\mathrm{cyc~} i,a,b}\varphi_{iap}(u^sN^p{}_{sb}-2N^{st}{}_bN^p{}_{st}).\label{nmuf}
\ee
\end{lemma}
{\it Proof of the Lemma}: By (\ref{torsion-box2}), we have
\bea
(N\boxtimes\mu)_{iab}&=&\sum_{\mathrm{cyc~} i,a,b}N^p{}_{ia}\mu_{pb}=e^{2u}\sum_{\mathrm{cyc~} i,a,b}N^p{}_{ia}(u^s\varphi_{spb}+2N^{st}{}_b\varphi_{stp})\nonumber\\
&=&e^{2u}\sum_{\mathrm{cyc~} i,a,b}N^p{}_{ia}(-u^s\varphi_{psb}+2N^{st}{}_b\varphi_{pst})\nonumber\\
&\overset{\textrm{(\ref{switch})}}{=}&e^{2u}\sum_{\mathrm{cyc~} i,a,b}\varphi_{iap}(u^sN^p{}_{sb}-2N^{st}{}_bN^p{}_{st}).\nonumber
\eea
Q.E.D.

\begin{lemma}~\\
For any real $2$-form $\mu$, we may write $\mu=\mu^++\mu^-$, where $\mu^+$ and $\mu^-$ are the $(1,1)$ and $(2,0)+(0,2)$ components of $\mu$ respectively. Then
\be
(U\boxtimes\mu)_{iab}=-\frac{1}{2}\tilde\o^{qp}u_q(\tilde\o\wedge\mu)_{piab}-\frac{1}{2}\sum_{\mathrm{cyc~} i,a,b}u_{Ji}(\mu^-)_{Ja,b}.
\ee
In particular, for the specific $\mu$ in (\ref{mu}), we get
\bea
(U\boxtimes\mu)_{iab}=\frac{1}{2}(\tilde\o\wedge\iota_W\mu)_{iab}+\frac{e^{2u}}{2}|du|_{\tilde g}^2\varphi_{iab}-e^{2u}\!\!\!\!\sum_{\mathrm{cyclic~} i,a,b}\!\!\!u_i(u^s\varphi_{sab}+N^{st}{}_b\varphi_{sta}).\label{umu}
\eea
\end{lemma}
{\it Proof of the Lemma}: By (\ref{torsion-box2}) and the definition of $U$ (\ref{uexp2}), we have
\bea
(U\boxtimes\mu)_{iab}&=&\tilde{g}^{pq}\sum_{\mathrm{cyc~} i,a,b}U_{qia}\mu_{pb}\nonumber\\
&=&\frac{1}{4}\tilde{g}^{pq}\sum_{\mathrm{cyc~} i,a,b}(2\alpha_{Jq}\tilde\o_{ia}+\alpha_{Ji}\tilde\o_{aq}+\alpha_{Ja}\tilde\o_{qi}+\alpha_i\tilde g_{aq}-\alpha_a\tilde g_{qi})\mu_{pb}\nonumber\\
&=&\frac{1}{4}\sum_{\mathrm{cyc~} i,a,b}2\tilde\o^{qp}\alpha_q\tilde\o_{ia}\mu_{pb}+\alpha_{Ji}\mu_{Ja,b}-\alpha_{Ja}\mu_{Ji,b} +\alpha_i\mu_{ab}-\alpha_a\mu_{ib}\nonumber\\
&=&\frac{1}{2}(\alpha\wedge\mu)_{iab}+\frac{1}{4}\sum_{\mathrm{cyc~} i,a,b}2\tilde\o^{qp}\alpha_q\tilde\o_{ia}\mu_{pb}+\alpha_{Ji}(\mu_{Ja,b}-\mu_{Jb,a}).\nonumber
\eea
Notice that $(\tilde\o\wedge\mu)_{piab}=\tilde\o_{ia}\mu_{pb}+\tilde\o_{ab}\mu_{pi}+\tilde\o_{bi}\mu_{pa}-\tilde\o_{pa}\mu_{ib} -\tilde\o_{pb}\mu_{ai}-\tilde\o_{pi}\mu_{ba}$, therefore
\be
\tilde\o^{qp}\alpha_q(\tilde\o\wedge\mu)_{piab}=(\alpha\wedge\mu)_{iab}+\tilde\o^{qp}\alpha_q\sum_{\mathrm{cyc~} i,a,b}\tilde\o_{ia}\mu_{pb}.\nonumber
\ee
Hence we find out that
\bea
(U\boxtimes\mu)_{iab}&=&\frac{1}{2}\tilde\o^{qp}\alpha_q(\tilde\o\wedge\mu)_{piab}+\frac{1}{4}\sum_{\mathrm{cyc~} i,a,b}\alpha_{Ji}(\mu_{Ja,b}-\mu_{Jb,a}).\nonumber
\eea
Write $\mu=\mu^++\mu^-$, where $\mu^+$ and $\mu^-$ are the $(1,1)$ and $(2,0)+(0,2)$ components. Then by definition, we have
\be
(\mu^+)_{Ja,b}=(\mu^+)_{Jb,a},\quad (\mu^-)_{Ja,b}=-(\mu^-)_{Jb,a}.\nonumber
\ee
So we conclude that
\be
(U\boxtimes\mu)_{iab}=-\frac{1}{2}\tilde\o^{qp}u_q(\tilde\o\wedge\mu)_{piab}-\frac{1}{2}\sum_{\mathrm{cyc~} i,a,b}u_{Ji}(\mu^-)_{Ja,b}.\nonumber
\ee
Now let us apply this to $\mu=\Lambda d(e^u\hat\varphi)$. It is clear from (\ref{mu}) that
\be
(\mu^-)_{ab}=e^{2u}u^s\varphi_{sab}=e^{2u}(\iota_{\tilde\nabla u}\varphi)_{ab}.\label{mu-}
\ee
Let $W=W^p\p_p$ be the vector field defined by $W^p=-\tilde\o^{qp}u_q=(\tilde\nabla u)^{Jp}$, we see that
\be
-\frac{1}{2}\tilde\o^{qp}u_q(\tilde\o\wedge\mu)_{piab}=\frac{1}{2}(\iota_W(\tilde\o\wedge\mu))_{iab}=\frac{1}{2}(\tilde\o\wedge\iota_W\mu)_{iab} -\frac{1}{2}(du\wedge\mu)_{iab},\nonumber
\ee
hence
\bea
(U\boxtimes\mu)_{iab}&\overset{\textrm{(\ref{mu-})}}{=}&\frac{1}{2}(\tilde\o\wedge\iota_W\mu)_{iab}-\frac{1}{2}(du\wedge\mu)_{iab}-\frac{1}{2}e^{2u}\sum_{\mathrm{cyclic~} i,a,b}u_{Ji}u^{Js}\varphi_{sab}\nonumber\\
&\overset{\textrm{(\ref{mu})}}{=}&\frac{1}{2}(\tilde\o\wedge\iota_W\mu)_{iab}-\frac{e^{2u}}{2}(du\wedge\iota_{\tilde\nabla u}\varphi+Jdu\wedge\iota_{J\tilde\nabla u}\varphi)-e^{2u}\!\!\!\!\sum_{\mathrm{cyclic~} i,a,b}\!\!\!u_iN^{st}{}_b\varphi_{sta}.\nonumber
\eea
From (\ref{rotate}) we know that
\bea
du\wedge\varphi=-Jdu\wedge J\varphi,\nonumber
\eea
by taking interior product with $\tilde\nabla u$, we get
\bea
|du|^2_{\tilde g}\varphi-du\wedge\iota_{\tilde\nabla u}\varphi=\iota_{\tilde\nabla u}(du\wedge\varphi)=-\iota_{\tilde\nabla u}(Jdu\wedge J\varphi)=-Jdu\wedge\iota_{J\tilde\nabla u}\varphi.\nonumber
\eea
Substitute the RHS of the above equation back to the previous one, we get
\bea
(U\boxtimes\mu)_{iab}=\frac{1}{2}(\tilde\o\wedge\iota_W\mu)_{iab}+\frac{e^{2u}}{2}|du|_{\tilde g}^2\varphi_{iab}-e^{2u}\!\!\!\!\sum_{\mathrm{cyclic~} i,a,b}\!\!\!u_i(u^s\varphi_{sab}+N^{st}{}_b\varphi_{sta}).\nonumber
\eea
Q.E.D.\\

Combining (\ref{nmuf}), (\ref{umu}) and (\ref{term3}), we conclude that
\bea
\p_t\varphi_{iab}&=&(d\mu)_{iab}\nonumber\\
&=&-\frac{1}{2}(\tilde\o\wedge\iota_W\mu)_{iab}-\frac{e^{2u}}{2}|du|_{\tilde g}^2\varphi_{iab}+e^{2u}\!\!\!\sum_{\mathrm{cyc~} i,a,b}\!\!\varphi_{iap}(u^sN^p{}_{sb}-2N^{st}{}_bN^p{}_{st})\nonumber\\
&&+e^{2u}\!\!\!\sum_{\mathrm{cyc~} i,a,b}\!\!(\varphi_{sab}(\tilde{\frak D}_i+2u_i)u^s+\varphi_{sta}(2\tilde{\frak D}_i+3u_i)N^{st}{}_b).\label{evolvarphi}
\eea

\subsection{The flow of $\tilde g_\varphi$: proof of Theorem \ref{th:g-varphi}(b) \label{step2}}


By definition of $\tilde g$ (\ref{tildemetric}), we know that
\bea
\p_t\tilde g_{ij}&=&-\p_t\varphi_{iab}\varphi_{jcd}\o^{ac}\o^{bd}-\varphi_{iab}\p_t\varphi_{jcd}\o^{ac}\o^{bd}\nonumber\\
&=&-\p_t\varphi_{iab}\varphi_{jcd}\o^{ac}\o^{bd}+(i\leftrightarrow j).\label{evog1}
\eea
We only need to compute the first term in (\ref{evog1}) as the full expression is the symmetrization of the first term there. This term can be calculated using (\ref{evolvarphi}). It is useful to observe the following:
\begin{lemma}\label{simpli}
Suppose $\lambda$ is a 3-form that can be factorized as the product of a 1-form with $\o$, i.e. $\lambda=\nu\wedge\o$ for some 1-form $\nu$. Then
\be
(\varphi_{iab}\lambda_{jcd}+\lambda_{iab}\varphi_{jcd})\omega^{ac}\omega^{bd}=0.\nonumber
\ee
\end{lemma}
{\it Proof of the lemma.} By our assumption, $\lambda_{iab}=\nu_i\o_{ab}+\nu_a\o_{bi}+\nu_b\o_{ia}$. Therefore
\bea
\lambda_{iab}\varphi_{jcd}\omega^{ac}\omega^{bd}&=&(\nu_i\o_{ab}+\nu_a\o_{bi}+\nu_b\o_{ia})\varphi_{jcd}\omega^{ac}\omega^{bd}\nonumber\\
&=&\nu_i\varphi_{jcd}\o^{dc}-\nu_a\varphi_{jci}\omega^{ac}+\nu_b\varphi_{jid}\o^{bd}\nonumber\\
&=&-2\nu_a\o^{ac}\varphi_{ijc},\nonumber
\eea
where the primitiveness of $\varphi$ is used. After symmetrization in $i$ and $j$, the outcome is zero. Q.E.D.\\

By Lemma \ref{simpli}, we do not need to worry about the first term in (\ref{evolvarphi}). For simplicity of notation, let $F$ be the 3-form defined by
\bea
F_{iab}&=&\!\!\!\!\sum_{\mathrm{cyc~} i,a,b}\!\!\!\!\left(\varphi_{iap}(u^sN^p{}_{sb}-2N^{st}{}_bN^p{}_{st})+\varphi_{sab}(\tilde{\frak D}_i+2u_i)u^s+\varphi_{sta}(2\tilde{\frak D}_i+3u_i)N^{st}{}_b\right)\nonumber\\
&=&\!\!\!\!\sum_{\mathrm{cyc~} i,a,b}\!\!\!\!\left(\varphi_{pab}(u^sN^p{}_{si}-2N^{st}{}_iN^p{}_{st}+(\tilde{\frak D}_i+2u_i)u^p)+\varphi_{sta}(2\tilde{\frak D}_i+3u_i)N^{st}{}_b\right)\!,\label{fterm}
\eea
where (\ref{evolvarphi}) can be rewritten as
\bea
\p_t\varphi=-\frac{1}{2}\tilde\o\wedge\iota_W\mu-\frac{e^{2u}}{2}|du|^2_{\tilde g}\varphi+e^{2u}F,\nonumber
\eea
and we have that
\bea
e^{-2u}\p_t\tilde g_{ij}+|du|^2_{\tilde g}\tilde g_{ij}=-F_{iab}\varphi_{jcd}\o^{ac}\o^{bd}+(i\leftrightarrow j).\label{evog2}
\eea
The goal is to compute $F_{iab}\varphi_{jcd}\o^{ac}\o^{bd}$. By (\ref{fterm}) we know
\bea
&&F_{iab}\varphi_{jcd}\o^{ac}\o^{bd}\nonumber\\
&=&\varphi_{jcd}\o^{ac}\o^{bd}\!\!\!\!\sum_{\mathrm{cyc~} i,a,b}\!\!\!\!\left(\varphi_{pab}(u^sN^p{}_{si}-2N^{st}{}_iN^p{}_{st}+(\tilde{\frak D}_i+2u_i)u^p)+\varphi_{sta}(2\tilde{\frak D}_i+3u_i)N^{st}{}_b\right)\nonumber\\
&=&\varphi_{pab}\varphi_{jcd}\o^{ac}\o^{bd}(u^sN^p{}_{si}-2N^{st}{}_iN^p{}_{st}+(\tilde{\frak D}_i+2u_i)u^p)\nonumber\\
&&+2\varphi_{iap}\varphi_{jcd}\o^{ac}\o^{bd}(u^sN^p{}_{sb}-2N^{st}{}_bN^p{}_{st}+(\tilde{\frak D}_b+2u_b)u^p)\nonumber\\
&&+\varphi_{stb}\varphi_{jcd}\o^{ac}\o^{bd}((4\tilde{\frak D}_a+6u_a)N^{st}{}_i-(2\tilde{\frak D}_i+3u_i)N^{st}{}_a))\label{123terms}\\
&=:&(A)+(B)+(C),\nonumber
\eea
where $(A)$, $(B)$, and $(C)$ denote the first, the second, and the third line in (\ref{123terms}) respectively. $(A)$ can be computed using definition of $\tilde g$ (\ref{tildemetric}) directly:
\bea
(A)&=&-(u^sN_{jsi}-2N^{st}{}_iN_{jst}+(\tilde{\frak D}_i+2u_i)u_j)\nonumber\\
&=&-u^sN_{jsi}+2(N^2_-)_{ij}-2(N^2_+)_{ij}-\tilde{\frak D}_iu_j-2u_iu_j.\label{A}
\eea
To compute $(B)$ and $(C)$, we need to invoke Lemma \ref{contract1}. It follows that
\bea
(B)&=&2\varphi_{iap}\varphi_{jcd}\o^{ac}\o^{bd}(u^sN^p{}_{sb}-2N^{st}{}_bN^p{}_{st}+(\tilde{\frak D}_b+2u_b)u^p)\nonumber\\
&=&\frac{1}{2}\o^{bd}(\o_{ij}\tilde g_{pd}-\o_{pj}\tilde g_{id}-\o_{id}\tilde g_{p_j}+\o_{pd}\tilde g_{ij})(u^sN^p{}_{sb}-2N^{st}{}_bN^p{}_{st}+(\tilde{\frak D}_b+2u_b)u^p)\nonumber\\
&=&\frac{1}{2}(\tilde\o_{ij}J^b{}_p-\tilde\o_{pj}J^b{}_i+\tilde g_{pj}\delta^b{}_i-\tilde g_{ij}\delta^b{}_p)(u^sN^p{}_{sb}-2N^{st}{}_bN^p{}_{st}+(\tilde{\frak D}_b+2u_b)u^p)\nonumber\\
&=&\frac{1}{2}(-4N^{st}{}_iN_{jst}+2N^{stp}N_{pst}\tilde g_{ij}+\tilde\o_{ij}\tilde\o^{bp}\tilde{\frak D}_bu_p+(\tilde{\frak D}_i+2u_i)u_j+(\tilde{\frak D}_{Ji}+2u_{Ji})u_{Jj}\nonumber\\
&&-(\tilde{\frak D}^su_s+2|du|^2_{\tilde g})\tilde g_{ij})\nonumber\\
&=&2(N^2_+)_{ij}-2(N^2_-)_{ij}+\frac{1}{2}((\tilde{\frak D}_i+2u_i)u_j+(\tilde{\frak D}_{Ji}+2u_{Ji})u_{Jj}+\tilde\o_{ij}\tilde\o^{bp}\tilde{\frak D}_bu_p)\nonumber\\
&&-\frac{1}{2}(|N|^2_{\tilde g}+\tilde{\frak D}^su_s+2|du|^2_{\tilde g})\tilde g_{ij}\label{B}
\eea
and
\bea
(C)&=&\varphi_{stb}\varphi_{jcd}\o^{ac}\o^{bd}((4\tilde{\frak D}_a+6u_a)N^{st}{}_i-(2\tilde{\frak D}_i+3u_i)N^{st}{}_a))\nonumber\\
&=&\frac{1}{4}\o^{ac}(\o_{sj}\tilde g_{tc}-\o_{tj}\tilde g_{sc}-\o_{sc}\tilde g_{tj}+\o_{tc}\tilde g_{sj})((4\tilde{\frak D}_a+6u_a)N^{st}{}_i-(2\tilde{\frak D}_i+3u_i)N^{st}{}_a))\nonumber\\
&=&\frac{1}{4}(\tilde\o_{sj}J^a{}_t-\tilde\o_{tj}J^a{}_s+\tilde g_{tj}\delta^a{}_s-\tilde g_{sj}\delta^a{}_t)((4\tilde{\frak D}_a+6u_a)N^{st}{}_i-(2\tilde{\frak D}_i+3u_i)N^{st}{}_a))\nonumber\\
&=&-(2\tilde{\frak D}^s+3u^s)N_{isj}.\label{C}
\eea
Combining (\ref{A}), (\ref{B}), (\ref{C}), and (\ref{evog2}), we get
\bea
\p_t\tilde g_{ij}&=&e^{2u}\bigg[2(\tilde{\frak D}^k+2u^k)(N_{ikj}+N_{jki})+\frac{1}{2}(\tilde{\frak D}_iu_j+\tilde{\frak D}_ju_i-\tilde{\frak D}_{Ji}u_{Jj}-\tilde{\frak D}_{Jj}u_{Ji})\nonumber\\
&&+2u_iu_j-2u_{Ji}u_{Jj}+(\tilde{\frak D}^su_s+|du|^2_{\tilde g}+|N|^2_{\tilde g})\tilde g_{ij}\bigg].\label{evog3}
\eea
Equation (\ref{evog3}) is self-contained in the sense that its RHS is entirely determined by the metric $\tilde g$ and the $\varphi$-dependence is fully eliminated. However, to study its analytic behavior, we need to rewrite it in a more familiar form as we have in the case of Ricci flow or the L\^e-Wang flow \cite{LW}. Moreover, we would like to replace all the $\tilde{\frak D}$-derivatives to $\tilde\nabla$-derivatives as it is more convenient for us to apply the conformal change technique.

Notice that
\bea
\tilde{\frak D}_iu_j=\tilde\nabla_iu_j-u^kN_{ikj}+\frac{1}{4}(u_iu_j+u_{Ji}u_{Jj}-|du|^2_{\tilde g}\tilde g_{ij}),\nonumber
\eea
therefore (\ref{evog3}) can be rephrased as
\bea
\p_t\tilde g_{ij}&=&e^{2u}\bigg[(2\tilde{\frak D}^k+3u^k)(N_{ikj}+N_{jki})+(\tilde\nabla^2u)_{ij}-(\tilde\nabla^2u)_{Ji,Jj}\nonumber\\
&&2u_iu_j-2u_{Ji}u_{Jj}+(\tilde\Delta u+|N|^2_{\tilde g})\tilde g_{ij}\bigg].\nonumber
\eea
Taking (\ref{tildericci}) into account, we obtain the desired formula
\bea \label{evog4}
\p_t\tilde g_{ij}&=&e^{2u}\bigg[-2\tilde R_{ij}-2(\tilde\nabla^2u)_{ij}+4u^k(N_{ikj}+N_{jki})+u_iu_j-u_{Ji}u_{Jj}-4(N^2_-)_{ij}\nonumber\\
&&+(|du|^2_{\tilde g}+|N|^2_{\tilde g})\tilde g_{ij}\bigg].
\eea
Q.E.D.\\

Recalling that $u= (1/6) \log \det \tilde{g}$, we can derive from (\ref{evog4}) that
\bea
\p_tu&=&\frac{e^{2u}}{3}\left(-\tilde\Delta u-\tilde R+3|du|^2_{\tilde g}+2|N|^2_{\tilde g}\right)\nonumber\\
&\overset{\textrm{(\ref{tildescalar})}}{=}&e^{2u}(\tilde\Delta u+|N|^2_{\tilde g}).\label{evolutilde}
\eea
By the maximum principle, we immediately prove the following estimate
\begin{lemma}
Suppose $\varphi(t)$ is a solution to the source-free Type IIA flow on $M\times[0,T]$. Then
\bea
|\varphi(t)|^2\geq\min_M|\varphi_0|^2\label{minimum}
\eea
for any $t\in[0,T]$.
\end{lemma}
This lemma has the important consequence that if $\varphi(t)$ is a
solution to the source-free Type IIA flow on $M \times [0,T]$ with primitive
closed initial data, then $\varphi(t)$ remains
positive on $M \times [0,T]$, which allows us to define the
almost-complex structure $J$ and the metric $g$. Indeed, the lemma
implies that
$\sqrt{-\lambda_\varphi} = {1 \over 2} |\varphi|^2 {\omega^3 \over
  3!}$ cannot pass through zero.


\subsection{Conformal transformation to a perturbed Ricci flow\label{step3}}

In this subsection, we wish to establish the uniqueness of the flow (\ref{evog4}). Besides the Ricci curvature, the right hand side of (\ref{evog4}) also contains the 2nd order term $\tilde\nabla^2u$, which cannot be reparametrized away since this would change the symplectic structure and the
reparametrized flow would be non-local. Therefore we need a different technique to deal with the Hessian term, namely we absorb it in the Ricci tensor by a conformal change of metric.

\smallskip
More specifically, we consider a family of conformal Hermitian metrics $g^{(s)}=e^{su}g$, where we have $g^{(0)}=g$ and $g^{(1)}=\tilde g$. Notice that all the metrics $g^{(s)}$ are equivalent except for $s=0$, in which case we need the pair $(g,u)$. This is because one can solve $u$ from $g^{(s)}$ when $s\neq0$ by
\bea
u=\frac{1}{6s}\log\det~g^{(s)}.\nonumber
\eea
Thus we only need to show the short-time existence and uniqueness of any of the flows satisfied by $g^{(s)}$ with $s\neq 0$, or that for the coupled flow $(g,u)$. To begin with, we first compute the evolution equation satisfied by the pair $(g,u)$.
\bea
&&\p_tg_{ij}\nonumber\\
&=&\p_t(e^{-u}\tilde g_{ij})=e^{-u}(\p_t\tilde g_{ij}-\p_tu\cdot\tilde g_{ij})\nonumber\\
&=&e^u\bigg[-2\tilde R_{ij}-2(\tilde\nabla^2u)_{ij}+4u^k(N_{ikj}+N_{jki})+u_iu_j-u_{Ji}u_{Jj}-4(N^2_-)_{ij}\nonumber\\
&&+(|du|^2_{\tilde g}-\tilde\Delta u)\tilde g_{ij}\bigg]\nonumber\\
&=&e^u\bigg[-2R_{ij}+2(\nabla^2u)_{ij}+u_iu_j-u_{Ji}u_{Jj}+4u^k(N_{ikj}+N_{jki})-4(N^2_-)_{ij}\bigg]\label{evolveg}\\
&\overset{\textrm{(\ref{ricciakg})}}{=}&2e^u\bigg[-(R^{-J})_{ij}+((\nabla^2u)^{-J})_{ij}+((du\otimes du)^{-J})_{ij}+2u^s(N_{isj}+N_{jsi})\bigg].\label{evolveg2}
\eea
Meanwhile (\ref{evolutilde}) can be rewritten as
\bea
\p_tu=e^u(\Delta u+2|du|^2+|N|^2).\label{evolu}
\eea
Using the same method, we can derive that
\bea
\p_tg^{(s)}_{ij}&=&e^{(s+1)u}\bigg[-2R^{(s)}_{ij}+(2-4s)((\nabla^{(s)})^2u)_{ij}+(1+2s-2s^2)u_iu_j-u_{Ji}u_{Jj}\nonumber\\
&&+4u^k(N_{ikj}+N_{jki})-4(N^2_-)_{ij}+s(|du|^2_{g^{(s)}}+|N|^2_{g^{(s)}})g^{(s)}_{ij}\bigg]\label{evols}\\
&=&e^{(s+1)u}\bigg[2\left((-R^{(s)}+(1-2s)(\nabla^{(s)})^2u+(1+s-s^2)du\otimes du)^{-J}\right)_{ij}\nonumber\\
&&+4u^k(N_{ikj}+N_{jki})+s\left(\Delta^{(s)}u+2(1-s)|du|^2_{g^{(s)}}+|N|^2_{g^{(s)}}\right)g^{(s)}_{ij}\bigg].
\eea
Formulae (\ref{ricciakgnabla}) and (\ref{divn}) now take the form
\bea
R^{(s)}_{ij}&=&-\nabla^{(s)}_k(N_i{}^k{}_j+N_j{}^k{}_i)+(\frac{1}{2}-2s)((\nabla^{(s)})^2u)_{ij}+\frac{1}{2}((\nabla^{(s)})^2u)_{Ji,Jj} -\frac{s}{2}\Delta^{(s)}ug^{(s)}_{ij}\nonumber\\
&&+\frac{s}{2}(1-2s)u_iu_j+\frac{s}{2}u_{Ji}u_{Jj}-\frac{s}{2}(1-2s)|du|^2_{g^{(s)}}g^{(s)}_{ij}+\frac{3s}{2}u^k(N_{ikj}+N_{jki})\nonumber\\
&&+2(N^2_-)_{ij}-2(N^2_+)_{ij},\label{riccis}\\
\nabla^{(s)}_kN^k{}_{ij}&=&\frac{5s-2}{2}u^kN_{kij}.\label{divns}
\eea
In particular for $s=\dfrac{1}{2}$, from (\ref{evols}) we know that the metric $\check{g}:=g^{(\frac{1}{2})}$ evolves by
\bea
\p_t\check g_{ij}&=&e^{\frac{3}{2}u}\bigg[-2\check R_{ij}+\frac{3}{2}u_iu_j-u_{Ji}u_{Jj}+4u^k(N_{ikj}+N_{jki})-4(N^2_-)_{ij}\nonumber\\
&&+\frac{1}{2}\left(|du|^2_{\check g}+|N|^2_{\check g}\right)\check g_{ij}\bigg],\label{evolucheck}
\eea
where the only 2nd order term on RHS is the Ricci curvature term. We stress that it is important to keep in mind the fact that $\check g_{ij}$ arises from a conformal change from a Type IIA geometry.

\subsection{An integrability condition: proof of Theorem \ref{th:integrability}}

We now prove Theorem \ref{th:integrability}, which provides the key integrability condition needed later to establish the uniqueness of the Type IIA solutions to the flow (\ref{evolucheck}) of the metrics $\check g_{ij}$.

\smallskip


We know that for the Ricci flow, the integrability operator $L$ comes from the contracted Bianchi identity. Since our flow (\ref{evolucheck}) can be viewed as a deformation of the Ricci flow, our $L$ should be a deformation of the contracted Bianchi operator.
Let us simplify our notation in (\ref{evolucheck}) by introducing the tensor $S$ defined as
\bea
S_{ij}:=\frac{3}{2}u_iu_j-u_{Ji}u_{Jj}+4u^k(N_{ikj}+N_{jki})-4(N^2_-)_{ij}+\frac{1}{2}\left(|du|^2_{\check g}+|N|^2_{\check g}\right)\check g_{ij},\label{sexp}
\eea
so (\ref{evolucheck}) can be written as
\bea
\p_t\check g_{ij}&=&e^{\frac{3}{2}u}(-2\check R_{ij}+S_{ij}),
\eea
where we can think of $S$ as the lower order deformation term of the Ricci curvature. Let $L_0$ denote the contracted Bianchi identity operator defined by
\bea
L_0(P)_j:=2\check g^{ik}\check\nabla_kP_{ij}-\check g^{ik}\check\nabla_jP_{ik}\nonumber
\eea
for any symmetric 2-tensor $P$. We know that $L_0(-2\check R_{ij})=0$. Now we would like to look for a zeroth order linear operator $Z$ such that $(L_0+Z)(-2\check R_{ij}+S_{ij})$ is of degree 1 in the metric $\check g$. To do so we need to compute $L_0(S_{ij})$ first.
\begin{proposition}
\bea
L_0(S)_j&=&4u^s\check R_{sj}-8\check R_{sk}N^{ks}{}_j-\frac{8}{3}u_j\check R-16u^s(N^2_-)_{sj}+2u_su_kN^{ks}{}_j\nonumber\\
&&+\frac{1}{3}u_j(2|du|^2_{\check g}-5|N|^2_{\check g}).\label{l0s}
\eea
\end{proposition}
{\it Proof of the Proposition}: We apply $L_0$ to each term of $S$ in (\ref{sexp}) to get
\bea
\frac{3}{2}L_0(u_iu_j)&=&3\check\Delta u\cdot u_j,\nonumber\\
-L_0(u_{Ji}u_{Jj})&=&2\check g^{ik}u_k(\check\nabla^2u)_{Ji,Jj}+\check\nabla_j|d u|^2_{\check g}-2\check g^{ik}u_su_t\check\nabla_i(J^s{}_jJ^t{}_k)\nonumber\\
&=&2\check g^{ik}u_k(\check\nabla^2u)_{Ji,Jj}+\check\nabla_j|d u|^2_{\check g}+4u_su_kN^{ks}{}_j,\nonumber\\
4L_0(u^k(N_{ikj}+N_{jki}))&=&8(\check\nabla^2u)_{sk}N^{ks}{}_j+8u^s\check\nabla^i(N_{isj}+N_{jsi}),\nonumber\\
-4L_0(N^2_-)_j&=&-8\check g^{ik}\check\nabla_k(N^2_-)_{ij}+2\check\nabla_j|N|^2_{\check g},\nonumber\\
\frac{1}{2}L_0\left((|du|^2_{\check g}+|N|^2_{\check g})\check g\right)_j&=&-2\check\nabla_j|du|^2_{\check g}-2\check\nabla_j|N|^2_{\check g}.\nonumber
\eea
To obtain these expressions, we need to use that
\bea
\nabla_iJ^s{}_t&=&-2N_{Ji}{}^s{}_t,\nonumber\\
\check\nabla_iJ^s{}_t&=&-2N_{Ji}{}^s{}_t+\frac{1}{4}(u_{Jt}\delta^s{}_i+g_{it}u^{Js}-\o_{it}u^s-u_tJ^s{}_i)\nonumber
\eea
and we raise and lower indices using the metric $\check g$. Combining the calculation above, we get
\bea
L_0(S)_j&=&3u_j\check\Delta u+2u^i(\check\nabla^2u)_{Ji,Jj}+8(\check\nabla^2u)_{sk}N^{ks}{}_j+8u^s\check\nabla^i(N_{isj}+N_{jsi})\nonumber\\
&&-8\check\nabla^i(N^2_-)_{ij}-\check\nabla_j|du|^2_{\check g}+4u_su_kN^{ks}{}_j.\label{l0s1}
\eea
Take $s=\dfrac{1}{2}$ in (\ref{riccis}) and (\ref{divns}), we get
\bea
\check R_{ij}&=&-(\check\nabla^k-\frac{3}{4}u^k)(N_{ikj}+N_{jki})-\frac{1}{2}(\check\nabla^2u)_{ij}+\frac{1}{2}(\check\nabla^2u)_{Ji,Jj}-\frac{1}{4}\check\Delta u\cdot\check g_{ij}\nonumber\\
&&+\frac{1}{4}u_{Ji}u_{Jj}+2(N^2_-)_{ij}-2(N^2_+)_{ij},\label{riccicheck}\\
\check R&=&-\frac{3}{2}\check\Delta u+\frac{1}{4}|du|^2_{\check g}-|N|^2_{\check g},\\
\check\nabla^kN_{kij}&=&\frac{1}{4}u^kN_{kij}.\label{divncheck}
\eea
Using (\ref{riccicheck}) and (\ref{divncheck}), Equation (\ref{l0s1}) can be rearranged as
\bea
L_0(S)_j&=&3u_j\check\Delta u+2u^i((\check\nabla^2u)_{Ji,Jj}-(\check\nabla^2u)_{ij})-4u^s\check\nabla^i(N_{jis}+N_{sij}+3N_{ijs})\nonumber\\
&&+8(\check\nabla^2u)_{sk}N^{ks}{}_j-8\check\nabla^i(N^2_-)_{ij}+4u_su_kN^{ks}{}_j\nonumber\\
&=&4u^s\check R_{sj}+4u_j\check\Delta u+8(\check\nabla^2u)_{sk}N^{ks}{}_j-8\check\nabla^i(N^2_-)_{ij}+4u_su_kN^{ks}{}_j\nonumber\\
&&+8u^s((N^2_+)_{sj}-(N^2_-)_{sj})\nonumber\\
&=&4u^s\check R_{sj}-\frac{8}{3}u_j\check R+8(\check\nabla^2u)_{sk}N^{ks}{}_j-8\check\nabla^i(N^2_-)_{ij}+4u_su_kN^{ks}{}_j\nonumber\\
&&+8u^s((N^2_+)_{sj}-(N^2_-)_{sj})+\frac{2}{3}u_j(|du|^2_{\check g}-4|N|^2_{\check g}).\label{l0s2}
\eea
By (\ref{riccicheck}), we also know that
\bea
(\check\nabla^2u_{sk})N^{ks}{}_j&=&-(\check\nabla^p-\frac{3}{4}u^p)(N_{spk}+N_{kps})N^{ks}{}_j-\check R_{sk}N^{ks}{}_j-\frac{1}{4}u^ku^sN_{ksj}\nonumber\\
&=&-\check R_{sk}N^{ks}{}_j-\check\nabla^p(2N_{kps}+N_{psk})N^{ks}{}_j-\frac{3}{4}u^p(N^2_++N^2_-)_{pj}-\frac{1}{4}u^ku^sN_{ksj}\nonumber\\
&=&-\check R_{sk}N^{ks}{}_j+2N^{ks}{}_j\check\nabla^pN_{ksp}-\frac{1}{2}u^p(N^2_-+2N^2_+)_{pj}-\frac{1}{4}u^ku^sN_{ksj}\nonumber\\
&=&-\check R_{sk}N^{ks}{}_j+2\check\nabla^p(N^2_+)_{pj}-2N^{ksp}\check\nabla_pN_{ksj}-\frac{1}{2}u^p(N^2_-+2N^2_+)_{pj}-\frac{1}{4}u^ku^sN_{ksj}\nonumber\\
&=&-\check R_{sk}N^{ks}{}_j+\check\nabla^i(N^2_-)_{ij}-2N^{ksp}\check\nabla_pN_{ksj}+\frac{1}{4}\check\nabla_j|N|^2_{\check g}-\frac{1}{4}u^ku^sN_{ksj}\nonumber\\
&&-\frac{1}{2}u^s(N^2_-+2N^2_+)_{sj}.\label{inter}
\eea
Plugging (\ref{inter}) in (\ref{l0s2}), we get
\bea
L_0(S)_j&=&4u^s\check R_{sj}-8\check R_{sk}N^{ks}{}_j-\frac{8}{3}u_j\check R-16N^{ksp}\check\nabla_pN_{ksj}+2\check\nabla_j|N|^2_{\check g}\nonumber\\
&&-12u^s(N^2_-)_{sj}+2u_su_kN^{ks}{}_j+\frac{2}{3}u_j(|du|^2_{\check g}-4|N|^2_{\check g}).\label{l0s3}
\eea
To deal with the remaining second order terms in (\ref{l0s3}), we rewrite (\ref{imp}) using $\check\nabla$-derivatives as
\bea
8N^{ksp}\check\nabla_pN_{ksj}+2u_j|N|^2_{\check g}-2u^i(N^2_-)_{ij}=\frac{3}{2}u_j|N|^2_{\check g}+\check\nabla_j|N|^2_{\check g}.\nonumber
\eea
Incorporating this identity, we see (\ref{l0s3}) becomes
\bea
L_0(S)_j&=&4u^s\check R_{sj}-8\check R_{sk}N^{ks}{}_j-\frac{8}{3}u_j\check R-16u^s(N^2_-)_{sj}+2u_su_kN^{ks}{}_j\nonumber\\
&&+\frac{1}{3}u_j(2|du|^2_{\check g}-5|N|^2_{\check g}).\nonumber
\eea
Q.E.D.\\

Let $Z$ be the zeroth order linear operator defined by
\bea
Z(P)_j:=2u^iP_{ij}-4N^{st}{}_jP_{st}-\frac{4}{3}u_j\check g^{st}P_{st},\nonumber
\eea
then (\ref{l0s}) says
\bea
L_0(S)_j=Z(2\check R_{**})_j-16u^s(N^2_-)_{sj}+2u_su_kN^{ks}{}_j+\frac{1}{3}u_j(2|du|^2_{\check g}-5|N|^2_{\check g}).\nonumber
\eea
Consider the first linear operator $L_1=L_0+Z$, then
\bea
L_1(-2\check R_{**}+S)_j&=&L_0(-2\check R_{**}+S)_j+Z(-2\check R_{**}+S)_j\nonumber\\
&=&L_0(S)_j-Z(2\check R_{**})_j+Z(S)_j\nonumber\\
&=&-16u^s(N^2_-)_{sj}+2u_su_kN^{ks}{}_j+\frac{1}{3}u_j(2|du|^2_{\check g}-5|N|^2_{\check g})\nonumber\\
&&+8u^s(2N^2_++N^2_-)_{sj}-2u_su_kN^{ks}{}_j-\frac{2}{3}u_j(|du|^2_{\check g}+2|N|^2_{\check g})\nonumber\\
&=&-u_j|N|^2_{\check g}\nonumber
\eea
is of first order in $\check g$. Therefore if we define the first order linear operator $L$ by
\bea
L(P)=L_1(e^{-\frac{3}{2}u}P),\nonumber
\eea
then $L$ is an integrability condition for the flow (\ref{evolucheck}). Theorem \ref{th:integrability} is proved. Q.E.D.

\subsection{Return to the proof of Theorem \ref{th:short}: uniqueness}
\label{7.5}

It is now easy to establish the uniqueness part in Theorem \ref{th:short}.

\smallskip

Assume that we have two closed, primitive, and positive solutions $\varphi(t)$ and $\varphi'(t)$ of the Type IIA flow on some time interval $[0,T)$ for some $T>0$, with the same initial data $\varphi(0)=\varphi'(0)$.
By Theorem \ref{th:g-varphi}, the corresponding pairs $(\varphi(t),\tilde g_\varphi(t))$ and $(\varphi'(t),\tilde g_{\varphi'}(t))$ satisfy the flows in Theorem \ref{th:g-varphi}.
Since the geometries $(\o,J_\varphi,g_\varphi)$ and $(\o,J_{\varphi'},g_{\varphi'})$ are by definition Type IIA geometries, the corresponding flows for $\check g_{\varphi}(t)$ and $\check g_{\varphi'}(t)$ satisfy the integrability condition in Theorem \ref{th:integrability}. Since the principal symbols in the flow of $\check g_{ij}$ and the integrability condition $L$ are the same up to a multiplicative factor as their counterparts in the Ricci flow, it follows that the flow of $\check g_{ij}$ together with the integrability condition $L$ satisfy all the conditions in the Hamilton-Nash-Moser theorem (\cite{Ha1}, Theorem 5.1). By the uniqueness part in this theorem, we conclude that $\check g_{ij}(t)$ and $\check g_{ij}'(t)$ must be equal. But then $\varphi$ and $\varphi'(t)$ satisfy the same ODE with the same initial data and hence must be equal. Q.E.D.

\subsection{Monotonicity formulas}

Recall that the function $u$ evolves by
\bea
\p_t u=e^u(\Delta u+2|\nabla u|^2+|N|^2).\label{uevolve}
\eea
From (\ref{uevolve}) one can derive a number of things.
\begin{proposition}
\be
\min_M u(t)\geq \min_M u(0),
\ee
\end{proposition}
{\it Proof of the Proposition}: Apply maximum principle to (\ref{uevolve}). Q.E.D.\\

This proposition can be interpreted that if $\varphi_0$ is initially positive, then $\varphi$ stays positive as long as the flow exists. This is because that the only possibility for $\varphi$ leaving the positive cone is that it first hits the wall of degeneracy defined $|\varphi|=0$, which contradicts the above $C^0$-estimate.

Like its Type IIB counterpart \cite{FeiPicard}, one has the following monotonicity formulas for the dilaton functional along the flow.
\begin{proposition}
\be
\p_t\int_M e^{pu}\frac{\o^3}{3!}=p\int_Me^{(p+1)u}\left((1-p)|\nabla u|^2+|N|^2\right)\frac{\o^3}{3!}.
\ee
If we denote $\int e^{pu}$ by $E_p$, then it follows that $E_p$ is monotonely non-increasing along the flow for $p<0$ and it is monotonely non-decreasing along the flow for $0<p\leq 1$. In particular, the Hitchin's functional $E_1$ \cite{Hitchin} is monotonely non-decreasing along the source-free Type IIA flow.
\end{proposition}

\section{Estimates for the Type IIA flow}
\setcounter{equation}{0}
\label{s:Shi}

Recall that the Type IIA flow becomes the following flow for the pair
$(g(t),u(t))$:
\bea \label{evolution-longtime}
\partial_t g_{ij} &=& e^u \left[-2R_{ij} + 2 \nabla_i \nabla_j u - 4 (N^2_-)_{ij} + u_i u_j - u_{Ji} u_{Jj} + 4 u_p (N_i{}^p{}_j + N_j{}^p{}_i )\right] \nonumber\\
(\p_t - e^u \Delta) u &=& e^u \left[ 2 |\nabla u |^2 + |N|^2 \right]
\eea
In this section, we show that if $|u| + |Rm(g)| \leq C$ remains
bounded on $[0,T)$, then the flow can be extended to $[0,T+\epsilon)$
    for some $\epsilon>0$.
    \smallskip
    \par To start, we examine some consequences of
    the boundedness of the Riemann curvature tensor. We note that by equation (\ref{scalarcurv}), a bound on $Rm$ implies a bound on $|N|^2$. Therefore we may assume
\be
|u|+ |N|^2 + |Rm| \leq C.
\ee
Next, since the Ricci curvature $\textrm{Ric}$ is also bounded, so are
its $J$-invariant and $J$-anti-invariant parts. From
(\ref{jinvricci2}), we know that the $J$-invariant part of the Ricci
curvature $(R^J)_{ij}$ is given by $(\mathrm{Ric}^J)= (\nabla^2 u)^J -2
N^2_-$. As $|N|^2$ is already bounded, we conclude that $\dfrac{1}{2}(\nabla^2u)_{ij}+\dfrac{1}{2}(\nabla^2u)_{Ji,Jj}$ is bounded, namely, the $J$-invariant part of $\nabla^2u$ is bounded. Consequently $\Delta u$ is also bounded.
\bigskip
\par Our goal will be to obtain bounds on all derivatives of $u$, $N$,
$Rm$. For this, we must first compute the evolution equations of
$\nabla^k u$, $\nabla^k N$ and $\nabla^k Rm$.

\subsection{The evolution of the derivatives of $u$}
In this section, we compute the evolution of $|\nabla u|^2$ and
$|\nabla \nabla u|^2$.

\subsubsection{The evolution of the gradient of $u$}
We start with
\be \label{evol-nabla-u-1}
\p_t |\nabla u|^2 = 2 g^{ij} \nabla_i \dot{u} \nabla_j u - g^{i a}
\dot{g}_{ab} g^{bj} u_i u_j.
\ee
The differentiated evolution of $u$ is
 \bea \label{diff-u-once}
\nabla_i \dot{u} = e^u \left( \nabla_i \Delta u + 2 \nabla_i |\nabla
  u|^2 + \nabla_i |N|^2 \right) + e^u \left( \Delta u + 2 |\nabla
  u|^2 + |N|^2 \right)\,  u_i
\eea
Commuting derivatives
\be
\nabla_i \Delta u = g^{pq} \nabla_i \nabla_p \nabla_q u = \Delta
\nabla_i u - g^{pq} R_{ip}{}^\lambda{}_q u_\lambda
\ee
Therefore the first term in (\ref{evol-nabla-u-1})
\bea
2 g^{ij} \nabla_i \dot{u} \nabla_j u &=& 2 e^u \left( g^{ij} \Delta
\nabla_i u \nabla_j u  - g^{pq} R_{p}{}^{i \lambda}{}_q u_\lambda u_i  +
2 g^{ij} \nabla_i |\nabla
u|^2 u_j + g^{ij} \nabla_i |N|^2 u_j \right) \nonumber\\
&&+ 2 e^u  \left(  \Delta u + 2 |\nabla
  u|^2 + |N|^2\right)  |\nabla u|^2 \nonumber\\
&=& e^u \big( \Delta |\nabla u|^2 - 2 |\nabla \nabla u|^2 - 2 R^{i \lambda} u_\lambda u_i  +
8 (\nabla \nabla u)^{ij} u_i u_j+ 2 g^{ij} \nabla_i |N|^2 u_j  \nonumber\\
&&+ 2 \Delta u |\nabla u|^2 + 4 |\nabla
  u|^4 + 2 |N|^2 |\nabla u|^2 \big)
\eea
The second term in (\ref{evol-nabla-u-1}) is
\bea
- g^{i a} \dot{g}_{ab} g^{bj} u_i u_j &=& e^u \bigg[ 2R^{ij} - 2 \nabla^i
\nabla^j u + 4 (N^2_-)^{ij} - 4 u_p (N^{ipj} + N^{jpi} ) \bigg] u_i u_j \nonumber\\
&&  - e^u |\nabla u|^4 + e^u [\omega^{ai} u_a u_i ] [ \omega^{bj}
  u_{b} u_j] \nonumber\\
&=& e^u \bigg[ 2R^{ij} - 2 (\nabla
\nabla u)^{ij} + 4 (N^2_-)^{ij} \bigg] u_i u_j - e^u |\nabla u|^4
\eea
The term $(N^{ipj}+N^{jpi})u_iu_j$
vanishes by $N^{ipj}=- N^{jip} - N^{pji}$ and $N^{jip} = -N^{pji} +
N^{ijp}$. Altogether, (\ref{evol-nabla-u-1})
becomes
\bea
(\p_t - e^u \Delta ) |\nabla u|^2 &=& e^u \bigg[- 2 |\nabla \nabla u|^2  + 2 g^{ij} \nabla_i |N|^2 u_j  + 2 \Delta u |\nabla u|^2 + 3 |\nabla
  u|^4 \nonumber\\
  &&+ 2 |N|^2 |\nabla u|^2  + 6 (\nabla^2 u)^{ij} u_i u_j + 4 (N^2_-)^{ij} u_i u_j \bigg]
\eea
The identity $N^2_- = 2 N^2_+ - {1 \over 4} |N|^2 g$ implies
\bea \label{evol-nabla-u-2}
(\p_t - e^u \Delta ) |\nabla u|^2 &=& e^u \bigg[- 2 |\nabla \nabla u|^2  + 2 g^{ij} \nabla_i |N|^2 u_j  + 2 \Delta u |\nabla u|^2 + 3 |\nabla
  u|^4 \nonumber\\
  &&+ |N|^2 |\nabla u|^2  + 6 (\nabla^2 u)^{ij} u_i u_j + 8
  (N^2_+)^{ij} u_i u_j  \bigg].
\eea

\subsubsection{The evolution of the Hessian of $u$}

We use as usual the notation $u_{ij} = \nabla_i \nabla_j u =
(\nabla^2 u)_{ij}$. The variation of the Hessian is
\be
\p_t u_{ip} = \nabla_i \nabla_p \dot{u} - \dot{\Gamma}^\lambda_{ip} u_\lambda.
\ee
Differentiating (\ref{diff-u-once})
\bea \label{diff-u-twice}
\nabla_p \nabla_i \dot{u} &=& e^u \bigg[ \nabla_p \nabla_i \Delta u +
  2 \nabla_p \nabla_i |\nabla
  u|^2 + \nabla_p \nabla_i |N|^2 \bigg]  + e^u \bigg[ \Delta u + 2 |\nabla u|^2 + |N|^2 \bigg] u_{ip} \nonumber\\
&&+ e^u \bigg[ \nabla_i \Delta u + 2 \nabla_i |\nabla
  u|^2 + \nabla_i |N|^2 \bigg] u_p + (i \leftrightarrow p)
\eea
Commuting derivatives, we see that
\bea
\nabla_p \nabla_i \Delta u &=& g^{ab} \nabla_p \nabla_a \nabla_i
\nabla_b u - g^{ab} \nabla_p (R_{ia}{}^\lambda{}_b u_\lambda)
\nonumber\\
&=& g^{ab} \nabla_a \nabla_p \nabla_i \nabla_b u - g^{ab}
R_{ap}{}^\lambda{}_i u_{\lambda b} - g^{ab} R_{ap}{}^\lambda{}_b u_{i
  \lambda} - g^{ab} \nabla_p (R_{ia}{}^\lambda{}_b u_\lambda)
\nonumber\\
&=& \Delta \nabla_p \nabla_i u - g^{ab} \nabla_a (R_{pb}{}^\lambda{}_i u_\lambda) - g^{ab}
R_{ap}{}^\lambda{}_i u_{\lambda b} - g^{ab} R_{ap}{}^\lambda{}_b u_{i \lambda} - g^{ab} \nabla_p (R_{ia}{}^\lambda{}_b u_\lambda) \nonumber
\eea
Therefore
\bea
\nabla_p \nabla_i \dot{u} &=& e^u \Delta u_{ip} + e^u \bigg[ 4g^{ab}
  u_{ia} u_{pb} \bigg] \nonumber\\
&&+ e^u \bigg[ \nabla Rm * \nabla u + \nabla^2 N * N + \nabla^3 u * \nabla u  \bigg] \nonumber\\
&&+ e^u \bigg[ \nabla^2 u *
  \mathcal{O}(\nabla u, Rm, N) + \nabla N * \nabla N + \nabla N
  * N * \nabla u \bigg] \nonumber
\eea
Here we used the identity $R = \Delta u + |N|^2$ on the term $e^u \Delta
u \, u_{ip}$.

Before further proceeding the computation, we explain the notations
used in the above formula. The terms written as $\alpha * \beta$
represent contractions of the tensors $\alpha$ and $\beta$ which are
linear in both $\alpha,\beta$. In later computation, we will also use $(\alpha+\gamma) * (\beta+ \eta)$ to represent the linear contractions among the tensors $\alpha, \beta, \gamma$ and $\eta$.
The notation $ \mathcal{O}(\nabla u, Rm, N)$ indicates terms which
only depend on $\nabla u, Rm$ and $N$ (but the dependence may be
nonlinear). We will soon prove a gradient estimate $|\nabla u| \leq
C$, so that $\mathcal{O}(\nabla u, Rm, N)$ will be treated as bounded terms.

\

Next,
\be
- \dot{\Gamma}^\lambda_{ip} u_\lambda = - {g^{\lambda \mu} \over 2}
(-\nabla_\mu \dot{g}_{ip} + \nabla_p \dot{g}_{\mu i} + \nabla_i \dot{g}_{\mu p}) u_\lambda
\ee
Since
\be
\nabla_\mu \dot{g}_{ij} = \nabla_\mu \bigg[ e^u (-2R_{ij} + 2 \nabla_i
  \nabla_j u - 4 (N^2_-)_{ij} + u_i u_j - u_{Ji} u_{Jj} + 4 u_p
  (N_i{}^p{}_j + N_j{}^p{}_i ) ) \bigg]
\ee
we get
\bea
- \dot{\Gamma}^\lambda_{ip} u_\lambda &=& e^u \bigg[ \nabla Rm *
  \nabla u + \nabla^3 u * \nabla u  \bigg] + \mathcal{O}(\nabla u, Rm, N) \nonumber\\
&&+ e^u \bigg[ \nabla^2 u *
  \mathcal{O}(\nabla u, Rm, N) + \nabla N
  * (N+\nabla u) * \nabla u \bigg]
\eea
Therefore
\bea \label{evol-D^2u}
(\p_t - e^u \Delta) u_{ip} &=& e^u \bigg[ 4g^{ab}
  u_{ia} u_{pb} \bigg] + \mathcal{O}(\nabla u, Rm, N) \nonumber\\
&&+ e^u \bigg[ \nabla Rm * \nabla u + \nabla^2 N * N + \nabla^3 u * \nabla u  \bigg] \nonumber\\
&&+ e^u \bigg[ \nabla^2 u *
  \mathcal{O}(\nabla u, Rm, N) + \nabla N * \nabla N + \nabla N
  * (N+\nabla u) * \nabla u \bigg] \nonumber
\eea
Next, we compute
\bea
(\p_t - e^u \Delta) |\nabla^2 u|^2 &=& 2 g^{ij} g^{pq} (\p_t - e^u
\Delta) u_{ip} u_{jq} - 2 e^u |\nabla^3 u|^2 \nonumber\\
&&- g^{ia} \dot{g}_{ab} g^{bj} g^{pq} u_{ip} u_{jq} - g^{ij} g^{pa} \dot{g}_{ab}
g^{bq} u_{ip} u_{jq}.
\eea
The first term is then
\bea
2 g^{ij} g^{pq} (\p_t-e^u\Delta) u_{ip} u_{jq} &=& e^u \bigg[ 8 g^{ab}
  u_{ia} u_{pb} (\nabla^2 u)^{ip} \bigg] + \mathcal{O}(\nabla u, Rm, N) *
\nabla^2 u \nonumber\\
&&+ e^u \bigg[ \nabla Rm * \nabla u + \nabla^2 N * N +
  \nabla^3 u * \nabla u  \bigg] *
\nabla^2 u \nonumber\\
&&+ e^u \bigg[ \nabla^2 u *
  \mathcal{O}(\nabla u, Rm, N) + \nabla N * \nabla N + \nabla N * (N+\nabla u) * \nabla u \bigg] * \nabla^2 u \nonumber
\eea
Since $\dot{g}_{ab}= 2 e^u u_{ab} + \mathcal{O}(\nabla u, Rm, N)$
\bea \label{evol-D^2u}
(\p_t - e^u \Delta) |\nabla^2 u|^2 &=& e^u \bigg[ 4 g^{ab}
  u_{ia} u_{pb} (\nabla^2 u)^{ip} \bigg] - 2e^{u}|\nabla^3 u|^2 + \mathcal{O}(\nabla u, Rm, N) *
\nabla^2 u \nonumber\\
&&+ e^u \bigg[ \nabla Rm * \nabla u + \nabla^2 N * N +
  \nabla^3 u * \nabla u  \bigg] *
\nabla^2 u \nonumber\\
&&+ e^u \bigg[ \nabla^2 u *
  \mathcal{O}(\nabla u, Rm, N) + \nabla N * \nabla N + \nabla N * (N+\nabla u) * \nabla u \bigg] * \nabla^2 u \nonumber
\eea

\subsection{The evolution of the Nijenhuis tensor: proof of Theorem \ref{th:Nijenhuis}(a)}

\subsubsection{Rewriting the flow of the complex structure}
The almost complex structure is given by $J^k{}_j = \omega^{ki}
g_{ij}$. Therefore, $\p_t J^k{}_j = \omega^{ki} \partial_t g_{ij}$. By substituting equation (\ref{ricciakg2}) for the Ricci
curvature $R_{ij}$ into the flow of metric $\p_t g_{ij}$ (\ref{evolution-longtime}), we obtain
\bea
\partial_t J^k{}_j &=&  e^u \omega^{ki} \bigg\{ 4 {\frak D}_p  N_{i}{}^{p}{}_{j} -  J^p{}_i J^q{}_j \nabla_q \nabla_p u + \nabla_i \nabla_j u \nonumber\\
&&  + u_i u_j  - u_{Ji} u_{Jj}  + 2 u_p N_j{}^p{}_i +6 u_p N_i{}^p{}_j  \bigg\}.
\eea
This simplifies to
\bea
\partial_t J^k{}_j &=&  e^u\bigg\{ 4 J^k{}_q {\frak D}_p  N^{qp}{}_{j} - J^q{}_j \nabla_q \nabla^k u + J^k{}_q \nabla^q \nabla_j u\nonumber\\
&&  + u^{Jk} u_j  - u^k u_{Jj}  +2 u_p N_j{}^{p,Jk} +6 u_p N^{Jk,p}{}_j  \bigg\}.
\eea
Converting covariant derivatives using $\nabla_\ell V^p = {\frak D}_\ell V^p + N_{\ell \lambda}{}^p V^\lambda$,
we obtain
\bea
 4 J^k{}_q {\frak D}_p  N^{qp}{}_{j}  &=& - 4 J^k{}_q {\frak D}_p N^{q}{}_{j}{}^p  \nonumber\\
&=& -4 J^k{}_q (\nabla_p N^{q}{}_{j}{}^p - N_{p \lambda}{}^q N^\lambda{}_j{}^p + N_{pj}{}^\lambda N^q{}_\lambda{}^p - N_{p \lambda}{}^p N^q{}_j{}^\lambda) \nonumber
 \eea
Using the symmetries of the Nijenhuis tensor, this is
 \bea \label{J-evol-id1}
 4 J^k{}_q {\frak D}_p  N^{qp}{}_{j} &=& -4 J^k{}_q \nabla_p N^{q}{}_{j}{}^p - 4 N_{p \lambda}{}^{Jk} N^{\lambda p}{}_j + 4 N^{p \lambda}{}_j N^{Jk}{}_{\lambda p} \nonumber\\
 &=& -4 J^k{}_q \nabla_p N^{q}{}_{j}{}^p - 4 N_{p \lambda}{}^{Jk} N^{\lambda p}{}_j + 4 N^{p \lambda}{}_j(- N_p{}^{Jk}{}_{\lambda}- N_{\lambda p}{}^{Jk}) \nonumber\\
 &=& -4 J^k{}_q \nabla_p N^{q}{}_{j}{}^p - 8 (N^2_-)^{Jk}{}_j + 4 (N^2_+)^{Jk}{}_j
 \eea
 using $(N^2_+)_{ij} = N^{p \lambda}{}_i N_{p \lambda j}$ and $(N^2_-)_{ij} = N^{p \lambda}{}_i N_{\lambda p j}$. Thus
\bea \label{evol-J-1}
\partial_t J^k{}_j &=&  e^u \bigg\{ -4 J^k{}_q \nabla_p N^{q}{}_{j}{}^p - J^q{}_j \nabla_q \nabla^k u  + J^k{}_q \nabla^q \nabla_j u \nonumber\\
&&  + u^{Jk} u_j  - u^k u_{Jj}  +2 u_p N_j{}^{p,Jk} + 6 u_p N^{Jk,p}{}_j - 8 (N^2_-)^{Jk}{}_j + 4 (N^2_+)^{Jk}{}_j \bigg\}.
\eea

\subsubsection{A first formulation of the evolution of the Nijenhuis tensor}
We start with the identity
\be \label{nabla-J}
\nabla_i J^k{}_j = - 2 N_{ij}{}^{Jk},
\ee
which follows from the formula relating $\nabla$ to ${\frak D}$ and
$N_{i,Jj}{}^k=-N_{ij}{}^{Jk}$. Indeed,
\be
\nabla_i J^k{}_j = {\frak D}_i J^k{}_j + N_{i \lambda}{}^k
J^\lambda{}_j - J^k{}_\lambda N_{ij}{}^\lambda = - 2 N_{ij}{}^{Jk}.
\ee
We can expand (\ref{nabla-J}) and obtain
\be
J^k{}_p N_{ij}{}^p = -{1 \over 2} \nabla_i J^k{}_j = -{1 \over 2} \partial_i J^k{}_j - {1 \over 2} (\Gamma^k_{i\lambda} J^\lambda{}_j - J^k{}_\lambda \Gamma^\lambda_{ij}).
\ee
Differentiating this gives
\be
\dot{J}^k{}_p N_{ij}{}^p + J^k{}_p \dot{N}_{ij}{}^p = -{1 \over 2} \nabla_i \dot{J}^k{}_j - {1 \over 2} (\dot{\Gamma}^k_{i\lambda} J^\lambda{}_j - J^k{}_\lambda \dot{\Gamma}^\lambda_{ij}),
\ee
which leads to
\be
\partial_t N_{ij}{}^\ell = {1 \over 2}J^\ell{}_k  \nabla_i \dot{J}^k{}_j   + J^\ell{}_k  \dot{J}^k{}_p N_{ij}{}^p + {1 \over 2} (J^\ell{}_k \dot{\Gamma}^k_{i\lambda} J^\lambda{}_j + \dot{\Gamma}^\ell_{ij}) .
\ee
We will introduce some notation to group terms. We first introduce the
tensor $Z$ given by
\be
Z_{ij}{}^{J\ell} = J^\ell{}_r \dot{\Gamma}^r_{in} J^n{}_j + \dot{\Gamma}^\ell{}_{ij}.
\ee
Next, we denote $ \dot{J}^k{}_j = e^u E^k{}_j$, where by (\ref{evol-J-1}),
\bea
E^k{}_j &=&   -4 J^k{}_q \nabla_p N^{q}{}_{j}{}^p - J^q{}_j \nabla_q \nabla^k u + J^k{}_q \nabla^q \nabla_j u\nonumber\\
&&  + J^k{}_p u^{p} u_j  -  J^p{}_j u^k u_{p}  + 2 J^k{}_\ell u_p N_j{}^{p \ell} + 6 J^k{}_\ell u_p N^{\ell p}{}_j \nonumber\\
&&- 8 J^k{}_\ell (N^2_-)^{\ell}{}_j + 4 J^k{}_\ell (N^2_+)^{\ell}{}_j .
\eea
We write
\bea
\nabla_i E^k{}_j &=&  -4 J^k{}_q \nabla_i \nabla_p N^{q}{}_{j}{}^p -
J^q{}_j \nabla_i \nabla_q \nabla^k u + J^k{}_q \nabla_i \nabla^q
\nabla_j u \nonumber\\
&&  + J^k{}_p \nabla_i (u^{p} u_j)  -  J^p{}_j \nabla_i(u^k u_{p}) + Y_i{}^k{}_j,
\eea
where
\bea
Y_r{}^k{}_j &=& -4 \nabla_r J^k{}_q \nabla_p N^{q}{}_{j}{}^p - \nabla_r J^q{}_j \nabla_q \nabla^k u + \nabla_r J^k{}_q \nabla^q \nabla_j u \nonumber\\
&&  + \nabla_r J^k{}_p u^{p} u_j  -  \nabla_r J^p{}_j u^k u_{p}  + 2\nabla_r ( J^k{}_\ell u_p N_j{}^{p \ell}) + 6\nabla_r( J^k{}_\ell u_p N^{\ell p}{}_j) \nonumber\\
&&- 8\nabla_r( J^k{}_\ell (N^2_-)^{\ell}{}_j) + 4 \nabla_r(J^k{}_\ell (N^2_+)^{\ell}{}_j)
\eea
Therefore
\bea
J^\ell{}_k \nabla_i \dot{J}^k{}_j &=& J^\ell{}_k \nabla_i e^u E^k{}_j + J^\ell{}_k e^u  \nabla_i  E^k{}_j \nonumber\\
&=& e^u \bigg[  4 \nabla_i \nabla_p N^{\ell}{}_{j}{}^p -  J^\ell{}_kJ^q{}_j \nabla_i \nabla_q \nabla^k u - \nabla_i \nabla^\ell \nabla_j u\nonumber\\
&&  - \nabla_i (u^{\ell} u_j)  -   J^\ell{}_kJ^p{}_j \nabla_i(u^k
  u_{p}) + J^\ell{}_k Y_i{}^k{}_j + u_i E^{J \ell}{}_j \bigg].
\eea
Substituting this,
\bea
\partial_t N_{ij}{}^\ell &=&  e^u \bigg[  -2 \nabla_i
  \nabla_p N^{\ell p}{}_{j} - {1 \over 2} J^\ell{}_kJ^q{}_j \nabla_i \nabla_q \nabla^k u - {1 \over 2} \nabla_i \nabla^\ell \nabla_j u \nonumber\\
&&  -{1 \over 2} \nabla_i (u^{\ell} u_j)  -  {1 \over 2} J^\ell{}_kJ^p{}_j \nabla_i(u^k u_{p}) +{1 \over 2} J^\ell{}_k Y_i{}^k{}_j \nonumber\\
&&+ {1 \over 2} u_i E^{J \ell}{}_j + J^\ell{}_k E^k{}_p N_{ij}{}^p
  \bigg] + {1 \over 2} Z_{ij}{}^{J \ell}.
\eea
To interpret the highest order terms, we will need the following
identity. We claim:
\bea \label{Delta-N-lemma}
\Delta N_{i j \ell} &=& - 2 \nabla_i \na_p N_{\ell}{}^p{}_j  - \nabla_i
R_{\ell j}- {1 \over 2} \nabla^p (R_{pi j \ell} - R_{p,i, Jj,J \ell}) \nonumber\\
&&+ [\nabla_i,\nabla_p] N^p{}_{\ell j} + {1 \over 2} \nabla_i
\nabla_\ell  \nabla_j u + {1 \over 2} \nabla_i(
J^p{}_\ell J^q{}_j \nabla_{p} \nabla_{q}u) \nonumber\\
&& +2 \bigg[ \nabla^p ( N_{pj}{}^{r} N_{i \ell r}) -  \nabla^p (
  N_{ij}{}^{r} N_{p \ell r}) +\nabla_i (N_-^2)_{\ell j} - \nabla_i
  (N_+^2)_{\ell j} \bigg].
\eea
We assume identity (\ref{Delta-N-lemma}) for now and give the proof in \S \ref{subsubsection-lemma-proofs}. The evolution of $N$ becomes
\bea
\partial_t N_{ij}{}^\ell &=& e^u \bigg[ \Delta
  N_{ij}{}^\ell - \nabla_i \nabla^\ell
  \nabla_j u  + \nabla_i
R_{\ell j}+ {1 \over 2} \nabla^p (R_{pi j \ell} - R_{p,i, Jj,J \ell})
\nonumber\\
&&- {1 \over 2} (u_i u_j u^\ell + u_i u_{Jj} u^{J \ell})+ ({\rm
  IIa})_{ij}{}^\ell + {e^{-u} \over 2} Z_{ij}{}^{J \ell} \nonumber\\
&&+ Rm * N + \nabla N * (N+ \nabla u) +N^3 + N^2 * \nabla u + N * (\nabla u)^2 \bigg]
\eea
where terms involving $\nabla \nabla u$ will need to be tracked for future
use, and are given explicitly by
\bea
({\rm IIa})_{ij}{}^\ell &=&  -{1 \over 2} \nabla_i (u^{\ell} u_j)  -  {1
  \over 2} J^\ell{}_kJ^p{}_j \nabla_i(u^k u_{p}) - {1 \over
  2} g^{r \ell} \nabla_i (J^p{}_r J^q{}_j) (\nabla^2 u)_{pq}
\nonumber\\
&&+ {1 \over 2} J^\ell{}_k \bigg[ - \nabla_i J^q{}_j (\nabla^2
  u)_q{}^k + \nabla_i J^k{}_q (\nabla^2 u)^q{}_j + 2 J^k{}_r
  (\nabla^2 u)_{ip} N_j{}^{p r} + 6 J^k{}_r (\nabla^2 u)_{ip} N^{r
    p}{}_j \bigg] \nonumber\\
&&+{1 \over 2} u_i J^\ell{}_k \bigg[ - J^q{}_j (\nabla^2 u)^k{}_q + J^k{}_q (\nabla^2 u)^q{}_j \bigg] \nonumber\\
&&+ J^\ell{}_k \bigg[ - J^q{}_p (\nabla^2 u)_q{}^k + J^k{}_q (\nabla^2
  u)^q{}_p \bigg] N_{ij}{}^p
\eea
which, using $\nabla_i J^k{}_j = - 2 N_{ij}{}^{Jk}$ and simplifying, become
\bea \label{N-IIa}
({\rm IIa})_{ij}{}^\ell &=&  -{1 \over 2} \nabla_i (u^{\ell} u_j)  -  {1
  \over 2} J^\ell{}_kJ^p{}_j \nabla_i(u^k u_{p}) + N_{i}{}^{\ell Jp}
(\nabla^2 u)_{p,Jj} -N_{ij}{}^{Jq} (\nabla^2 u)^{J \ell}{}_{q}
\nonumber\\
&&+ N_{ij}{}^{Jq} (\nabla^2
  u)_q{}^{J \ell} +  N_{iq}{}^{\ell} (\nabla^2 u)^q{}_j -
  (\nabla^2 u)_{ip} N_j{}^{p \ell} - 3  (\nabla^2 u)_{ip} N^{\ell
    p}{}_j \nonumber\\
&&-{1 \over 2} u_i \bigg[ (\nabla^2 u)^{J \ell}{}_{Jj} +(\nabla^2
  u)^\ell{}_j \bigg] - \bigg[ (\nabla^2 u)_{Jp}{}^{J \ell} +(\nabla^2
  u)^\ell{}_p \bigg] N_{ij}{}^p.
\eea
Next, we claim that
\bea \label{lemma-Z-terms}
Z_{ij}{}^{J p} &=&  \nabla_i \dot{g}^p{}_{j}   + {1 \over 2} (-\nabla^p \dot{g}_{ij}  + \nabla_j \dot{g}^p{}_{i})+{1 \over 2} (\omega^{rp} J^n{}_j -  \omega^{np} J^r{}_j) \nabla_r \dot{g}_{in}  \nonumber\\
&&+ (N_{i}{}^p{}^{r} \dot{g}_{jr} +   N_{ij}{}^{r} \dot{g}^p{}_{r}).
\eea
We assume identity (\ref{lemma-Z-terms}) for now and give the proof later in \S
\ref{subsubsection-lemma-proofs}. Substituting the evolution of $g_{ij}$
(\ref{evolution-longtime}) into this expression for $Z_{ij}{}^{Jp}$
and then in our expression for $\p_t N_{ij}{}^\ell$, we obtain
\bea \label{p-t-N-ijl}
\partial_t N_{ij}{}^\ell &=& e^u \bigg[ \Delta
  N_{ij}{}^\ell - \nabla_i \nabla^\ell
  \nabla_j u + \nabla_i
R^\ell{}_{j}+ {1 \over 2} \nabla^p (R_{pi j}{}^\ell - g^{\ell r} R_{p,i, Jj,J r})
\nonumber\\
&&- \nabla_i R^\ell{}_{j}   - {1 \over 2} (-\nabla^\ell
R_{ij}  + \nabla_j R^\ell{}_{i}) - {1 \over 2}
(\omega^{r\ell} J^n{}_j -  \omega^{n\ell} J^r{}_j) \nabla_r
R_{in} \nonumber\\
&&+ \nabla_i \nabla^\ell \nabla_j u   + {1 \over 2} (-\nabla^\ell
\nabla_i \nabla_j u  +\nabla_j \nabla^\ell \nabla_i u) + {1 \over 2}
(\omega^{r\ell} J^n{}_j -  \omega^{n\ell} J^r{}_j) \nabla_r \nabla_i
\nabla_n u\nonumber\\
&&+ Ric * \nabla u + ({\rm IIa})_{ij}{}^\ell + ({\rm IIb})_{ij}{}^\ell\nonumber\\
&&+ Rm * N + \nabla N * (N+ \nabla u)+ N^3 + N^2 * \nabla u + N * (\nabla u)^2 \bigg]
\eea
where terms of order $(\nabla u)^3$ (e.g. $u_i u_j u^\ell$) have
cancelled and the additional terms involving $\nabla^2 u$ are
\bea \label{N-IIb}
({\rm IIb})_{ij}{}^\ell &=& {1 \over 2} [\nabla_i (u^\ell u_j - u^{J \ell}
u_{Jj} ) + 4 (\nabla^2 u)_{pi} (N^{\ell p}{}_j + N_j{}^{p \ell}) + 2 u_i (\nabla^2 u)^\ell{}_j] \nonumber\\
&&-  {1 \over 4} [\nabla^\ell(u_i u_j - u_{Ji} u_{Jj}) + 4 (\nabla^2
  u)_p{}^\ell (N_i{}^p{}_j + N_j{}^p{}_i) + 2 u^\ell (\nabla^2 u)_{ij}]  \nonumber\\
&&+ {1 \over 4}[\nabla_j(u_i u^\ell - u_{Ji} u^{J \ell}) + 4 (\nabla^2
  u)_{pj} (N_i{}^{p \ell} + N^{\ell p}{}_i) + 2 u_j (\nabla^2 u)_{i}{}^\ell ] \nonumber\\
  &&+{1 \over 4} (\omega^{r \ell} J^n{}_j -  \omega^{n \ell} J^r{}_j)
    [\nabla_r (u_i u_n - u_{Ji} u_{Jj}) + 4 (\nabla^2
      u)_{pr}(N_i{}^p{}_n + N_n{}^p{}_i) + 2 u_r(\nabla^2 u)_{in} ] \nonumber\\
&&+ (N_{i}{}^\ell{}^{r} (\nabla^2 u)_{jr}+   N_{ij}{}^{r} (\nabla^2 u)^\ell{}_r).
\eea
Since we can commute $\nabla_\ell \nabla_j \nabla_i u = \nabla_j
\nabla_\ell \nabla_i u - R_{\ell j}{}^p{}_i u_p$, the terms of order
$\nabla^3 u$ in (\ref{p-t-N-ijl}) cancel. We are left with
\bea\label{1.43}
\partial_t N_{ij}{}^\ell &=& e^u \bigg[ \Delta
  N_{ij}{}^\ell + {1 \over 2} \nabla^p R_{pi j}{}^\ell - {1 \over 2}
  \omega^{n \ell} J^r{}_j \nabla^p R_{pi r n}
\nonumber\\
&&  - {1 \over 2} (-\nabla^\ell
R_{ij}  + \nabla_j R^\ell{}_{i}) - {1 \over 2}
(\omega^{r\ell} J^n{}_j -  \omega^{n\ell} J^r{}_j) \nabla_r
R_{in} \nonumber\\
&&+ {1 \over 2} (-R_j{}^{\ell p}{}_i u_p + R^{J \ell}{}_i{}^p{}_{Jj}
u_p - R_{Jj,i}{}^{p,J\ell} u_p)+ Ric * \nabla u+ ({\rm
  IIa})_{ij}{}^\ell + ({\rm IIb})_{ij}{}^\ell \nonumber\\
&& + Rm * N + \nabla N * (N+ \nabla u)+ N^3 + N^2 *
\nabla u + N * (\nabla u)^2 \bigg] .
\eea

The terms of order $\nabla Rm$ also cancel. Indeed, the Bianchi
identity is
\bea\nonumber
\nabla^p R_j{}^\ell{}_{pi} + \nabla^\ell R^p{}_{jpi} +
\nabla_j R^{\ell p}{}_{pi} = 0
\eea
and hence
\be
{1 \over 2} \nabla^p R_{pi j}{}^\ell = {1 \over 2} \nabla^p
R_{j}{}^\ell{}_{pi} = {1 \over 2} (-\nabla^\ell R_{ij} + \nabla_j R^\ell{}_i).
\ee
For the terms
involving $\omega,J$, the same argument gives
\be
- {1 \over 2}  \omega^{n \ell} J^r{}_j \nabla^p R_{pi r n} = -{1 \over
   2}  \omega^{n \ell} J^r{}_j \nabla^p R_{nr ip} =   {1 \over
   2}  \omega^{n \ell} J^r{}_j (\nabla_n R_{ri} - \nabla_r R_{in})
\ee
This is the same thing as
\be
- {1 \over 2}  \omega^{n \ell} J^r{}_j \nabla^p R_{pi r n} = {1 \over
  2}  (\omega^{r \ell} J^n{}_j - \omega^{n \ell} J^r{}_j) \nabla_r R_{in}
\ee
Putting these identities back to (\ref{1.43}), we get the cancellation of the $\nabla Rm$ terms. Thus
\bea \label{evol-Nijk}
\partial_t N_{ij}{}^k &=& e^u \bigg[ \Delta N_{ij}{}^k + (Rm * \nabla u)_{ij}{}^k + ({\rm
  IIa})_{ij}{}^k + ({\rm IIb})_{ij}{}^k \nonumber\\
&&+ Rm * N + \nabla N * (N+ \nabla u)+ N^3 + N^2 *
\nabla u + N * (\nabla u)^2 \bigg],
\eea
where the $({\rm II})$ terms involve $\nabla^2 u$ and are explicitly given in (\ref{N-IIa})
and (\ref{N-IIb}) and $Rm * \nabla u$ is of the form
\be \label{RMstarDu}
(Rm * \nabla u)_{ij}{}^k =  {1 \over 2} (-R_j{}^{k p}{}_i u_p + R^{J k}{}_i{}^p{}_{Jj}
u_p - R_{Jj,i}{}^{p,Jk} u_p)+ (Ric * \nabla u)_{ij}{}^k.
\ee

\subsubsection{Proof of Identity (\ref{Delta-N-lemma}) and Identity
  (\ref{lemma-Z-terms})} \label{subsubsection-lemma-proofs}

\bigskip
\par {\it Proof of identity (\ref{Delta-N-lemma}):} The starting point is
the identity (\ref{Rij-JkJl}) for the action of $J$ on the Riemann curvature
tensor. Recall that in the case $d \omega = 0$, (\ref{christ}) and
(\ref{torsion}) specialize to $A=N=\mathfrak T$, and so the identity (\ref{Rij-JkJl}) becomes
\bea \label{Rij-JkJl2}
R_{i,j,Jk,Jl}-R_{ijkl} &=& 2(\mathfrak D_i N_{jkl}-\mathfrak D_j N_{ikl}+
N^r{}_{ij} N_{r kl}) \nonumber\\
&=& 2(\mathfrak D_i N_{jkl}-\mathfrak D_j N_{ikl}+ N_{ji}{}^r N_{r kl} - N_{ij}{}^r N_{r kl}),
\eea
using $N_{ijk}+N_{kij}+N_{jki}=0$. We can convert ${\frak D} N$ to $\nabla N$. For example,
\be
{\mathfrak D}_i N_{jkl} = \nabla_i N_{jkl} + N_{ij}{}^r
N_{r k l} + N_{ik}{}^r N_{j r l} + N_{i l}{}^r
N_{jk r}.
\ee
After converting ${\frak D}$ to $\nabla$, (\ref{Rij-JkJl2}) becomes
\be
\na_j N_{i \ell k} =  \na_i N_{j \ell k} + {1 \over 2} (R_{ji k \ell} - R_{j,i, Jk,J \ell}) -2 N_{jk}{}^{r} N_{i \ell r} + 2 N_{ik}{}^{r} N_{j \ell r}.
\ee
Differentiating this identity, we obtain
\bea
\nabla_q \na_p N_{k i j} &=&  \nabla_q \na_k N_{p i j} + {1 \over 2} \nabla_q (R_{pk j i} - R_{p,k, Jj,J i}) \nonumber\\
&&-2 \nabla_q( N_{pj}{}^{r} N_{k i r}) + 2 \nabla_q( N_{kj}{}^{r} N_{p i r}).
\eea
By the Bianchi identity,
\bea
\Delta N_{ijk} &=& - \Delta N_{kij} + \Delta N_{jik} \nonumber\\
&=&  -\nabla_p \na_k N^p{}_{i j} - {1 \over 2} \nabla^p (R_{pk j i} - R_{p,k, Jj,J i}) \nonumber\\
&&+2 \nabla^p ( N_{pj}{}^{r} N_{k i r}) - 2 \nabla^p ( N_{kj}{}^{r} N_{p i r}) + \Delta N_{jik}.
\eea
Therefore
\bea \label{DeltaN-1}
-\Delta N_{jik} + \Delta N_{ijk} &=&  -\nabla_k \na_p N^p{}_{i j} - {1 \over 2} \nabla^p (R_{pk j i} - R_{p,k, Jj,J i}) \nonumber\\
&&+ [\nabla_k,\nabla_p] N^p{}_{ij} +2 \nabla^p ( N_{pj}{}^{r} N_{k i r}) - 2 \nabla^p ( N_{kj}{}^{r} N_{p i r}) \nonumber\\
 &=&  \nabla_k \na_p N_j{}^p{}_{i} - \nabla_k \na_p N_{i}{}^p{}_j - {1 \over 2} \nabla^p (R_{pk j i} - R_{p,k, Jj,J i}) \nonumber\\
&&+ [\nabla_k,\nabla_p] N^p{}_{ij} +2 \nabla^p ( N_{pj}{}^{r} N_{k i r}) - 2 \nabla^p ( N_{kj}{}^{r} N_{p i r})
\eea
The formula for Ricci curvature in Type IIA geometry given in
(\ref{ricciakgnabla}) is
\be
\nabla_p N_{j}{}^p{}_i
=  -\nabla_p N_{i}{}^p{}_j - R_{ij} + 2(N^2_-)_{ij} - 2(N^2_+)_{ij} + \frac{1}{2}(\nabla^2u)_{ij}
+ \frac{1}{2}(\nabla^2u)_{Ji,Jj}.
\ee
Substituting this into (\ref{DeltaN-1}),
\bea
-\Delta N_{jik} + \Delta N_{ijk} &=& - 2 \nabla_k \na_p N_{i}{}^p{}_j  - \nabla_k R_{ij}- {1 \over 2} \nabla^p (R_{pk j i} - R_{p,k, Jj,J i}) \nonumber\\
&&+ [\nabla_k,\nabla_p] N^p{}_{ij} + {1 \over 2} \nabla_k \nabla_i
\nabla_j u+ {1 \over 2} \nabla_k( J^p{}_iJ^q{}_j
\nabla_{p} \nabla_{q} u) \nonumber\\
&& +2 \nabla^p ( N_{pj}{}^{r} N_{k i r}) - 2 \nabla^p ( N_{kj}{}^{r} N_{p i r}) +2 \nabla_k (N_-^2)_{ij} - 2 \nabla_k (N_+^2)_{ij}
\eea
This proves the identity after using $N_{ijk} + N_{kij} + N_{jki} =0$. Q.E.D.
\bigskip

\par {\it Proof of Identity (\ref{lemma-Z-terms}):} The variation of
the Christoffel symbol is given by
\be
\dot{\Gamma}^p_{in} = {g^{ps} \over 2} (-\nabla_s \dot{g}_{in} + \nabla_i \dot{g}_{sn} + \nabla_n \dot{g}_{is}).
\ee
Thus
\bea
Z_{ij}{}^p = \dot{\Gamma}^p_{in} J^n{}_j - J^p{}_n \dot{\Gamma}^n{}_{ij} &=& {g^{ps} \over 2} (-\nabla_s \dot{g}_{in}J^n{}_j + \nabla_i \dot{g}_{sn}J^n{}_j + J^n{}_j\nabla_n \dot{g}_{is}) \nonumber\\
&&- {g^{ns} \over 2} J^p{}_n(-\nabla_s \dot{g}_{ij} + \nabla_i \dot{g}_{sj} + \nabla_j \dot{g}_{is}).
\eea
Therefore
\bea
Z_{ij}{}^{J\ell} &=& {g^{ps} \over 2} J^\ell{}_p (-\nabla_s \dot{g}_{in}J^n{}_j + \nabla_i \dot{g}_{sn}J^n{}_j + J^n{}_j\nabla_n \dot{g}_{is}) \nonumber\\
&&- {g^{ns} \over 2} J^\ell{}_p J^p{}_n(-\nabla_s \dot{g}_{ij} + \nabla_i \dot{g}_{sj} + \nabla_j \dot{g}_{is}),
\eea
which becomes
\bea
Z_{ij}{}^{J p} &=& -{g^{p s} \over 2}  J^r{}_s \nabla_i \dot{g}_{rn}J^n{}_j + {1 \over 2} (-\nabla^p \dot{g}_{ij} + \nabla_i \dot{g}^p{}_{j} + \nabla_j \dot{g}^p{}_{i}) \nonumber\\
&&  +{1 \over 2} (\omega^{rp} J^n{}_j \nabla_r \dot{g}_{in} -  \omega^{np} J^r{}_j \nabla_r \dot{g}_{in}  ).
\eea
From the evolution equation of $g$, we have the identity $J^r{}_i \dot{g}_{rn} J^n{}_j = - \dot{g}_{ij}$. Therefore the first term
\be
-{g^{p s} \over 2}  J^r{}_s \nabla_i \dot{g}_{rn}J^n{}_j = - {g^{ps} \over 2} \nabla_i (J^r{}_s \dot{g}_{rn} J^n{}_j) +{g^{ps} \over 2} \nabla_i J^r{}_s \dot{g}_{rn} J^n{}_j +{g^{ps} \over 2} J^r{}_s \dot{g}_{rn} \nabla_i J^n{}_j
\ee
simplifies to
\be
-{g^{p s} \over 2}  J^r{}_s \nabla_i \dot{g}_{rn}J^n{}_j =  {1 \over 2} \nabla_i \dot{g}^p{}_{j}  + N_{i}{}^p{}^{r} \dot{g}_{rj} +  \dot{g}^p{}_{n} N_{ij}{}^{n} .
\ee
Therefore
\bea
Z_{ij}{}^{J p} &=&  \nabla_i \dot{g}^p{}_{j}   + {1 \over 2} (-\nabla^p \dot{g}_{ij}  + \nabla_j \dot{g}^p{}_{i}) +{1 \over 2} (\omega^{rp} J^n{}_j \nabla_r \dot{g}_{in} -  \omega^{np} J^r{}_j \nabla_r \dot{g}_{in}  ) \nonumber\\
&&+ N_{i}{}^p{}^{r} \dot{g}_{rj} +  \dot{g}^p{}_{n} N_{ij}{}^{n}.
\eea
Q.E.D.

\subsubsection{The evolution of the norm of $N$}
The norm of the Nijenhuis tensor, which is $
|N|^2 = g^{ij} g^{pq} g^{k \ell} N_{ipk} N_{jq \ell}$, evolves by
\bea
(\p_t - e^u \Delta) |N|^2 &=& 2 N^{ij}{}_k (\p_t - e^u \Delta) N_{ij}{}^k
 - 2 e^u |\nabla N|^2 \nonumber\\
&&- g^{ir} \dot{g}_{rs} g^{sj} N_{ipk} N_{j}{}^{pk} -2 g^{k r} \dot{g}_{rs} g^{s \ell} N_{ipk} N^{ip}{}_{\ell}.
\eea
By the equation for the evolution of $g_{ij}$ (\ref{evolution-longtime}), this is of the form
\bea
(\p_t - e^u \Delta) |N|^2 &=& 2 N^{ij}{}_k (\p_t - e^u \Delta)
N_{ij}{}^k - 2 e^u |\nabla N|^2 \nonumber\\
&&+ e^u \bigg[- 2 (\nabla^2 u)^{ij} N_{ipk} N_{j}{}^{pk} -4 (\nabla^2 u)^{k \ell} N_{ipk} N^{ip}{}_{\ell} \nonumber\\
&&+ N^2 * (Rm +N * N+ \nabla u * \nabla u + \nabla u
  * N) \bigg]
\eea
Substituting (\ref{evol-Nijk}), we obtain
\bea
& \ & (\p_t - e^u \Delta) |N|^2 \nonumber\\
&=& e^u \bigg[- 2 |\nabla N|^2 + 2 N^{ij}{}_k ({\rm IIa} +{\rm
  IIb})_{ij}{}^k - 2 (\nabla^2 u)^{ij} N_{ipk} N_{j}{}^{pk} -4 (\nabla^2 u)^{k \ell} N_{ipk} N^{ip}{}_{\ell} \nonumber\\
  &&+ 2 N^{ij}{}_k (Rm * \nabla u)_{ij}{}^k + Rm * N^2 \nonumber\\
  &&+ \nabla N * N * (N+\nabla u) + N^4 + N^3 * \nabla u + N^2 *
  (\nabla u)^2 \bigg].
\eea
The expressions for the terms $({\rm IIa} +{\rm IIb})$ are given in (\ref{N-IIa})
and (\ref{N-IIb}). A calculation can be done to verify that each term $N^{ij}{}_k ({\rm IIa})_{ij}{}^k$ and $N^{ij}{}_k ({\rm
  IIb})_{ij}{}^k$ only contributes terms of the type $(\nabla^2
u) * N^2$ since the others vanish by symmetry.
For example, if we
denote the terms on each line of $({\rm IIb})$ (\ref{N-IIb}) by
$(i)+(ii)+(iii)+(iv)+(v)$, we have
\bea
(i) &=& {1 \over 2} N^{ij}{}_\ell [\nabla_i (u^\ell u_j)  -
  J^\ell{}_r J^s{}_j \nabla_i (u^{r}
u_{s} ) + 4 u_{pi} (N^{\ell p}{}_j + N_j{}^{p \ell}) + 2 u_i
u^\ell{}_j] \nonumber\\
&=&-2  u_{pi} N^{ij}{}_\ell N^{\ell}{}_{jp} - 2 u_{pi} N^{ij}{}_\ell
N_j{}^\ell{}_p \nonumber\\
&=& 2  u^p{}_{i} (N^{\ell ij}+N^{j \ell i}) N_{\ell jp} + 2
u_{pi} (N^{\ell ij} + N^{j \ell i})
N_{j \ell p} =0
\eea
where we used the symmetry $N^{ij}{}_\ell = -
N^{i}{}_\ell{}^j$, the identity $N^{i,Jj}{}_{Jr} = N^{ij}{}_r$, and
the Bianchi identity $N_{ijk} + N_{kij}+N_{jki} = 0$. The
symmetry $N^{ij}{}_\ell = -
N^{i}{}_\ell{}^j$ allows us to combine the $(ii)+(iii)$ terms:
\bea
(ii)+(iii) &=& -{1 \over 2} N^{ij}{}_\ell [\nabla^\ell(u_i u_j) - J^r{}_i J^s{}_j \nabla^\ell (u_{r} u_{s}) + 4 u_p{}^\ell (N_i{}^p{}_j + N_j{}^p{}_i) + 2 u^\ell u_{ij}] \nonumber\\
&=& -{1 \over 2} N^{ij}{}_\ell [ 2u^\ell{}_i u_j + 4 u_p{}^\ell
  (N_i{}^p{}_j + N_j{}^p{}_i) + 2 u^\ell u_{ij}] \nonumber\\
&=&  2 u_p{}^\ell (N^{ij}{}_\ell N_{ij}{}^p) + 2 u_p{}^\ell
(N^{ij}{}_\ell N_{ji}{}^p )
\eea
where we used $N^{Ji,Jj}{}_\ell = - N^{ij}{}_\ell$. Next,
\bea
(iv) &=& {1 \over 2} N^{ij}{}_\ell \omega^{r \ell} J^n{}_j [\nabla_r (u_i u_n) - J^p{}_i J^q{}_n \nabla_r( u_{p} u_{q}) + 4 u_{pr}(N_i{}^p{}_n + N_n{}^p{}_i) + 2 u_r u_{in} ]
\nonumber\\
&=& N^{ij}{}_\ell \omega^{r \ell} J^n{}_j [ u_{ri} u_n +
  u_i u_{rn} + 2 u_{pr}(N_i{}^p{}_n + N_n{}^p{}_i) +  u_r u_{in} ]
\nonumber\\
&=& - N^{i,Jn}{}_{J\ell} g^{r \ell} [ u_{ri} u_n +
  u_i u_{rn} + 2 u_{pr}(N_i{}^p{}_n + N_n{}^p{}_i) +  u_r u_{in} ]
\nonumber\\
&=& -N^{in}{}_\ell [ u^\ell{}_{i} u_n +
  u_i u^\ell{}_{n} + 2 u^\ell{}_{p}(N_i{}^p{}_n + N_n{}^p{}_i) +
  u^\ell u_{in} ] \nonumber\\
&=&  2 u^\ell{}_p (N^{in}{}_\ell N_{in}{}^p) + 2 u^\ell{}_p
(N^{in}{}_\ell N_{ni}{}^p).
\eea
The computations of the other terms is similar. The result is then
\bea
(\p_t - e^u\Delta) |N|^2 &=& e^u \bigg[ - 2 |\nabla N|^2 + (\nabla^2
  u) * N^2 + 2 N^{ij}{}_k (Rm * \nabla u)_{ij}{}^k \\
  && + Rm * N^2 + \nabla N * N * (N+\nabla u) + N^4 + N^3 * \nabla
  u + N^2 * (\nabla u)^2 \bigg] \nonumber
\eea
where $(\nabla^2 u) \star N^2$ is an expression involving $(\nabla
u)^{ij} (N^2_+)_{ij}$ and $(\nabla
u)^{ij} (N^2_-)_{ij}$. The expression for the term $ (Rm * \nabla
u)_{ij}{}^k$ is given in (\ref{RMstarDu}), and using $N^{i,Jj}{}_{Jk}
= N^{ij}{}_k$ and symmetries of the curvature tensor, it becomes
\be
N^{ij}{}_k (Rm * \nabla u)_{ij}{}^k =  -R^{p}{}_{ijk}
u_p N^{ijk} + N^{ij}{}_k(Ric * \nabla u)_{ij}{}^k.
\ee
By (\ref{curvrel2}), we can convert $R^p{}_{ijk} = {\frak
  R}^p{}_{ijk} + \nabla N + N^2$. Since ${\frak
  R}^p{}_{i,Jj,Jk} = {\frak R}^p{}_{ijk}$ and $N^{i,Jj,Jk} = -
N^{ijk}$, the term ${\frak R}^p_{ijk} N^{ijk} = 0$ by
symmetry. Similarly, we can reduce terms of the form $ N^{ij}{}_k u^k R_{ij} = N^{ij}{}_k
u^k R^{-J}_{ij}$ where $R^{-J}_{ij}$ is the $J$-anti-invariant part of the Ricci tensor
given in (\ref{jinvricci2}) by $Ric^{-J} = \nabla N + N*N$. Absorbing
these terms gives the expression
\bea
(\p_t - e^u\Delta) |N|^2 &=& e^u \bigg[ - 2 |\nabla N|^2 + (\nabla^2
  u) * N^2  + Rm * N^2 \nonumber\\
  && + \nabla N * N * (N+\nabla u) + N^4 + N^3 * \nabla
  u + N^2 * (\nabla u)^2 \bigg] .
\eea
We remark that from this expression, we see that if $|N|^2=0$ at the
initial time, that $|N|^2 \equiv 0$ along the flow. To see this, we
assume $|N|^2 \leq 1$ and $|Rm|+|\nabla^2 u| + |\nabla u| \leq C$ on
$[0,\epsilon)$ and
apply the maximum principle to $e^{-At} |N|^2$ for $A \gg 1$.

\subsubsection{The evolution of the gradient of $N$}
We compute the evolution of $\nabla N$.
\be
\p_t \nabla_\ell N_{ij}{}^k = \nabla_\ell \dot{N}_{ij}{}^k -
N_{\lambda j}{}^k \dot{\Gamma}^\lambda_{\ell i} - N_{i \lambda}{}^k
\dot{\Gamma}^\lambda_{\ell j} + \dot{\Gamma}^k_{\ell \lambda} N_{ij}{}^\lambda
\ee
Using the equation (\ref{evolution-longtime}) for $\dot{g}_{ij}$, we can compute the time derivative of the Christoffel symbol as
\bea \label{evol-christoffel}
\dot{\Gamma}^k_{ij} &=&
- {g^{k \mu} \over 2}
(-\nabla_\mu \dot{g}_{ij} + \nabla_j \dot{g}_{\mu i} + \nabla_i \dot{g}_{\mu p}) \nonumber\\
&=& e^u ( \nabla Rm * g + \nabla^3 u * g + (\nabla N +
\nabla^2 u+1) * \mathcal{O}(\nabla u, N, Rm)).
\eea
Substituting (\ref{evol-Nijk}),
\bea
\partial_t \nabla_\ell N_{ij}{}^k &=& e^u \bigg[ \nabla_\ell \Delta
  N_{ij}{}^k + (\nabla Rm + \nabla^2 N + \nabla^3 u) * (N+\nabla u)\nonumber\\
&&+(\nabla N + \nabla^2 u) * (\nabla N + \nabla^2 u)+ (\nabla N
  +\nabla^2 u+1) * \mathcal{O}(\nabla u, N, Rm)\bigg] \nonumber
\eea
We can commute the derivatives $\nabla_\ell \Delta N_{ij}{}^k$ up to
lower order terms, and so
\bea \label{evol-DN0}
(\partial_t- e^u \Delta) \nabla_\ell N_{ij}{}^k &=& e^u \bigg[ (\nabla Rm + \nabla^2 N + \nabla^3 u) * (N+\nabla u)\nonumber\\
&&+(\nabla N + \nabla^2 u) * (\nabla N + \nabla^2 u)+ (\nabla N
  +\nabla^2 u+1) * \mathcal{O}(\nabla u, N, Rm)\bigg] \nonumber
\eea
The evolution of the norm $|\nabla N|^2$ is
\bea
(\p_t - e^u \Delta) |\nabla N|^2 &=& 2 \langle (\p_t - e^u \Delta)
\nabla N, \nabla N \rangle - 2 e^u | \nabla^2 N|^2 + \partial_t g
* \nabla N * \nabla N \nonumber
\eea
Altogether,
\bea \label{evol-DN}
(\p_t - e^u \Delta) |\nabla N|^2 &\leq& - 2 e^u | \nabla^2 N|^2 + C
e^u \bigg[ |\nabla N|^3 + |\nabla^2 u| |\nabla N|^2 + |\nabla^2 u|^2
  |\nabla N| \nonumber\\
&&+ |\nabla N| (|\nabla Rm| + |\nabla^2 N| + |\nabla^3 u|) (|N|+|\nabla u|)
  \nonumber\\
  &&+ \mathcal{O}(\nabla u, N, Rm) (|\nabla N|^2 + |\nabla N| |\nabla^2
  u| + 1)\bigg].
\eea

\subsection{The evolution of the curvature tensor: proof of Theorem \ref{th:Nijenhuis}(b)}
The general formula for the
variation of the curvature tensor (see e.g. \cite{Ha1}) is
\bea \label{evol-Rm-1}
   {d \over dt} R_{jikl} &=& {1 \over 2} (\nabla_i \nabla_k \dot{g}_{jl} - \nabla_i \nabla_l \dot{g}_{jk} - \nabla_j \nabla_k \dot{g}_{il} + \nabla_j \nabla_l \dot{g}_{ik}) \nonumber\\
   &&+ {1 \over 2} \dot{g}_{k \lambda} R_{ji}{}^\lambda{}_l - {1 \over 2} \dot{g}_{l \lambda} R_{ji}{}^\lambda{}_k.
\eea
In our case, we write the evolution of $g_{ij}$ (\ref{evolution-longtime}) as
\bea
\dot{g}_{ij} &=& e^u (-2R_{ij} + 2 \nabla_i \nabla_j u + E_{ij}), \nonumber\\
E_{ij} &=& - 4 (N^2_-)_{ij} + u_i u_j - u_{Ji} u_{Jj} + 4 u_p (N_i{}^p{}_j + N_j{}^p{}_i ) .
\eea
Differentiating once gives
\be
\nabla_k \dot{g}_{j \ell} = e^u (-2 \nabla_k R_{j \ell} + 2 \nabla_k \nabla_j
\nabla_\ell u + \nabla_k E_{j \ell}) + e^u (-2R_{j \ell} + 2 \nabla_j \nabla_\ell u
+ E_{j \ell}) \nabla_k u.
\ee
Differentiating twice gives
\bea
\nabla_i \nabla_k \dot{g}_{j \ell} &=& e^u (-2 \nabla_i \nabla_k R_{j
  \ell} + 2 \nabla_i \nabla_k \nabla_j \nabla_\ell u + \nabla_i
\nabla_k E_{j \ell}) \nonumber\\
&&+ e^u (-2 \nabla_k R_{j \ell} + 2 \nabla_k \nabla_j
\nabla_\ell u + \nabla_k E_{j \ell}) \nabla_i u + (i \leftrightarrow
k) \nonumber\\
&&+ e^u (-2R_{j \ell} + 2 \nabla_j \nabla_\ell u + E_{j \ell})
\nabla_i \nabla_k u + e^u (-2R_{j \ell} + 2 \nabla_j \nabla_\ell u +
E_{j \ell}) \nabla_k u \nabla_i u. \nonumber
\eea
We can group this as
\bea
\nabla_i \nabla_k \dot{g}_{j \ell} &=& e^u (-2 \nabla_i \nabla_k R_{j
  \ell} + 2 \nabla_i \nabla_k \nabla_j \nabla_\ell u) \nonumber\\
&&+ e^u (\nabla Rm + \nabla^3 u + \nabla^2 N) * \mathcal{O}(\nabla
u, N ) + \nabla N * (\nabla^2 u + \mathcal{O}(\nabla u,N))+ \nabla N * \nabla N \nonumber\\
&&+ \mathcal{O}(Rm,N,\nabla^2 u, \nabla u , u).
\eea
We have
\be
\nabla_i \nabla_k \nabla_j \nabla_\ell u - \nabla_i \nabla_j \nabla_k
\nabla_\ell u = \nabla_i (- R_{kj}{}^\lambda{}_\ell u_\lambda)
\ee
and (Lemma 7.2 in \cite{Ha1})
\be
\nabla_i \nabla_k R_{j \ell} - \nabla_i \nabla_\ell R_{jk} - \nabla_j
\nabla_k R_{i \ell} + \nabla_j \nabla_\ell R_{ik} = -\Delta R_{jik \ell} + Rm * Rm.
\ee
Substituting all this into (\ref{evol-Rm-1}), we obtain
\bea \label{evol-Rm-2}
 (\p_t - e^u \Delta) Rm &=& e^u \bigg[ (\nabla Rm + \nabla^3 u + \nabla^2 N) * \mathcal{O}(\nabla
u, N ) + \mathcal{O}(Rm,N, \nabla u) \nonumber\\
&& + \nabla N * (\nabla^2 u + \mathcal{O}(\nabla u,N))+ \nabla N
* \nabla N + \nabla^2 u * \nabla^2 u \bigg].
\eea
The norm evolves by
\be
(\p_t - e^u \Delta) |Rm|^2 = 2 \langle (\p_t - e^u \Delta) Rm, Rm
\rangle - 2 e^u |\nabla Rm|^2 + \p_t g * Rm * Rm.
\ee
Therefore
\bea
 (\p_t - e^u \Delta) |Rm|^2 &=& e^u \bigg[ - 2 |\nabla Rm|^2 + (\nabla
  Rm + \nabla^3 u + \nabla^2 N) *
  \mathcal{O}(Rm, \nabla u,N) \nonumber\\
  &&+ (\nabla N * \nabla N + \nabla^2 u * \nabla^2 u +1) * \mathcal{O}(Rm, \nabla u,N) \bigg].
\eea

\subsection{Lower order estimates}

\subsubsection{Gradient estimate}
In this section, we estimate the gradient of $u$.

\begin{proposition}
Suppose over a finite interval $[0,T)$ the flow exists and that
  $|u|+|Rm| \leq C$. then there exists a constant $C'$ such that
  $|\nabla u |^2\leq C'$ in the time interval $[0,T)$.
\end{proposition}

We recall that by our work so far, we know that $|N|^2$ is bounded and
the $J$-invariant part of $\nabla^2 u$ is bounded.

\smallskip
\par Equation (\ref{evol-nabla-u-2}) for the evolution of $|\nabla u|^2$
together with the evolution $(\p_t - e^u \Delta)u = e^u(2|\nabla u|^2 +
|N|^2)$ of $u$ imply
\bea
&&(\p_t-e^u\Delta)(e^{pu} |\nabla u |^2)\nonumber\\
&=&e^{(p+1)u}\bigg(-2 |\nabla^2
u|^2+(6-4p)(\nabla^2u)_{ij}u^iu^j+8(N^2_+)_{ij}u^iu^j+2u^s \nabla_s |N |^2\nonumber\\
&&+|\nabla u |^2(2\Delta u+(3+2p-p^2) |\nabla u |^2+(1+p) |N |^2)\bigg).\label{evolu2}
\eea
Let $V=e^{pu}(|N |^2+ |\nabla u|^2)$ for some constant $p$. From (\ref{evol-norm-N}) and (\ref{evolu2}), we see that
\bea
&&e^{-(p+1)u}(\p_t-e^u\Delta)V\nonumber\\
&=&\bigg(-2 |\nabla N|^2 + (\nabla^2 u) * N^2 + N * (Rm+\nabla
N) * (N+\nabla u)+ N * (N+\nabla u)^3 \nonumber\\
&&- 2p u^s \nabla_s |N|^2-(p^2-2p) | \nabla u |^2
|N |^2+ p |N |^4\bigg)\nonumber\\
&&+\bigg(-2 |\nabla^2
u |^2+(6-4p)(\nabla^2u)_{ij}u^iu^j+8(N^2_+)_{ij}u^iu^j+2u^s \nabla_s
|N |^2\nonumber\\
&&+ |\nabla u|^2(2\Delta u+(3+2p-p^2)|\nabla u|^2+(1+p)|N|^2)\bigg).\nonumber
\eea
We note that $(\nabla^2u) * N^2$ represent terms of the form $a (\nabla^2u)^{ij} (N^2_+)_{ij}+ b  (\nabla^2u)^{ij} (N^2_-)_{ij}$ for some constant $a$ and $b$. Those terms are also bounded since $N^2_\pm$ is $J$-invariant hence only the $J$-invariant part of $\nabla^2u$ contributes to this
term, which is bounded. Also, we can control all the terms linear in $\nabla N$ by the good term $-|\nabla N |^2$. Therefore
\bea
&&e^{-(p+1)u}(\p_t-e^u\Delta)V\nonumber\\
&\leq&-2 |\nabla^2
u |^2+(6-4p)(\nabla^2u)_{ij}u^iu^j+(3+2p-p^2) |\nabla u |^4
\nonumber\\
&&+C(p) |\nabla
u |^3+C(p).
\eea
To handle the term $(\nabla^2 u)_{ij} u^i u^j$, we need to make use of the fact that the $J$-invariant part of $\nabla^2u$ is bounded. To do so, let us denote the $J$-invariant and the $J$-anti-invariant parts of $\nabla^2u$ by $\nabla^2_Ju$ and $\nabla^2_{-J}u$ respectively. Under this notation, we see that
\bea
&&-2|\nabla^2 u|^2+(6-4p)(\nabla^2u)_{ij}u^iu^j+(3+2p-p^2)|\nabla u|^4\nonumber\\
&=&-2|\nabla^2_Ju|^2-2|\nabla^2_{-J}u|^2+(6-4p)(\nabla^2_Ju+\nabla^2_{-J}u,\nabla u\otimes\nabla u)+(3+2p-p^2)|\nabla u|^4\nonumber\\
&\leq&-2|\nabla^2_{-J}u|^2+(6-4p)(\nabla^2_{-J}u,\nabla u\otimes\nabla u)+(3+2p-p^2)|\nabla u|^4+C|\nabla u|^2.
\eea
The advantage of this consideration is that only the $J$-anti-invariant part of $\nabla u\otimes\nabla u$, namely $\frac{1}{2}(u_iu_j-u_{Ji}u_{Jj})$, contributes to the inner product term. Therefore
\bea
&&-2|\nabla^2 u|^2+(6-4p)(\nabla^2u)_{ij}u^iu^j+(3+2p-p^2)|\nabla u|^4\nonumber\\
&\leq&-2|\nabla^2_{-J}u|^2+(3-2p)(\nabla^2_{-J}u)^{ij}(u_iu_j-u_{Ji}u_{Jj})+(3+2p-p^2)|\nabla u|^4+C|\nabla u|^2\nonumber\\
&=&-2|(\nabla^2_{-J}u)_{ij}+\frac{1}{2}(p-\frac{3}{2})(u_iu_j-u_{Ji}u_{Jj})|^2+\frac{1}{2}(p-\frac{3}{2})^2|u_iu_j-u_{Ji}u_{Jj}|^2\nonumber\\
&&+(3+2p-p^2)|\nabla u|^4+C|\nabla u|^2\nonumber\\
&\leq&\frac{1}{2}(p-\frac{3}{2})^2|u_iu_j-u_{Ji}u_{Jj}|^2+(3+2p-p^2)|\nabla u|^4+C|\nabla u|^2.\nonumber
\eea
The $J$-anti-invariant part of $\nabla u\otimes\nabla u$ has half of the norm square compared to the full $\nabla u\otimes\nabla u$:
\bea
|u_iu_j-u_{Ji}u_{Jj}|^2=2|\nabla u\otimes\nabla u|^2=2|\nabla u|^4.
\eea
So the conclusion is that
\bea
&&-2|\nabla^2 u|^2+(2-4p)(\nabla^2u)_{ij}u^iu^j+(3+2p-p^2)|\nabla u|^4\nonumber\\
&\leq&((p-\frac{3}{2})^2-p^2+2p+3)|\nabla u|^4+C|\nabla u|^2\nonumber\\
&=&(-p+\frac{9}{4}+3)|\nabla u|^4+C|\nabla u|^2\nonumber\\
&\leq&-|\nabla u|^4+C|\nabla u|^2\nonumber
\eea
for $p=(9/4)+4$. Thus
\bea
e^{-(p+1)u}(\p_t-e^u\Delta)V\leq -|\nabla u|^4+C(p) |\nabla u|^3+C(p).
\eea
Then by maximum principle and the boundedness of $u$, we prove the proposition.

\subsubsection{Second order estimate}
In this section, we obtain estimates on $|\nabla N| + |\nabla^2 u|$. We refer to
$\nabla N$ and $\nabla^2 u$ as second order terms since they involve
two derivatives of $\varphi$.

\begin{proposition}
Let $(g_{ij}(t),u(t))$ evolve by Type IIA flow on $M \times
[0,T]$. Suppose
\be
\sup_{M \times [0,T]} \bigg( |Rm| + |N| +
|\nabla u| + |u| \bigg) \leq \Lambda.
\ee
Then there exists a constant $C$
depending on $\Lambda$ and $(g_{ij}(0),u(0))$ such that
\be
\sup_{M \times [0,T]} \bigg( |\nabla N| + |\nabla^2 u| \bigg) \leq C.
\ee
\end{proposition}

A basic building block in the
construction of our test function for this estimate will be
\be
\tau(z) = |N|^2 + |\nabla u|^2.
\ee
It satisfies $\tau \leq C$, and using our work so far, its evolution can be estimated by
\be \label{evol-tau}
(\p_t - e^u \Delta) \tau \leq - |\nabla N|^2 - |\nabla^2 u|^2 + C.
\ee
We start with the test function
\be
Q = {|\nabla^2 u|^2 +|\nabla N|^2 \over K - \tau}
\ee
where $K$ is a large constant to be determined.
We can compute its evolution
\bea
(\p_t - e^u \Delta) Q &=& {1 \over K- \tau} (\p_t - e^u \Delta)
(|\nabla^2 u|^2+ |\nabla N|^2)
+ {|\nabla^2 u|^2+|\nabla N|^2 \over (K- \tau)^2} (\p_t - e^u \Delta) \tau
\nonumber\\
&&- {2 e^u \over (K- \tau)^2} \nabla_i (|\nabla^2 u|^2+|\nabla N|^2) \nabla^i
  \tau -2 e^u
  {|\nabla^2 u|^2+ |\nabla N|^2 \over (K- \tau)^3} |\nabla \tau |^2. \nonumber
  \eea
  By our evolution equations (\ref{evol-DN}), (\ref{evol-D^2u}),
  (\ref{evol-tau}) for $|\nabla N|^2$, $|\nabla^2 u|^2$, and $\tau$,
  we obtain the estimate
  \bea
(\p_t - e^u \Delta) Q &\leq & {e^u \over (K-\tau)} \bigg[ - |\nabla^3 u|^2 - |\nabla^2
  N|^2 + C |\nabla N|^3 + C |\nabla^2 u|^3 \nonumber\\
&&+ C |\nabla Rm| |\nabla N| + C |\nabla Rm| |\nabla^2 u|+ C|\nabla N|
  |\nabla^3 u| \nonumber\\
  && + C|\nabla N|^2 |\nabla^2 u| + C|\nabla^2 N| |\nabla^2 u| +C \nonumber\\
&& - {|\nabla N|^4 \over (K-\tau)} - {|\nabla^2 u|^4 \over (K-\tau)} -
  2 {|\nabla^2 u|^2 |\nabla N|^2 \over (K-\tau)}
  \nonumber\\
  &&+ C {|\nabla^2 u|^2+|\nabla N|^2 \over (K- \tau)}
\nonumber\\
&&- {2 \over (K- \tau)} \nabla_i (|\nabla^2 u|^2+|\nabla N|^2) \nabla^i
  \tau -2
  {|\nabla^2 u|^2+|\nabla N|^2 \over (K- \tau)^2} |\nabla \tau |^2 \bigg] \nonumber
\eea
By the bound on $|N|$ and $|\nabla u|$, we can choose $K - \tau \geq {K \over 2}$ large. The terms $|\nabla
N|^4$ and $|\nabla^2 u|^4$ can absorb lower order terms. We also drop
the last term.
\bea
(\p_t - e^u \Delta) Q &\leq & {e^u \over (K-\tau)} \bigg[- |\nabla^3 u|^2 - |\nabla^2
  N|^2 + {1 \over 10} |\nabla Rm|^2 +C(K) \nonumber\\
&& - {|\nabla N|^4 \over 2K} - {|\nabla^2 u|^4 \over 2K} - {|\nabla^2 u|^2 |\nabla N|^2 \over K}
\nonumber\\
&&- {2 \over (K- \tau)} \nabla_i (|\nabla^2 u|^2+|\nabla N|^2) \nabla^i\tau\bigg]
\eea
Using $|\nabla \tau| \leq C (|\nabla^2 u| + |\nabla N| + 1)$, we can estimate
\bea
& \ & - {2 \over (K- \tau)} \nabla_i (|\nabla^2 u|^2+|\nabla N|^2) \nabla^i
\tau \nonumber\\
     &\leq& {C \over K} |\nabla |\nabla^2 u|^2 + \nabla |\nabla N|^2| |\nabla
\tau| \nonumber\\
&\leq& {C \over K} \bigg( |\nabla^3 u| (|\nabla^2 u|^2 + |\nabla^2 u| |\nabla
N|) + |\nabla^2 N| (|\nabla N|^2 + |\nabla^2 u||\nabla N|) \bigg) \nonumber\\
&\leq& {1 \over 2} |\nabla^3 u|^2 + {1 \over 2} |\nabla^2 N|^2 + {C_0 \over K^2} |\nabla^2 u|^4 + {C_0
  \over K^2} |\nabla^2 u|^2 |\nabla N|^2 + {C_0 \over K^2} |\nabla N|^4
\eea
Choose $K$ large such that $K \geq 4C_0 \gg 1$. Then the main
inequality becomes
\bea \label{evol-Q}
(\p_t - e^u \Delta) Q &\leq & {e^u \over (K-\tau)} \bigg[- {1 \over 2}
  |\nabla^3 u|^2 - {1 \over 2} |\nabla^2
  N|^2 - {|\nabla N|^4 \over 4K} - {|\nabla^2 u|^4 \over 4K}
  \nonumber\\
  &&+ {1 \over 10} |\nabla Rm|^2 +C(K) \bigg].
\eea
We can now prove that if $|u|+|\nabla u|+|N| + |Rm| \leq C$ along the flow, then we
can bound $|\nabla N|$ and $|\nabla^2 u|$. Consider the test function
\be
S = Q + |Rm|^2.
\ee
By (\ref{evol-Q}) and (\ref{evol-Rm-3}), we can estimate the evolution
of $Q$ and $|Rm|^2$.
\bea
(\p_t - e^u \Delta) S &\leq & {e^u \over (K-\tau)} \bigg[- {1 \over 2}
  |\nabla^3 u|^2 - {1 \over 2} |\nabla^2
  N|^2 - {|\nabla N|^4 \over 4K} - {|\nabla^2 u|^4 \over 4K}
  \nonumber\\
  &&+ {1 \over 10} |\nabla Rm|^2 +C(K) \bigg] - e^u |\nabla Rm|^2 \nonumber\\
  &&+ C |\nabla^3 u| + C
|\nabla^2 N| + C |\nabla N|^2 + C |\nabla^2 u|^2 + C
\eea
As long as $K$ is large enough such that $K - \tau \geq 1$, it follows
that at a maximum point $(x_0,t_0)$ of $S$ with $t_0>0$, then
\be
|\nabla N|^2(x_0,t_0) + |\nabla^2 u|^2(x_0,t_0) \leq C(K).
\ee
Since $Q = {|\nabla^2 u|^2 + |\nabla N|^2 \over K - \tau}$, it follows
that $Q(x_0,t_0) \leq C$ and hence $S(x_0,t_0) \leq C$.
\smallskip
\par Therefore $S$ is bounded on $M \times [0,T]$. It follows that if $|u|
+ |\nabla u| + |N|+ |Rm| \leq C_0$ on $M \times [0,T]$, then
\be
|\nabla N| + |\nabla^2 u| \leq C
\ee
where $C$ depends on $C_0$ and the initial data.

\subsection{Higher order estimates}
In this section, we prove the following estimate.
\begin{proposition}
Let $(g_{ij}(t),u(t))$ evolve by Type IIA flow on $M \times
[0,T]$. Suppose
\be
\sup_{M \times [0,T]} \bigg( |Rm| + |\nabla N| + |N| + |\nabla^2 u| +
|\nabla u| + |u| \bigg) \leq \Lambda.
\ee
Then for each integer $k \geq 1$, there exists a constant $C_k$
depending on $k$, $\Lambda$ and $(g_{ij}(0),u(0))$ such that
\be
\sup_{M \times [0,T]} \bigg( |\nabla^k Rm| + |\nabla^{k+1} N| + |\nabla^{k+2} u| \bigg) \leq C_k.
\ee

\end{proposition}

Note: in earlier work, we showed that the estimate $|u|+|Rm| \leq C$
implies the estimate $|\nabla u| + |\nabla^2 u| + |N| + |\nabla N|
\leq C$. Combining these two results, we conclude that if  $|u|+|Rm|
\leq C$ remains bounded along the flow, then all geometric terms
remain bounded.

\medskip
\par Let $I_k$ denote any combination of geometric terms of derivative order $\leq k$
in the metric. For example,
\bea
I_2 &=& f(u,\nabla u, \nabla^2 u, N, \nabla N, Rm), \nonumber\\
I_3 &=& f(u,\nabla u, \nabla^2 u, \nabla^3 u, N, \nabla N, \nabla^2 N, Rm,
\nabla Rm).
\eea
In this section, we will evolve all higher order geometric terms
appearing in the equation of the metric Type IIA flow.

\subsubsection{The evolution of $|\na^kRm|^2$}
We write the evolution of the curvature as
\be
(\p_t - e^u \Delta) Rm = E(Rm),
\ee
where
\be
E(Rm) = (\nabla^3 u + \nabla^2 N + \nabla Rm) * I_1 + I_2.
\ee
We have for example
\be
\nabla E(Rm) = (\nabla^4 u + \nabla^3 N + \nabla^2 Rm) * I_1 +
(\nabla^3 u + \nabla^2 N + \nabla Rm) * I_2 + I_2
\ee
and in general
\bea
\nabla^k E(Rm) &=& (\nabla^{k+3} u + \nabla^{k+2} N + \nabla^{k+1} Rm) * I_1
\nonumber\\
&&+ (\nabla^{k+2} u + \nabla^{k+3} N + \nabla^{k+1} Rm) * I_2 + I_{k+1}
\eea
Then
\bea
\p_t \nabla^k Rm &=& \p_t (\p + \Gamma)^k Rm \nonumber\\
&=& \nabla^k (\p_t Rm) + \sum_{i=0}^{k-1} \nabla^i \p_t \Gamma \,
\nabla^{k-1-i} Rm \nonumber\\
&=& \nabla^k E(Rm) + \nabla^k(e^u \Delta Rm) +  \sum_{i=0}^{k-1}
\nabla^i \p_t \Gamma \,
\nabla^{k-1-i} Rm
\eea
We have the general commutator formula
\be \label{commute-laplace-nablak}
\nabla^k \Delta A= \Delta \nabla^k A+ \nabla^k (Rm * A)
\ee
which implies
\bea
(\p_t - e^u \Delta) \nabla^k Rm &=& \nabla^k (Rm * Rm) + \nabla^k
E(Rm) + \sum_{i=1}^k \nabla^i e^u * \nabla^{k-i}
\Delta Rm \nonumber\\
&&+  \sum_{i=0}^{k-1}
\nabla^i \p_t \Gamma *
\nabla^{k-1-i} Rm
\eea
We note
\be \label{evol-Gamma}
\p_t \Gamma = (\nabla^3 u + \nabla^2 N + \nabla Rm) * I_1 + I_2,
\ee
and
\bea
\nabla^k \p_t \Gamma &=& (\nabla^{k+3} u + \nabla^{k+2} N + \nabla^{k+1} Rm) * I_1
\nonumber\\
&&+ (\nabla^{k+2} u + \nabla^{k+3} N + \nabla^{k+1} Rm) * I_2 + I_{k+1}.
\eea
Therefore
\bea
(\p_t - e^u \Delta) \nabla^k Rm
&=& I_1 * (\nabla^{k+3} u + \nabla^{k+2} N + \nabla^{k+1} Rm)
\nonumber\\
&&+ I_2 * (\nabla^{k+2} u + \nabla^{k+1} N + \nabla^k Rm) + I_{k+1}
\eea
The norm is evolving by
\bea
(\p_t - e^u \Delta) |\nabla^k Rm|^2 &=& 2 \langle (\p_t - e^u \Delta)
\nabla^k Rm, \nabla^k Rm \rangle \nonumber\\
&&- 2 e^u |\nabla^{k+1} Rm|^2 + \p_t g
* \nabla^k Rm * \nabla^k Rm.
\eea
Thus
\bea
(\p_t - e^u \Delta) |\nabla^k Rm|^2 &=& - 2 e^u |\nabla^{k+1} Rm|^2
\nonumber\\
&&+ I_1 * (\nabla^{k+3} u + \nabla^{k+2} N + \nabla^{k+1} Rm)
* \nabla^k Rm \nonumber\\
&&+ I_2 * (\nabla^{k+2} u + \nabla^{k+1} N + \nabla^k Rm)^2 +
I_{k+1} * \nabla^k Rm.
\eea

\subsubsection{The evolution of $|\na^kN|^2$}
We will evolve $\nabla^k N$ in this section for all $k \geq 2$. We write $(\p_t - e^u \Delta) N = E(N)$. Higher order terms evolve by
\bea
\p_t \nabla^k N &=& \p_t (\p + \Gamma)^k N\nonumber\\
&=& \nabla^k (\p_t N) + \sum_{i=0}^{k-1} \nabla^i \p_t \Gamma *
\nabla^{k-1-i} N \nonumber\\
&=& \nabla^k E(N) + e^u \nabla^k \Delta N + \sum_{i=1}^k
\nabla^i e^u * \nabla^{k-i} \Delta N  \nonumber\\
&&+  \sum_{i=0}^{k-1} \nabla^i \p_t \Gamma *
\nabla^{k-1-i} N.
\eea
Using (\ref{commute-laplace-nablak}) to commute derivatives gives
\bea \label{evol-nablak-N}
(\p_t - e^u \Delta) \nabla^k N &=& \nabla^k E(N) + \sum_{i=1}^k
\nabla^i e^u * \nabla^{k-i} \Delta N \nonumber\\
&&+ \nabla^k (Rm * N) +  \sum_{i=0}^{k-1} \nabla^i \p_t \Gamma *
\nabla^{k-1-i} N.
\eea
By (\ref{evol-Nijk}),
\be
E(N) = (\nabla^2 u + \nabla N + Rm) * I_1 + I_1.
\ee
Differentiating this once gives
\be
\nabla E(N) = (\nabla^3 u + \nabla^2 N + \nabla Rm) * I_1 + I_2.
\ee
Differentiating again, we obtain
\bea
\nabla^2 E(N) &=& (\nabla^4 u + \nabla^3 N + \nabla^2 Rm) * I_1
\nonumber\\
&&+ (\nabla^3 u + \nabla^2 N + \nabla Rm) * I_2 + I_2.
\eea
Higher order derivatives are
\bea
\nabla^k E(N) &=& (\nabla^{k+2} u + \nabla^{k+1} N + \nabla^k Rm)
* I_1 \nonumber\\
&&+ (\nabla^{k+1} u + \nabla^{k} N + \nabla^{k-1} Rm) * I_2 + I_{k}
\eea
for $k \geq 2$. Substituting this and (\ref{evol-Gamma}) into (\ref{evol-nablak-N})
\bea
(\p_t - e^u \Delta) \nabla^k N &=& (\nabla^{k+2} u + \nabla^{k+1} N + \nabla^k Rm)
* I_1 \nonumber\\
&&+ (\nabla^{k+1} u + \nabla^{k} N + \nabla^{k-1} Rm) * I_2 + I_{k}.
\eea
The norm evolves by
\bea
(\p_t - e^u \Delta) |\nabla^k N|^2 &=& - 2 e^u |\nabla^{k+1} N|^2 + 2 \langle  (\p_t - e^u \Delta)
\nabla^k N, \nabla^k N \rangle \nonumber\\
&&+ \p_t g * \nabla^k N * \nabla^k N.
\eea
Therefore
\bea
(\p_t - e^u \Delta) |\nabla^k N |^2 &=& - 2 e^u |\nabla^{k+1} N|^2 +
I_1 * (\nabla^{k+2} u + \nabla^{k+1} N + \nabla^k Rm) *
\nabla^k N \nonumber\\
&&+ I_2 * (\nabla^{k+1} u + \nabla^{k} N + \nabla^{k-1} Rm)^2 +
I_{k} * \nabla^k N.
\eea

\subsubsection{The evolution of $|\na^ku|^2$}
Denote as before $(\p_t - e^u \Delta) u = E(u)$. We will compute the
evolution of $\nabla^k u$ for $k \geq 3$.
\bea \label{evol-nablak-u}
\p_t \nabla^k u &=& \p_t (\partial+\Gamma)^{k-1} \partial u \nonumber\\
&=& \nabla^k (\p_t u) + \sum_{i=0}^{k-2} \nabla^i \p_t
\Gamma * \nabla^{k-1-i} u \nonumber\\
&=& e^u \nabla^k \Delta u+ \nabla^k E(u) + \sum_{i=1}^k \nabla^i e^u
\nabla^{k-i} \Delta u + \sum_{i=0}^{k-2} \nabla^i \p_t
\Gamma * \nabla^{k-1-i} u \nonumber\\
&=& e^u \Delta \nabla^k u + \sum_{i=0}^{k-1} \nabla^i Rm * \nabla^{k-i}
u + \nabla^k E(u) \nonumber\\
&&+ \sum_{i=1}^k \nabla^i e^u
\nabla^{k-i} \Delta u + \sum_{i=0}^{k-2} \nabla^i \p_t
\Gamma * \nabla^{k-1-i} u
\eea
The evolution of $u$ is of the form $E(u) = I_1$. We will
differentiate this 3 times before it becomes linear enough to use in
our general argument. Differentiating once
\be
\nabla E(u) = (\nabla^2 u + \nabla N + Rm) * I_1 + I_1,
\ee
twice
\be
\nabla^2 E(u) = (\nabla^3 u + \nabla^2 N + \nabla Rm) * I_1 + I_2
\ee
and three times
\be
\nabla^3 E(u) = (\nabla^4 u + \nabla^3 N + \nabla^2 Rm) * I_1 +
(\nabla^3 u + \nabla^2 N + \nabla Rm) * I_2 + I_2.
\ee
Higher order derivatives are
\be
\nabla^k E(u) = (\nabla^{k+1} u + \nabla^k N + \nabla^{k-1} Rm) * I_1 +
(\nabla^k u + \nabla^{k-1} N + \nabla^{k-2} Rm) * I_2 + I_{k-1}.
\ee
for $k \geq 3$. Substituting this and (\ref{evol-Gamma}) into (\ref{evol-nablak-u})
\bea
(\p_t - e^u \Delta) \nabla^k u &=& (\nabla^{k+1} u + \nabla^{k} N + \nabla^{k-1} Rm)
* I_1 \nonumber\\
&&+ (\nabla^{k} u + \nabla^{k-1} N + \nabla^{k-2} Rm) * I_2 + I_{k-1}.
\eea
Using the evolution of the norm
\bea
(\p_t - e^u \Delta) |\nabla^k u|^2 &=& - 2 e^u |\nabla^{k+1} u|^2 + 2
\langle (\p_t - e^u \Delta) \nabla^k u, \nabla^k u \rangle \nonumber\\
&&+ \partial_t g * \nabla^k u * \nabla^k u,
\eea
we conclude
\bea
(\p_t - e^u \Delta) |\nabla^k u|^2 &=& - 2 e^u |\nabla^{k+1} u|^2  + I_1 * (\nabla^{k+1} u + \nabla^{k} N + \nabla^{k-1} Rm)
* \nabla^{k} u \nonumber\\
&&+ I_2 * (\nabla^{k} u + \nabla^{k-1} N + \nabla^{k-2} Rm)^2 +
I_{k-1} * \nabla^k u.
\eea

\subsubsection{Estimates: proof of Theorem \ref{th:Shi}}
Putting everything together, we obtain
\bea
& \ & (\p_t - e^u \Delta) (|\nabla^k u|^2 + |\nabla^{k-1} N|^2 +
|\nabla^{k-2} Rm|^2) \nonumber\\
&=& - 2 e^u (|\nabla^{k+1} u|^2 + |\nabla^k N|^2 + |\nabla^{k-1}
Rm|^2) \nonumber\\
&&+ I_1 * (\nabla^{k+1} u + \nabla^{k} N + \nabla^{k-1} Rm)
* (\nabla^{k} u + \nabla^{k-1} N + \nabla^{k-2} Rm)\nonumber\\
&&+ I_2 * (\nabla^{k} u + \nabla^{k-1} N + \nabla^{k-2} Rm)^2 +
I_{k-1} * (\nabla^k u + \nabla^{k-1} N + \nabla^{k-2} Rm)
\eea
Let $k \geq 3$. Suppose $I_{k-1} \leq C$. Then
\bea
& \ & (\p_t - e^u \Delta) (|\nabla^k u|^2 + |\nabla^{k-1} N|^2 +
|\nabla^{k-2} Rm|^2) \nonumber\\
&\leq& - e^u (|\nabla^{k+1} u|^2 + |\nabla^k N|^2 + |\nabla^{k-1}
Rm|^2) \nonumber\\
&&+ C |\nabla^k u|^2 +C |\nabla^{k-1} N|^2 + C |\nabla^{k-2} Rm|^2 + C
\eea
and
\bea
& \ & (\p_t - e^u \Delta) (|\nabla^{k-1} u|^2 + |\nabla^{k-2} N|^2 +
|\nabla^{k-3} Rm|^2) \nonumber\\
&\leq& - e^u (|\nabla^{k} u|^2 + |\nabla^{k-1} N|^2 + |\nabla^{k-2}
Rm|^2) + C.
\eea
It follows that the test function
\bea
& \ & (\p_t - e^u \Delta) \bigg[ |\nabla^{k} u|^2 + |\nabla^{k-1} N|^2 +
|\nabla^{k-2} Rm|^2 + \Lambda (|\nabla^{k-1} u|^2 + |\nabla^{k-2} N|^2 +
|\nabla^{k-3} Rm|^2) \bigg] \nonumber\\
&\leq& -|\nabla^{k} u|^2 - |\nabla^{k-1} N|^2 - |\nabla^{k-2}
Rm|^2 + \Lambda C
\eea
for $ \Lambda \gg 1$ large. By the maximum principle, we conclude that
if $I_{k-1} \leq K$ then
\be
|\nabla^{k} u|^2 + |\nabla^{k-1} N|^2 + |\nabla^{k-2}
Rm|^2 \leq C(K,g(0))
\ee
and hence $I_k$ is bounded along the flow. This argument shows that if
$I_2$ is bounded, then $I_k$ is bounded for all $k$.

\subsection{Long-time existence}
Let $(u(t),g(t))$ be a solution to the Type IIA flow on
$[0,T)$. Suppose $|u|+|Rm| \leq C$ on $M \times [0,T)$. We have shown
that in this case $|\nabla^k u| + |\nabla^k N| + |\nabla^k Rm| \leq C$
for all $k \geq 1$. We now give the standard argument (see
e.g. \cite{Ha1}) which shows that the flow can be extended past
$t=T$. Denote $\p_t g_{ij} = E_{ij}$. Our estimates imply that
\be
|E| + |\nabla^k E| \leq C
\ee
where $\nabla$ is with respect to the evolving metric $g(t)$. If we
take $x \in M$, $v \in T_x M$ and $t_1,t_2 \in (0,T)$, then
\be
\bigg| \log g(t_2)(v,v) - \log g(t_1)(v,v) \bigg| = \bigg|
\int_{t_1}^{t_2} {\dot{g}(\tau)(v,v) \over g(\tau)(v,v)} d \tau \bigg|
\leq C |t_2-t_1|.
\ee
It follows that $g(t)$ is a Cauchy sequence as $t \rightarrow T$ and
\be
e^{-CT} g(0) \leq g(t) \leq e^{CT} g(0)
\ee
and the metrics $g_{ij}$ do not degenerate on $[0,T)$. Let
$\bar{\nabla}$ denote the covariant derivative with respect to $\bar g=g(0)$. We have
\be
\p_t \bar{\nabla}_k g_{ij} = \bar{\nabla}_k \dot{g}_{ij} = \nabla_k
E_{ij} + (\bar{\Gamma} - \Gamma) * E_{ij}
\ee
The difference between two connections is
\be
\Gamma^k_{ij} - \bar{\Gamma}^k_{ij} = {1 \over 2}
g^{kp}(-\bar{\nabla}_p g_{ij} + \bar{\nabla}_i g_{pj} + \bar{\nabla}_j g_{pi}),
\ee
and hence
\be
\p_t \bar{\nabla} g = \bar{\nabla} g * E + \mathcal{O}(1).
\ee
Therefore
\be
\p_t |\bar{\nabla} g|^2_{\bar{g}} \leq C |\bar{\nabla} g|^2_{\bar{g}} + C |\bar{\nabla}g|_{\hat{g}}
\ee
and hence $|\bar{\nabla} g|_{\bar{g}} \leq C(T)$ on $[0,T)$. Higher order derivatives are similar: indeed, let $k \geq 1$ and suppose that $|\bar{\nabla}^\ell g|_{\bar{g}}\leq C$ for all $\ell \leq k$. Then a similar calculation gives
\be
\p_t \bar{\nabla}^{k+1} g = \bar{\nabla}^{k+1} g * E + \mathcal{O}(1),
\ee
from which it follows that
\be
|\bar{\nabla}^{k+1} g|_{\bar{g}} \leq C(T).
\ee
Therefore, the evolving metrics $g_{ij}$ and all their derivatives are
bounded on $[0,T)$. Since we showed $g(t)$ is Cauchy, it follows that
$g(t) \rightarrow g(T)$ smoothly as $t \rightarrow T$. A similar
  argument shows that $u(t) \rightarrow u(T)$.
\smallskip
\par This produces a limiting pair $(g(T),u(T))$. The linear ODE for
$\varphi$ given in Theorem \ref{th:g-varphi} has coefficients which
only depend on $\tilde{g}_{ij} = e^u g_{ij}$, thus these coefficients
are smoothly defined on $[0,T]$. It follows that $\varphi(t)$ has a
smooth solution on $[0,T]$. By the non-degeneracy estimate
(\ref{minimum}), we have that $\varphi(T)$ is closed, primitive, and
in the positive cone $(-\lambda_\varphi)>0$. By the short-time existence theorem, the
flow can be extended to $[0,T+\epsilon)$ for some $\epsilon>0$. The
  discussion here also implies that $|\nabla^\alpha \varphi| \leq C$ on
  $[0,T)$ for any multi-index $\alpha$. This completes the proof of Theorem \ref{th:Shi}.

\section{Examples and Applications}
\setcounter{equation}{0}
\label{s: examples}

In this section, we discuss a range of examples and applications of the Type IIA flow with no sources.

\subsection{The stationary points: proof of Theorem \ref{th:stationary}}

We begin with the proof of Theorem \ref{th:stationary}. First, we note that
$-d^\Lambda(|\varphi|^2\hat\varphi)=\Lambda d(|\varphi|^2\hat\varphi)=\p_-(|\varphi|^2\hat\varphi)$, where $\p_-$ is a first order differential operator introduced in \cite{TY1} such that
\be
d(|\varphi|^2\hat\varphi)=\omega\wedge\p_-(|\varphi|^2\hat\varphi)\nonumber
\ee
and $\p_-(|\varphi|^2\hat\varphi)$ is a primitive 2-form. Therefore
the stationary point equation $d \Lambda d (|\varphi|^2 \star \varphi)
= 0$ can be expressed as $d\p_-(|\varphi|^2\hat\varphi)=0$. In particular,
\be
0=\int_Md\p_-(|\varphi|^2\hat\varphi)\wedge\hat\varphi=-\int_M\p_-(|\varphi|^2\hat\varphi)\wedge d\hat\varphi.\label{ibp}
\ee
Combining (\ref{dstarvarphi}), (\ref{dconnection}), and $\beta=0$
implies $d \hat{\varphi} = \alpha \wedge \hat{\varphi} + \tilde{{\frak T}}
\boxtimes \hat{\varphi}$. By (\ref{torsiontilde}), we have
$\tilde{{\frak T}} = N - {1 \over 4} (d^c \tilde{\omega} + {\cal M}(d^c \tilde{\omega}))$. By Lemma \ref{uvarphi},
we see that
\be
d\hat\varphi=N\boxtimes\hat\varphi
\ee
which is a (2,2)-form (by the argument in the proof of Lemma \ref{22}). Therefore
\be
d(|\varphi|^2\hat\varphi)=|\varphi|^2(-\alpha\wedge\hat\varphi+N\boxtimes\hat\varphi),\label{rhs}
\ee
where the first term is a $(3,1)+(1,3)$-form, the second is a $(2,2)$-form, and $\alpha=-d|\varphi|^2$. It follows that
\bea
\p_-\hat\varphi&=&\Lambda(N\boxtimes\hat\varphi),\nonumber\\
\p_-(|\varphi|^2\hat\varphi)&=&|\varphi|^2(-\Lambda(\alpha\wedge\hat\varphi)+\p_-\hat\varphi),\nonumber
\eea
where $\Lambda(N\boxtimes\hat\varphi)$ is a $(1,1)$-form and $\Lambda(\alpha\wedge\hat\varphi)$ is a $(2,0)+(0,2)$-form. Thus (\ref{ibp}) becomes
\bea
0&=&\int_M|\varphi|^2(-\Lambda(\alpha\wedge\hat\varphi)+\p_-\hat\varphi)\wedge\omega\wedge\p_-\hat\varphi =\int_M|\varphi|^2\omega\wedge\p_-\hat\varphi\wedge\p_-\hat\varphi\nonumber\\
&=&-\int_M|\varphi|^2|\p_-\hat\varphi|^2\frac{\omega^3}{3!}.
\eea
Consequently we conclude that $\p_-\hat\varphi=0$ and $d\hat\varphi=\omega\wedge\p_-\hat\varphi=0$. In view of Lemma \ref{lm:Omega}, the almost-complex structure $J$ is integrable and the form $\varphi$ is harmonic. Now we use the integration by parts argument again to get
\bea
0&=&\int_Md\p_-(|\varphi|^2\hat\varphi)\wedge|\varphi|^2\hat\varphi = -\int_M\omega\wedge\p_-(|\varphi|^2\hat\varphi)\wedge\p_-(|\varphi|^2\hat\varphi)\nonumber\\
&=&-\int_M|\varphi|^4|\Lambda(\alpha\wedge\hat\varphi)|^2\frac{\omega^3}{3!}.\nonumber
\eea
So we deduce that $\Lambda(\alpha\wedge\hat\varphi)=0$, hence $\alpha\wedge\hat\varphi=0$ and $\alpha=0$, so $|\varphi|$ is a constant. Q.E.D.

\subsection{Integrable almost-complex structures: proof of Theorem \ref{th:integrable}}

Next, we give the proof of Theorem \ref{th:integrable}.
The following identity for any smooth function $f$ and any differential form is well-known:
\bea
\label{Lie}
d^\dagger(f\mu)=f d^\dagger\mu-\iota_{\nabla f}\mu.
\eea
Indeed, it can be quickly verified by using $d^\dagger\mu=-\iota_k(\nabla^k\mu)$.

\if
Indeed, the adjoint $d^\dagger$ on an $n$-dimensional manifold is given by
\be
d^\dagger \alpha = (-1)^{n(r+1)+1} \star d \star \alpha, \ \ \ \alpha \in \Lambda^r(M).
\ee
Therefore
\be
d^\dagger (f \mu) = f d^\dagger \mu + (-1)^{n(r+1) +1} \star (d f \wedge \star \mu).
\ee
Using the fact that $\star^2 \alpha = (-1)^{r(n-r)} \alpha$ and the identity for the Hodge star operator
\be
\star (\eta \wedge \alpha) = (-1)^r  \iota_{\eta^{\sharp}} \star \alpha
\ee
for all $ \alpha \in \Lambda^r(X)$, $\eta \in \Lambda^1(X)$, and $(\eta^\sharp)^i = g^{ik} \eta_k$, we obtain the desired identity.
\fi

\bigskip
Back to the proof of Theorem \ref{th:integrable}, we apply Theorem \ref{th:Laplacian} and the identity (\ref{Lie}) to rewrite the Type IIA flow without sources as
\bea
\p_t\varphi&=&-d(|\varphi|^2d^\dagger\varphi-\iota_{\nabla|\varphi|^2}\varphi)+2d(|\varphi|^2 N^\dagger\cdot\varphi)
\nonumber\\
&=&{\cal L}_{\nabla|\varphi|^2}\varphi-d(|\varphi|^2d^\dagger\varphi)+2d(|\varphi|^2 N^\dagger\cdot\varphi).
\eea
On any orbit of the diffeomorphism group which contains a form $\varphi$ with an integrable almost-complex structure, we have $N=0$ and, in view of Lemma \ref{lm:Omega}, $d^\dagger\varphi=0$. Thus the flow reduces to
\bea
\p_t\varphi={\cal L}_{\na|\varphi|^2}\varphi
\eea
and a solution is given by the reparametrizations of $\varphi$ along the time-dependent vector field $\na|\varphi|^2$. By the uniqueness part of Theorem \ref{th:short}, this is the unique solution.

\medskip

We now re-express the Type IIA flow in an equivalent formulation, but with a fixed complex structure. For this, let $f_t$ be the flow generated by the time-dependent vector field $-\nabla|\varphi|^2$ in the sense that
\be
\frac{d}{dt}f_t(x)=-\nabla|\varphi|^2(t,x)
\ee
for any $x\in X$ and time $t$. It follows that
\be
\frac{d}{dt}(f_t^*\varphi_t)=f_t^*\p_t\varphi_t+f_t^*{\cal L}_{-\nabla|\varphi|^2}\varphi_t=0,\nonumber
\ee
hence $f_t^*\varphi_t\equiv\varphi_0$ is a constant 3-form on $X$. We see immediately that if we reparametrize $M$ by the time-dependent diffeomorphism $f_t$, then along the flow, the 3-form $f_t^*\varphi_t$ and hence the complex structure $f_t^*J_t$ are fixed. In this new gauge, the K\"ahler metric $\omega_t=f_t^*\o$ evolves by the equation
\be
\p_t\o_t=f_t^*{\cal L}_{-\nabla|\varphi|^2}\omega.
\ee
Notice that
\bea
{\cal L}_{-\nabla|\varphi|^2}\omega=-d\iota_{\nabla|\varphi|^2}\omega=d J d|\varphi|^2=-dd^c|\varphi|^2
\eea
so the flow of $\omega_t$ can be written as
\be
\p_t\omega_t=-dd^c|\varphi|^2_{\o_t}.\label{eq: gauged IIA flow}
\ee
We remark that this equation can be viewed as a T-dual of the Anomaly flow for conformally K\"ahler data. In that case, we are given a fixed holomorphic $(3,0)$ form $\check{\Omega}$ on a Calabi-Yau threefold and the evolving K\"ahler metrics $\check{\omega}_t$ satisfy (see equation (4.10) in \cite{FeiPicard})
\be
\p_t \check{\omega}_t = d d^c | \check{\Omega} |_{\check{\omega}_t}^{-2}.
\ee
We see that $\check{R} = | \check{\Omega} |_{\check{\omega}}^{-2}$ is exchanged with $1/R = | \varphi |^2_{\o}$.
\smallskip
\par This type of duality was observed in \cite{FeiPicard} on semi-flat Calabi-Yau threefolds. To connect with the work there, we can consider the conformally changed metric $\eta_t=|\varphi|^{-2}_{\o_t}\o_t$. It follows that $|\varphi|_{\eta_t}=|\varphi|^4_{\o_t}$, hence $\eta_t$ satisfies the conformally balanced equation $d(|\varphi|_{\eta_t}\eta_t^2)=0$ since
\be
|\varphi|_{\eta_t}\eta_t^2=|\varphi|^4_{\o_t}(|\varphi|^{-2}_{\o_t}\o_t)^2=\o_t^2\nonumber
\ee
is closed. Moreover its evolution equation is
\bea
\p_t(|\varphi|_{\eta_t}\eta_t^2)&=&\p_t(\omega_t^2)=2\o_t\wedge\p_t\o_t\nonumber\\
&=&2\o_t\wedge(-dd^c|\varphi|^2_{\o_t})=-4i\p\bar\p(|\varphi|^2_{\o_t}\o_t)\nonumber\\
&=&-4i\p\bar\p(|\varphi|_{\eta_t}\eta_t),
\eea
which is exactly (up to a positive constant) the dual Anomaly flow in complex dimension 3 firstly introduced in \cite{FeiPicard}. Since we are in the conformally K\"ahler case, by the results of \cite{FeiPicard}, we know the flow (\ref{eq: gauged IIA flow}) is equivalent to the inverse MA-flow introduced by Cao-Keller \cite{CK} and Collins-Hisamoto-Takahashi \cite{CHT}, which converges to the unique Ricci-flat K\"ahler metric in the cohomology class $[\o_0]$.

\subsection{Symplectic manifolds with non-integrable almost complex structures}

Next, we work out the Type IIA flow on some model symplectic manifolds with non-integrable almost-complex structures, more specifically tori, symplectic half-flat manifolds, and nilmanifolds.

\subsubsection{The Type IIA flow on a torus}
\label{ex:torus}

Consider the 6-torus $M=(\R/\Z)^6$, with coordinates $\{x^j\}_{j=1}^6$ and the standard symplectic form $\o=d x^{12}+dx^{34}+dx^{56}$. Let $\alpha,\beta,\gamma,\delta:\R/\Z\to\R$ be smooth functions depending only on the variable $x^1$. Consider
\bea
\varphi&=&e^\alpha dx^{135}-e^\beta dx^{146}-dx^{245}-dx^{236}+\gamma dx^{136}+\delta dx^{145}.\label{torusans}
\eea
Clearly $\varphi$ is closed and primitive. It is straightforward to find out that
\be
|\varphi|^2=2\sqrt{4e^{\alpha+\beta}-(\gamma-\delta)^2}.\nonumber
\ee
and
\bea
\hat{\varphi}&=&\frac{2}{|\varphi|^2}(-e^\alpha(\gamma+\delta)dx^{135}+(2e^{\alpha+\beta}-\gamma(\gamma-\delta))dx^{136}+(2e^{\alpha+\beta} +\delta(\gamma-\delta))dx^{145}\nonumber\\
&&+e^\beta(\gamma+\delta)dx^{146}+2e^\alpha dx^{235}+(\gamma-\delta)dx^{236}-(\gamma-\delta)dx^{245}-2e^\beta dx^{246}).\nonumber
\eea
Consequently we see that
\bea
d(|\varphi|^2\hat\varphi)&=&2dx^{12}\wedge(2(e^\alpha)'dx^{35}+(\gamma-\delta)'(dx^{36}-dx^{45})-2(e^\beta)'dx^{46}),\nonumber\\
\Lambda d(|\varphi|^2\hat\varphi)&=&2(2(e^\alpha)'dx^{35}+(\gamma-\delta)'(dx^{36}-dx^{45})-2(e^\beta)'dx^{46}),\nonumber\\
d\Lambda d(|\varphi|^2\hat\varphi)&=&4(e^\alpha)''dx^{135}+2(\gamma-\delta)''(dx^{136}-dx^{145})-4(e^\beta)''dx^{146}.\nonumber
\eea
So the Type IIA flow in this case reduces to
\bea
&&
\p_t(e^\alpha)=4(e^\alpha)'',\qquad
\p_t(e^\beta)=4(e^\beta)'',\\
&&
\p_t\gamma=2(\gamma-\delta)'',\qquad
\p_t\delta=-2(\gamma-\delta)''.
\eea
For calculations, it is convenient to introduce $a=2e^\alpha,b=2e^\beta,c=\gamma-\delta,d=\gamma+\delta$ and
\bea
|\varphi|^2=2\sqrt{4e^{\alpha+\beta}-(\gamma-\delta)^2}=2\sqrt{ab-c^2}.\nonumber
\eea
It follows that $d$ is a constant along the flow, while $a$, $b$, $c$ satisfy the standard heat equation:
\bea
\p_t\begin{bmatrix}a & c\\ c & b\end{bmatrix}=4\begin{bmatrix}a & c\\ c & b\end{bmatrix}'',
\eea

Obviously the matrix $\begin{bmatrix}a & c\\ c & b\end{bmatrix}$ converges to a constant matrix as $t$ goes to infinity. Moreover along the flow the positive-definiteness is preserved and the limiting matrix is also positive definite. Thus $\varphi_t$ converges to a positive primitive harmonic form.

Now let us analyze the behavior of $|N|^2$ along the flow. The easiest way to find $|N|^2$ is to use (\ref{evolu}), which says that
\be
|N|^2=e^{-u}\p_t u-(\Delta u+2|du|^2).
\ee
The key is to compute $\Delta u$. Observe that the metric $g$ can be expressed as
\be
g=e^{-u}\begin{bmatrix}ab+d^2-c^2 & -2d & & & &\\ -2d & 4 & & & &\\ & & 2a & -2c & & \\ & & -2c & 2b & & \\ & & & & 2a & 2c \\ & & & & 2c & 2b \end{bmatrix}
\ee
As $\Delta u=g^{ij}(u_{ij}-\Gamma^k_{ij}u_k)=g^{11}u''-g^{ij}\Gamma^1_{ij}u'$, and the Christoffel symbol term can be simplified to
\bea
g^{ij}\Gamma^1_{ij}&=&\frac{1}{2}g^{ij}g^{1l}(\p_ig_{jl}+\p_jg_{il}-\p_lg_{ij})\nonumber\\
&=&g^{1j}g^{1l}g_{jl}'-\frac{1}{2}g^{11}g^{ij}g_{ij}'\nonumber\\
&=&\frac{1}{2}(g^{11})^2g'_{11}+g^{11}g^{12}g'_{12}+((g^{12})^2-\frac{1}{2}g^{11}g^{22})g'_{22}\nonumber\\
&&-g^{11}(g^{33}g'_{33}+2g^{34}g'_{34}+g^{44}g'_{44})\nonumber\\
&=&8e^{-3u}(4u'(ab-c^2)-(a'b+ab'-2cc'))\nonumber\\
&=&8e^{-3u}v',
\eea
where $v=ab-c^2$ and $u=\log 2+\dfrac{1}{2}\log v$. Therefore
\be
\Delta u=4e^{-u}u''-8e^{-3u}u'v'=4e^{-u}(u''-(u')^2).
\ee
Consequently
\bea
|N|^2&=&e^{-u}\p_t u-(\Delta u+2|du|^2)\nonumber\\
&=&2e^{-u}\frac{a''b+ab''-2cc''}{v}-4e^{-u}(u''+(u')^2)\nonumber\\
&=&16e^{-5u}\left(ab\left(2c'-\frac{cb'}{b}-\frac{ca'}{a}\right)^2+\frac{ab-c^2}{ab}(ab'-a'b)^2\right).
\eea
This calculation suggests that $J$ is integrable if and only if $a,b,c$ are proportional to each other. In summary we have proved
\begin{proposition}
Under our ansatze (\ref{torusans}), the Type IIA flow on $(\R/\Z)^6$ reduces to the standard heat equation on $\R/\Z$. If initially $\varphi$ is of the form (\ref{torusans}) whose associated almost complex structure is not integrable, the Type IIA flow still converges to K\"ahler Calabi-Yau geometry.
\end{proposition}

\subsubsection{The Type IIA flow on homogeneous symplectic half-flat manifolds}
\label{ex:half-flat}

Because of Theorem \ref{th:stationary}, the convergence of the Type IIA flow is only possible when the underlying manifold is K\"ahler. We shall see in this subsection and the next that the Type IIA flow can be used to find optimal almost complex structures compatible with a given symplectic form, even when the underlying manifold does not admit any K\"ahler structure.

\medskip

In order to run the Type IIA flow, we first need compact symplectic 6-manifolds with Type IIA structures. A special case of Type IIA structures can be found on the so-called symplectic half-flat manifolds (firstly introduced by de Bartolomeis \cite{dB}, also known as special generalized Calabi-Yau manifolds \cite{dBT}). In our terminology, a symplectic half-flat manifold is simply a symplectic manifold with Type IIA structure $(M,\o,\varphi)$ and the extra condition that $|\varphi|^2$ is constant. Many compact symplectic half-flat manifolds can be constructed as quotients of Lie groups by co-compact lattices, where all the structures are homogeneous under the natural group action. Therefore we shall call symplectic half-flat manifolds constructed in this way \emph{homogeneous}. It is clear that for homogeneous symplectic half-flat manifolds, their geometry up to covering is fully characterized by the underlying Lie algebra, or equivalently the exterior differential system defined by invariant 1-forms. Moreover, homogeneous symplectic half-flat structures have been fully classified by \cite{CT} and \cite{FMOU} when the Lie group is nilpotent or solvable respectively.

\medskip

It is clear that if we run the Type IIA flow on a homogeneous symplectic half-flat manifold with homogeneous initial data, the homogeneity is preserved and the Type IIA flow reduces to a polynomial ODE system. Moreover, in the homogeneous setting, the function $u$ and $|N|^2$ are constants on the manifold, therefore we have the following monotonicity formulas

\begin{proposition} Along the Type IIA flow on homogenous symplectic half-flat manifolds, the following monotonicity formulae hold
\bea
\p_t u &=& e^u|N|^2\geq 0,\label{mf1}\\
\p_t |N|^2 &=& -2e^u|(R^{-J})_{ij}|^2\leq 0.\label{mf2}
\eea
\end{proposition}
{\it Proof}: The first formula follows directly from (\ref{evolu}). For the second formula, we note that Blair-Ianus \cite{BI} proved that
\bea
\p_t\int_M|N|^2\frac{\o^3}{3!}=\int_M(\p_t g_{ij}, (R^{-J})_{ij})\frac{\o^3}{3!}.\label{biv}
\eea
In our case $u$ is a constant, hence (\ref{evolveg2}) becomes $\p_tg_{ij}=-2e^u(R^{-J})_{ij}$, and (\ref{biv}) simplifies to
\bea
\p_t\int_M|N|^2\frac{\o^3}{3!}=-2\int_Me^u|(R^{-J})_{ij}|^2\frac{\o^3}{3!}.\label{hbiv}
\eea
As everything is homogeneous, so all the scalars must be constant, consequently (\ref{hbiv}) still holds without integration, and (\ref{mf2}) is proved. Q.E.D.

\begin{corollary}\label{smallnijenhuis}
Let $(M,\omega)$ be a compact 6-dimensional homogeneous symplectic manifold. If $(M,\o)$ admits a homogeneous symplectic half-flat structure $(M,\o,\varphi_0)$ with which the Type IIA flow exists for all time, then there exist homogenous almost complex structures compatible with $\omega$ and with arbitrary small Nijenhuis tensor.
\end{corollary}
{\it Proof of the Corollary}: We run the Type IIA flow with initial data $\varphi_0$. By monotonicity formulas above, we know that
\bea
\frac{d}{dt}e^{-u}&=&-|N|^2\leq0,\\
\frac{d^2}{dt^2}e^{-u}&=&2e^u|(R^{-J})_{ij}|^2\geq 0.
\eea
So $e^{-u}$ is a monotone non-increasing and convex function with lower bound. If the flow exists for all time, then we must have
$\lim_{t\to\infty}|N|^2\nonumber=-\lim_{t\to\infty}\frac{d}{dt}e^{-u}=0$, as was to be shown. Q.E.D.

\subsubsection*{The Type IIA flow on a nilmanifold}

Now let us consider some explicit examples.

\medskip

Consider the homogeneous symplectic half-flat structure in \cite[Example 5.2]{dBT}, where the Lie algebra of the nilpotent Lie group is characterized by invariant 1-forms $\{e^1,\dots,e^6\}$ satisfying
\bea
&&de^1=de^2=de^3=de^5=0,\nonumber\\
&&de^4=e^{15},\qquad de^6=e^{13}.\nonumber
\eea
Clearly $\o=e^{12}+e^{34}+e^{56}$ defines an invariant symplectic structure. Moreover, this nilpotent Lie group admits co-compact lattices so all the constructions descend to compact nilmanifolds. Consider the ansatze
\be
\label{data:nil}
\varphi=\varphi_{a,b}=(1+a)e^{135}-e^{146}-e^{245}-e^{236}+b(e^{134}-e^{156}),
\ee
it is straightforward to check that $\varphi_{a,b}$ is primitive and closed for any $a,b$. The positivity condition for $\varphi_{a,b}$ is that $\dfrac{1}{16}|\varphi|^4=1+a-b^2>0$. By straightforward calculations, we get
\bea
\hat\varphi=4|\varphi|^{-2}((1+a-b^2)e^1\wedge(e^{36}+e^{45})+e^2\wedge(be^{34}+(1+a)e^{35}-e^{46}-be^{56})).\nonumber
\eea
It follows that
\bea
&&
d(|\varphi|^2\hat\varphi)=4e^{12}(e^{34}+2be^{35}-e^{56}),\nonumber\\
&&
\Lambda d(|\varphi|^2\hat\varphi)=4(e^{34}+2be^{35}-e^{56}),\quad
d\Lambda d(|\varphi|^2\hat\varphi)=8e^{135}.\nonumber
\eea
Therefore under our ansatze the Type IIA flow reduces to the following ODE system
\bea
\dot a(t) = 8,\qquad
\dot b(t) = 0.\nonumber
\eea
Hence the unique solution to the Type IIA flow is
\be
\varphi(t)=(1+a_0+8t)e^{135}-e^{146}-e^{245}-e^{236}+b_0(e^{134}-e^{156}),
\ee
which exists for all time $t\geq0$.

One can easily verify that $\lim_{t\to\infty} J_t$ does not exist and
\bea
|N|^2=(1+a-b^2)^{-3/2}=(1+a_0+8t-b_0^2)^{-3/2}
\eea
is decreasing to zero as $t\to\infty$. This is an explicit example where Corollary \ref{smallnijenhuis} applies.

\subsubsection*{The Type IIA flow on a solvmanifold}

Consider the symplectic half-flat structure on the solvmanifold $M$ constructed by Tomassini and Vezzoni in \cite[Theorem 3.5]{TV}. The geometry of this solvmanifold is characterized by  invariant 1-forms $\{e^j\}_{j=1}^6$ satisfying
\bea
&&
d e^1 = -\lambda e^{15},\quad de^2= \lambda e^{25},\quad
d e^3 = -\lambda e^{36},
\nonumber\\
&&
de^4= \lambda e^{46},\quad
d e^5=0,\quad de^6=0,\nonumber
\eea
where $\lambda=\log \dfrac{3+\sqrt{5}}{2}$. One can easily check that $\o=e^{12}+e^{34}+e^{56}$ is an invariant symplectic form on $M$. A particular symplectic half-flat structure on $M$ takes the form
\bea
\varphi&=&\frac{\sqrt{2}}{2}(e^{135}+e^{136}+e^{145}-e^{146}+e^{235}-e^{236}-e^{245}-e^{246})\nonumber\\
&=&\frac{\sqrt{2}}{2\lambda}d(e^{13}+e^{14}-e^{23}+e^{24}),\nonumber
\eea
so $[\varphi]=0\in H^3(M;\R)$.

Consider the ansatze
\be
\varphi=\alpha(e^{135}+e^{136})+\beta(e^{145}-e^{146})+\gamma(e^{235}-e^{236})-\delta(e^{245}+e^{246}).\label{solvans}
\ee
A direct calculation gives
\be
|\varphi|^2\hat\varphi=8(-\alpha\beta\gamma(e^{135}-e^{136})+\alpha\beta\delta(e^{145}+e^{146})+\alpha\gamma\delta(e^{235}+e^{236}) +\beta\gamma\delta(e^{245}-e^{246})).\nonumber
\ee
The nondegenerate condition is that $|\varphi|^4=64\alpha\beta\gamma\delta>0$. It follows that
\bea
d(|\varphi|^2\hat\varphi)&=&16\lambda(\alpha\beta\gamma e^{1356}+\alpha\beta\delta e^{1456}-\alpha\gamma\delta e^{2356}+\beta\gamma\delta e^{2456}),\nonumber\\
\Lambda d(|\varphi|^2\hat\varphi)&=&16\lambda(\alpha\beta\gamma e^{13}+\alpha\beta\delta e^{14}-\alpha\gamma\delta e^{23}+\beta\gamma\delta e^{24}),\nonumber\\
d\Lambda d(|\varphi|^2\hat\varphi)&=&16\lambda^2(\alpha\beta\gamma(e^{135}+e^{136})+\alpha\beta\delta(e^{145}-e^{146})\nonumber\\
&&+\alpha\gamma\delta(e^{235}-e^{236})- \beta\gamma\delta(e^{245}+e^{246})).\nonumber
\eea
After time rescaling, the Type IIA flow under our ansatze reduces to
\bea
&&
\p_t\alpha=\alpha\beta\gamma,\qquad
\p_t\beta=\alpha\beta\delta,\nonumber\\
&&
\p_t\gamma=\alpha\gamma\delta,\qquad
\p_t\delta=\beta\gamma\delta.\nonumber
\eea
It is easy to see that there exist time-independent nonzero constants $C_1$ and $C_2$ such that $\alpha(t)=C_1\delta(t)$ and $\beta(t)=C_2\gamma(t)$. The ODE system simplifies to
\bea
\p_t\gamma=C_1\gamma\delta^2,\qquad
\p_t\delta=C_2\gamma^2\delta.\nonumber
\eea
Integrate these equations we know there is a constant $C$ such that $C_2\gamma^2-C_1\delta^2=C$, hence
\be
\p_t\gamma=\gamma(C_2\gamma^2-C).\nonumber
\ee
One can solve explicitly
\bea
\gamma^2(t)=\frac{C\gamma_0^2}{C_2\gamma_0^2(1-e^{2Ct})+Ce^{2Ct}},\qquad
\delta^2(t)=\frac{C\delta_0^2e^{2Ct}}{C_1\delta_0^2+C-C_1\delta_0^2e^{2Ct}}.
\eea
Assuming $\varphi$ is initially positive, we know that $C_1,C_2$ are positive constants and $\gamma_0,\delta_0$ are initial values of $\gamma$ and $\delta$ satisfying $\beta\gamma-\alpha\delta=C_2\gamma_0^2-C_1\delta_0^2=C$. When $C=0$, the above formula should be understood as
\bea
\gamma^2=\frac{\gamma_0^2}{1-2C_2\gamma_0^2t},\qquad
\delta^2=\frac{\delta_0^2}{1-2C_1\delta_0^2t}.
\eea
From the above explicit formulas, we can deduce that, no matter what $C$ is, the flow has finite time singularity. A more symmetric expression for the solution (without time rescaling) is
\bea
&&
\alpha(t)=\alpha_0\sqrt{\frac{(\beta_0\gamma_0-\alpha_0\delta_0)e^{32\lambda^2\beta_0\gamma_0t}}{\beta_0\gamma_0e^{32\lambda^2\alpha_0\delta_0t} -\alpha_0\delta_0e^{32\lambda^2\beta_0\gamma_0t}}},\quad
\beta(t)=\beta_0\sqrt{\frac{(\beta_0\gamma_0-\alpha_0\delta_0)e^{32\lambda^2\alpha_0\delta_0t}}{\beta_0\gamma_0e^{32\lambda^2\alpha_0\delta_0t} -\alpha_0\delta_0e^{32\lambda^2\beta_0\gamma_0t}}},\nonumber\\
&&
\gamma(t)=\gamma_0\sqrt{\frac{(\beta_0\gamma_0-\alpha_0\delta_0)e^{32\lambda^2\alpha_0\delta_0t}}{\beta_0\gamma_0e^{32\lambda^2\alpha_0\delta_0t} -\alpha_0\delta_0e^{32\lambda^2\beta_0\gamma_0t}}},\quad
\delta(t)=\delta_0\sqrt{\frac{(\beta_0\gamma_0-\alpha_0\delta_0)e^{32\lambda^2\beta_0\gamma_0t}}{\beta_0\gamma_0e^{32\lambda^2\alpha_0\delta_0t} -\alpha_0\delta_0e^{32\lambda^2\beta_0\gamma_0t}}},\nonumber
\eea
and in the critical case when $\alpha_0\delta_0=\beta_0\gamma_0=S>0$, one has
\bea
&&
\alpha(t)=\frac{\alpha_0}{\sqrt{1-32\lambda^2St}},\qquad
\beta(t)=\frac{\beta_0}{\sqrt{1-32\lambda^2St}},\nonumber\\
&&
\gamma(t)=\frac{\gamma_0}{\sqrt{1-32\lambda^2St}},\qquad
\delta(t)=\frac{\delta_0}{\sqrt{1-32\lambda^2St}}.\nonumber
\eea
From these explicit expressions, we see that the maximal existence time $T$ is given by
\bea
T=\frac{1}{32\lambda^2}\frac{\log(\alpha_0\delta_0)-\log(\beta_0\gamma_0)}{\alpha_0\delta_0-\beta_0\gamma_0}.
\eea
All of $\alpha,\beta,\gamma,\delta$ tend to infinity as $t\to T$, therefore $|\varphi|^2=8(\alpha\beta\gamma\delta)^{1/2}\to\infty$. To compute $|N|^2$, the quickest way is to use (\ref{mf1}):
\bea
|N|^2&=&e^{-u}\p_tu=-\p_te^{-u}=\frac{(\alpha\beta\gamma\delta)^{-3/2}}{16}\p_t(\alpha\beta\gamma\delta)= 2\lambda^2\frac{\alpha\delta+\beta\gamma}{(\alpha\beta\gamma\delta)^{1/2}}\label{Nnorm}\\
&=&\frac{2\lambda^2}{(\alpha_0\beta_0\gamma_0\delta_0)^{1/2}} \left(\alpha_0\delta_0e^{16\lambda^2(\beta_0\gamma_0-\alpha_0\delta_0)t}+\beta_0\gamma_0e^{-16\lambda^2(\beta_0\gamma_0-\alpha_0\delta_0)t}\right) \nonumber\\ &\geq&4\lambda^2.\nonumber
\eea

Now let us analyze the behavior of $\varphi$ in detail. We shall see that the Type IIA flow naturally leads us to optimal almost-complex structures compatible with $\o$.
\begin{enumerate}
\item No matter what the integral constant $C=\beta_0\gamma_0-\alpha_0\delta_0$ is, the flow of $\varphi$ blows up when $t\to T$, and the same is true for the metric $\tilde g$. However the expressions of $g$, $J$, and $N$ extend smoothly to $t=T$. In fact, the limit $|\varphi|^{-1}\varphi$ as $t\to T$ exists.

\item In the critical case $C=\beta_0\gamma_0-\alpha_0\delta_0=0$, the flow of $\varphi$ (as well as $\tilde g$) is a self-expander in the sense that
\be
\varphi(t)=\frac{\varphi_0}{\sqrt{1-32\lambda^2St}}
\ee
for a positive constant $S$ determined by initial data. In this case, all of $g$, $J$ and $N$ are stationary with $|N|^2=4\lambda^2$. In fact such $J$ provide examples of harmonic almost-complex structures in the sense of Blair-Ianus \cite{BI}, namely almost-complex structures compatible with $\o$ and satisfying
\be
(R^{-J})_{ij}=\frac{1}{2}(R_{ij}-R_{Ji,Jj})=0.
\ee
These harmonic almost-complex structures are critical points of the energy functional studied by Blair-Ianus \cite{BI} and L\^e-Wang \cite{LW}.

\item When the integral constant $C=\beta_0\gamma_0-\alpha_0\delta_0$ is not zero, all of $g$, $J$ and $N$ are evolving. When $t$ approaches $T$, the limit $\lim_{t\to T}J(t)$ exists and is a harmonic almost-complex structure, which is also a minimizer of $|N|^2$ among all almost complex structures associated to our ansatze (\ref{solvans}).
\end{enumerate}

\bigskip
\noindent
{\bf Acknowledgements} The authors would like to thank F. Fong, R. Rosso, A. Sun, and L.-S. Tseng for very helpful communications. The second-named author would also like to thank the Galileo Galilei Institute for Theoretical Physics, the Banff Research Station, and the American Institute of Mathematics for their hospitality.

\bigskip

\noindent Department of Mathematics $\&$ Computer Science, Rutgers, Newark, NJ 07102, USA

\smallskip

\noindent teng.fei@rutgers.edu

\bigskip

\noindent Department of Mathematics, Columbia University, New York, NY 10027, USA

\smallskip

\noindent phong@math.columbia.edu

\bigskip

\noindent \noindent Mathematics Department, University of British Columbia, Vancouver, BC V6T 1Z2, CAN

\smallskip

\noindent spicard@math.ubc.ca

\bigskip

\noindent Department of Mathematics, University of California, Irvine, CA 92697, USA

\smallskip
\noindent xiangwen@math.uci.edu

\end{document}